\documentclass{article}

\usepackage{overpic}
\usepackage{amsmath}
\usepackage{xfrac}
\usepackage{blindtext}
\usepackage{comment}
\usepackage{faktor}
\usepackage[all]{xy}
\usepackage[french, main=english]{babel}
\usepackage[utf8x]{inputenc}
\usepackage{float}
\usepackage{graphicx}
\usepackage[colorinlistoftodos]{todonotes}
\usepackage{url}
\usepackage{hyperref}
\usepackage{array}
\usepackage{tabularx}
\usepackage{setspace}
\usepackage[T1]{fontenc}
\usepackage{subfig}
\usepackage[margin=1in]{geometry} 
\usepackage{fancyhdr}
\usepackage{rotating}
\usepackage{amsfonts,amssymb,amsthm,scrextend}
\usepackage{pinlabel}
\usepackage[nottoc,numbib]{tocbibind}
\usepackage{appendix}
\usepackage{quiver}
\usepackage{lscape}
\usepackage{mathtools}
\usepackage{tikz}

\newcommand{\Z}{\mathbb{Z}}

\newcommand{\R}{\mathbb{R}}

\newcommand{\X}{\mathbb{X}}

\newcommand{\E}{\mathcal{E}}
\newcommand{\Ai}{\mathcal{A}_{\infty}}

\newcommand{\Tree}{\overline{\mathcal{T}}^{\X}}

\newcommand{\deffct}[5]{#1\colon \begin{array}{ccc}
    #2 & \to & #3  \\
     #4 & \mapsto & \displaystyle #5
\end{array}}
\newcommand{\trajb}[1]{\overline{\mathcal{L}}(#1)}
\newcommand{\trajbi}[2]{\overline{\mathcal{L}}_{#1}(#2)}

\newcommand{\wb}[2]{\overline{W}^{#1}(#2)}
\newcommand{\wbi}[3]{\overline{W}^{#1}_{#2}(#3)}

\newcommand{\fibration}{E \overset{\pi}{\to} X}

\newcommand{\Crit}{\textup{Crit}}
\newcommand{\Or}{\textup{Or}}
\newcommand{\Id}{\textup{Id}}

\newcommand{\Hom}{\textup{Hom}}
\newcommand{\End}{\textup{End}}
\newcommand{\EZ}{\textup{EZ}}
\newcommand{\ev}{\textup{ev}}

\newcommand{\blue}[1]{\textcolor{blue}{#1}}

\newcommand{\Sk}{\textup{Sk}}

\newcommand{\Y}{\mathcal{Y}}

\newcommand{\Xq}{X/\Y}

\newcommand{\ftimes}[2]{\ _{#1}\!\times_{#2}}

\newcommand{\ls}[1]{\mathcal{L} #1}
\newcommand{\CS}{\textup{CS}_{DG}}

\newcommand{\m}{\mathbf{\widetilde{m}}}
\newcommand{\mL}{\mathbf{\widetilde{m}^L}}

\renewcommand{\qed}{\hfill$\blacksquare$}

\renewenvironment{proof}{\begin{addmargin}[1em]{0em}\begin{newproof}}{\end{newproof}\end{addmargin}\qed}

\definecolor{darkgreen}{RGB}{0,168,0}

\newtheorem{thm}{Theorem}[section]
\newtheorem{prop}[thm]{Proposition}
\newtheorem{cor}[thm]{Corollary}
\newtheorem{lemme}[thm]{Lemma}

\newtheorem{theorem}{Theorem}

\newtheorem{corollary}[theorem]{Corollary}

\theoremstyle{definition}
\newtheorem{notation}[thm]{Notation}
\newtheorem{defi}[thm]{Definition}
\newtheorem{rem}[thm]{Remark}
\newtheorem{ex}[thm]{Example}

\newcommand{\drawingperturbedgradienttree}{
\tikzset{every picture/.style={line width=0.75pt}} 

\begin{tikzpicture}[x=0.75pt,y=0.75pt,yscale=-1,xscale=1]

\draw    (331.72,141.73) -- (332.44,252) ;
\draw [shift={(332.11,201.87)}, rotate = 269.63] [fill={rgb, 255:red, 74; green, 78; blue, 226 }  ][line width=0.08]  [draw opacity=0] (8.93,-4.29) -- (0,0) -- (8.93,4.29) -- cycle    ;
\draw [color={rgb, 255:red, 0; green, 0; blue, 0 }  ,draw opacity=1 ]   (215.38,25.48) -- (256.1,65) -- (331.72,141.73) ;
\draw [shift={(239.33,48.72)}, rotate = 224.15] [fill={rgb, 255:red, 74; green, 78; blue, 226 }  ,fill opacity=1 ][line width=0.08]  [draw opacity=0] (8.93,-4.29) -- (0,0) -- (8.93,4.29) -- cycle    ;
\draw [shift={(297.42,106.93)}, rotate = 225.42] [fill={rgb, 255:red, 74; green, 78; blue, 226 }  ,fill opacity=1 ][line width=0.08]  [draw opacity=0] (8.93,-4.29) -- (0,0) -- (8.93,4.29) -- cycle    ;
\draw    (453.38,25.15) -- (415.38,64.15) -- (331.72,141.73) ;
\draw [shift={(430.89,48.23)}, rotate = 314.26] [fill={rgb, 255:red, 74; green, 78; blue, 226 }  ][line width=0.08]  [draw opacity=0] (8.93,-4.29) -- (0,0) -- (8.93,4.29) -- cycle    ;
\draw [shift={(369.89,106.34)}, rotate = 317.16] [fill={rgb, 255:red, 74; green, 78; blue, 226 }  ][line width=0.08]  [draw opacity=0] (8.93,-4.29) -- (0,0) -- (8.93,4.29) -- cycle    ;
\draw    (295.38,26.15) -- (256.1,65) ;
\draw [shift={(272.19,49.09)}, rotate = 315.32] [fill={rgb, 255:red, 74; green, 78; blue, 226 }  ][line width=0.08]  [draw opacity=0] (8.93,-4.29) -- (0,0) -- (8.93,4.29) -- cycle    ;
\draw    (374.38,27.15) -- (415.38,64.15) ;
\draw [shift={(398.6,49)}, rotate = 222.06] [fill={rgb, 255:red, 74; green, 78; blue, 226 }  ][line width=0.08]  [draw opacity=0] (8.93,-4.29) -- (0,0) -- (8.93,4.29) -- cycle    ;

\draw (302.98,84.36) node [anchor=north west][inner sep=0.75pt]    {$e$};
\draw (354.73,82.36) node [anchor=north west][inner sep=0.75pt]    {$f$};
\draw (237,24.4) node [anchor=north west][inner sep=0.75pt]    {$e_{1}$};
\draw (416,24.4) node [anchor=north west][inner sep=0.75pt]    {$e_{4}$};
\draw (394,24.4) node [anchor=north west][inner sep=0.75pt]    {$e_{3}$};
\draw (260,24.4) node [anchor=north west][inner sep=0.75pt]    {$e_{2}$};
\draw (310,192.4) node [anchor=north west][inner sep=0.75pt]    {$e_{0}$};
\draw (203,4.4) node [anchor=north west][inner sep=0.75pt]    {$x_{1}$};
\draw (287,4.4) node [anchor=north west][inner sep=0.75pt]    {$x_{2}$};
\draw (365,4.4) node [anchor=north west][inner sep=0.75pt]    {$x_{3}$};
\draw (445,4.4) node [anchor=north west][inner sep=0.75pt]    {$x_{4}$};
\draw (328,257.4) node [anchor=north west][inner sep=0.75pt]    {$y$};
\draw (211,45.4) node [anchor=north west][inner sep=0.75pt]  [color={rgb, 255:red, 74; green, 78; blue, 226 }  ,opacity=1 ]  {$\gamma _{e_{1}}$};
\draw (277.74,48.98) node [anchor=north west][inner sep=0.75pt]  [color={rgb, 255:red, 74; green, 78; blue, 226 }  ,opacity=1 ]  {$\gamma _{e_{2}}$};
\draw (370,48.4) node [anchor=north west][inner sep=0.75pt]  [color={rgb, 255:red, 74; green, 78; blue, 226 }  ,opacity=1 ]  {$\gamma _{e_{3}}$};
\draw (436.38,48.05) node [anchor=north west][inner sep=0.75pt]  [color={rgb, 255:red, 74; green, 78; blue, 226 }  ,opacity=1 ]  {$\gamma _{e_{4}}$};
\draw (268,104.4) node [anchor=north west][inner sep=0.75pt]  [color={rgb, 255:red, 74; green, 78; blue, 226 }  ,opacity=1 ]  {$\gamma _{e}$};
\draw (375.55,106.34) node [anchor=north west][inner sep=0.75pt]  [color={rgb, 255:red, 74; green, 78; blue, 226 }  ,opacity=1 ]  {$\gamma _{f}$};
\draw (339,192.4) node [anchor=north west][inner sep=0.75pt]  [color={rgb, 255:red, 74; green, 78; blue, 226 }  ,opacity=1 ]  {$\gamma _{e_{0}}$};

\end{tikzpicture}
}

\newcommand{\drawingbrokentree}{

\tikzset{every picture/.style={line width=0.75pt}} 

\begin{tikzpicture}[x=0.5pt,y=0.5pt,yscale=-1,xscale=1]
\draw    (584.3,172.33) -- (584.78,240) ;
\draw    (507.75,100.98) -- (584.3,172.33) ;
\draw    (662.38,99.55) -- (584.3,172.33) ;
\draw    (552.92,101.26) -- (534.55,125.24) ;
\draw    (622.58,99.55) -- (638.65,122.38) ;
\draw    (71.68,172.87) -- (72.12,241) ;
\draw    (0.38,101.05) -- (71.68,172.87) ;
\draw    (36.75,100.33) -- (21.77,122.6) ;
\draw    (144.41,99.61) -- (71.68,172.87) ;
\draw    (129.44,62.55) -- (144.41,85.53) ;
\draw    (159.38,63.27) -- (144.41,85.53) ;
\draw    (144.41,85.53) -- (144.41,99.61) ;
\draw    (257.67,171.81) -- (257.67,239) ;
\draw    (186.38,100.97) -- (257.67,171.81) ;
\draw    (330.38,99.55) -- (257.67,171.81) ;
\draw    (237.38,100.15) -- (216.38,131.15) ;
\draw    (388.44,144.55) -- (403.41,167.53) ;
\draw    (418.38,145.27) -- (403.41,167.53) ;
\draw    (403.41,167.53) -- (403.41,181.61) ;

\draw (475.67,162.48) node [anchor=north west][inner sep=0.75pt]    {$\displaystyle \lim _{l_{2}\rightarrow +\ \infty }$};
\draw (548.04,150.78) node [anchor=north west][inner sep=0.75pt]    {$l_{1}$};
\draw (607.75,152.2) node [anchor=north west][inner sep=0.75pt]    {$l_{2}$};
\draw (37.46,151.24) node [anchor=north west][inner sep=0.75pt]    {$l_{1}$};
\draw (151.01,165.81) node [anchor=north west][inner sep=0.75pt]    {$=$};
\draw (223.45,150.35) node [anchor=north west][inner sep=0.75pt]    {$l_{1}$};
\draw (449.01,165.81) node [anchor=north west][inner sep=0.75pt]    {$=$};
\draw (333,163.4) node [anchor=north west][inner sep=0.75pt]    {$\circ _{3}$};

\end{tikzpicture}
}

\newcommand{\flowtree}{

\tikzset{every picture/.style={line width=0.75pt}} 

\begin{tikzpicture}[x=0.75pt,y=0.75pt,yscale=-1,xscale=1]

\draw    (324.02,179.81) -- (324.86,300) ;
\draw    (187.38,19.36) -- (324.02,179.81) ;
\draw    (462.38,17.15) -- (324.02,179.81) ;
\draw    (310.38,20.15) -- (255.7,99.59) ;
\draw    (341.38,20.15) -- (393.7,97.98) ;
\draw  [color={rgb, 255:red, 74; green, 78; blue, 226 }  ,draw opacity=1 ][fill={rgb, 255:red, 74; green, 78; blue, 226 }  ,fill opacity=1 ] (216.31,57) .. controls (216.31,55.51) and (217.51,54.31) .. (219,54.31) .. controls (220.49,54.31) and (221.69,55.51) .. (221.69,57) .. controls (221.69,58.49) and (220.49,59.69) .. (219,59.69) .. controls (217.51,59.69) and (216.31,58.49) .. (216.31,57) -- cycle ;
\draw  [color={rgb, 255:red, 208; green, 2; blue, 27 }  ,draw opacity=1 ][fill={rgb, 255:red, 208; green, 2; blue, 27 }  ,fill opacity=1 ] (282.35,56.18) .. controls (282.35,54.69) and (283.56,53.48) .. (285.04,53.48) .. controls (286.53,53.48) and (287.73,54.69) .. (287.73,56.18) .. controls (287.73,57.66) and (286.53,58.87) .. (285.04,58.87) .. controls (283.56,58.87) and (282.35,57.66) .. (282.35,56.18) -- cycle ;
\draw  [color={rgb, 255:red, 65; green, 117; blue, 5 }  ,draw opacity=1 ][fill={rgb, 255:red, 65; green, 117; blue, 5 }  ,fill opacity=1 ] (363.31,56) .. controls (363.31,54.51) and (364.51,53.31) .. (366,53.31) .. controls (367.49,53.31) and (368.69,54.51) .. (368.69,56) .. controls (368.69,57.49) and (367.49,58.69) .. (366,58.69) .. controls (364.51,58.69) and (363.31,57.49) .. (363.31,56) -- cycle ;
\draw  [color={rgb, 255:red, 245; green, 166; blue, 35 }  ,draw opacity=1 ][fill={rgb, 255:red, 245; green, 166; blue, 35 }  ,fill opacity=1 ] (426.31,55) .. controls (426.31,53.51) and (427.51,52.31) .. (429,52.31) .. controls (430.49,52.31) and (431.69,53.51) .. (431.69,55) .. controls (431.69,56.49) and (430.49,57.69) .. (429,57.69) .. controls (427.51,57.69) and (426.31,56.49) .. (426.31,55) -- cycle ;
\draw  [color={rgb, 255:red, 74; green, 78; blue, 226 }  ,draw opacity=1 ][fill={rgb, 255:red, 74; green, 78; blue, 226 }  ,fill opacity=1 ] (253.01,96.89) .. controls (253.01,95.41) and (254.21,94.2) .. (255.7,94.2) .. controls (257.19,94.2) and (258.39,95.41) .. (258.39,96.89) .. controls (258.39,98.38) and (257.19,99.59) .. (255.7,99.59) .. controls (254.21,99.59) and (253.01,98.38) .. (253.01,96.89) -- cycle ;
\draw  [color={rgb, 255:red, 74; green, 78; blue, 226 }  ,draw opacity=1 ][fill={rgb, 255:red, 74; green, 78; blue, 226 }  ,fill opacity=1 ] (321.33,177.12) .. controls (321.33,175.63) and (322.53,174.43) .. (324.02,174.43) .. controls (325.5,174.43) and (326.71,175.63) .. (326.71,177.12) .. controls (326.71,178.61) and (325.5,179.81) .. (324.02,179.81) .. controls (322.53,179.81) and (321.33,178.61) .. (321.33,177.12) -- cycle ;
\draw  [color={rgb, 255:red, 245; green, 166; blue, 35 }  ,draw opacity=1 ][fill={rgb, 255:red, 245; green, 166; blue, 35 }  ,fill opacity=1 ] (391.01,95.29) .. controls (391.01,93.8) and (392.21,92.6) .. (393.7,92.6) .. controls (395.19,92.6) and (396.39,93.8) .. (396.39,95.29) .. controls (396.39,96.78) and (395.19,97.98) .. (393.7,97.98) .. controls (392.21,97.98) and (391.01,96.78) .. (391.01,95.29) -- cycle ;
\draw  [color={rgb, 255:red, 0; green, 0; blue, 0 }  ,draw opacity=1 ][fill={rgb, 255:red, 0; green, 0; blue, 0 }  ,fill opacity=1 ] (321.31,228) .. controls (321.31,226.51) and (322.51,225.31) .. (324,225.31) .. controls (325.49,225.31) and (326.69,226.51) .. (326.69,228) .. controls (326.69,229.49) and (325.49,230.69) .. (324,230.69) .. controls (322.51,230.69) and (321.31,229.49) .. (321.31,228) -- cycle ;

\draw (201,55) node [anchor=north west][inner sep=0.75pt]  [color={rgb, 255:red, 74; green, 78; blue, 226 }  ,opacity=1 ]  {$a_{1}$};
\draw (192,98.4) node [anchor=north west][inner sep=0.75pt]  [font=\footnotesize,color={rgb, 255:red, 74; green, 78; blue, 226 }  ,opacity=1 ]  {$\textcolor[rgb]{0.29,0.31,0.89}{\phi _{e_{1} ,X}( a_{1})}\textcolor[rgb]{0.82,0.01,0.11}{\ } =\textcolor[rgb]{0.82,0.01,0.11}{\ \phi _{e_{2} ,X}( a_{2})}$};
\draw (213,23.4) node [anchor=north west][inner sep=0.75pt]  [color={rgb, 255:red, 0; green, 0; blue, 0 }  ,opacity=1 ]  {$e_{1}$};
\draw (274,24.4) node [anchor=north west][inner sep=0.75pt]    {$e_{2}$};
\draw (355,23.4) node [anchor=north west][inner sep=0.75pt]    {$e_{3}$};
\draw (456,30.4) node [anchor=north west][inner sep=0.75pt]    {$e_{4}$};
\draw (270,135.4) node [anchor=north west][inner sep=0.75pt]    {$e$};
\draw (367,139.4) node [anchor=north west][inner sep=0.75pt]    {$f$};
\draw (297,261.4) node [anchor=north west][inner sep=0.75pt]    {$e_{0}$};
\draw (288.04,55) node [anchor=north west][inner sep=0.75pt]  [color={rgb, 255:red, 208; green, 2; blue, 27 }  ,opacity=1 ]  {$a_{2}$};
\draw (348,55) node [anchor=north west][inner sep=0.75pt]  [color={rgb, 255:red, 65; green, 117; blue, 5 }  ,opacity=1 ]  {$a_{3}$};
\draw (432,55) node [anchor=north west][inner sep=0.75pt]  [color={rgb, 255:red, 245; green, 166; blue, 35 }  ,opacity=1 ]  {$a_{4}$};
\draw (336,99.4) node [anchor=north west][inner sep=0.75pt]  [font=\footnotesize,color={rgb, 255:red, 65; green, 117; blue, 5 }  ,opacity=1 ]  {$\textcolor[rgb]{0.25,0.46,0.02}{\phi _{e_{3} ,X}( a_{3})}\textcolor[rgb]{0.82,0.01,0.11}{\ } =\textcolor[rgb]{0.97,0.91,0.11}{\ }\textcolor[rgb]{0.96,0.65,0.14}{\phi _{e_{4} ,X}( a_{4})}$};
\draw (245,176.4) node [anchor=north west][inner sep=0.75pt]  [font=\footnotesize,color={rgb, 255:red, 74; green, 78; blue, 226 }  ,opacity=1 ]  {$\textcolor[rgb]{0.29,0.31,0.89}{\phi _{e}( \phi }\textcolor[rgb]{0.29,0.31,0.89}{_{e_{1} ,\X}}\textcolor[rgb]{0.29,0.31,0.89}{(}\textcolor[rgb]{0.29,0.31,0.89}{a}\textcolor[rgb]{0.29,0.31,0.89}{_{1}}\textcolor[rgb]{0.29,0.31,0.89}{))}\textcolor[rgb]{0.82,0.01,0.11}{\ } =\textcolor[rgb]{0.82,0.01,0.11}{\ }\textcolor[rgb]{0.96,0.65,0.14}{\phi _{f}( \phi _{e_{4} ,\X}( a_{4}))}$};
\draw (326,228.71) node [anchor=north west][inner sep=0.75pt]    {$ \begin{array}{l}
\Phi_{T,\X}(\textcolor[rgb]{0.29,0.31,0.89}{a_{1}} ,\textcolor[rgb]{0.82,0.01,0.11}{a_{2}} ,\textcolor[rgb]{0.25,0.46,0.02}{a_{3}} ,\textcolor[rgb]{0.96,0.65,0.14}{a_{4}})\\
\in \Delta _{4}\overline{W}^{s}(y)
\end{array}$};
\draw (178,1.4) node [anchor=north west][inner sep=0.75pt]    {$x_{1}$};
\draw (300,1.4) node [anchor=north west][inner sep=0.75pt]    {$x_{2}$};
\draw (334,1.4) node [anchor=north west][inner sep=0.75pt]    {$x_{3}$};
\draw (455,1.4) node [anchor=north west][inner sep=0.75pt]    {$x_{4}$};
\draw (321,305.4) node [anchor=north west][inner sep=0.75pt]    {$y$};

\end{tikzpicture}
}
\usepackage{bbm}
\usepackage{makecell}
\usepackage{csquotes}
\usepackage{enumerate}
\usepackage{ulem}
\usepackage{cancel}
\newcommand{\ux}{\underline{x}}
\newcommand{\uz}{\underline{z}}
\newcommand{\uw}{\underline{w}}
\newcommand{\id}{\operatorname{id}}
\newcommand{\dec}{\operatorname{dec}}

\newcommand{\PG}{\mathcal{PG}}
\renewcommand{\wb}{\overline{W}}

\title{An \texorpdfstring{$\Ai$}{Ainfty}-algebra structure for enriched Morse complexes}
\author{Jonathan Clivio \footnote{email: jclivio@math.unistra.fr}\\ IRMA, University of Strasbourg \vspace{1em} \\  Robin Riegel \footnote{email: r.riegel@math.unistra.fr} \\ IRMA, University of Strasbourg}
\date{}

\begin{document}

\setlength{\parskip}{5pt}
\setlength{\parindent}{0pt}

\maketitle

\begin{abstract}
    We construct a homotopy coherent product structure extending the Chas-Sullivan product on chains of the free loop space of a manifold. Our construction uses the framework of Morse homology with DG coefficients and perturbed Morse gradient trees. We define an operad given by moduli spaces of Morse gradient trees with a twisted differential. This Morse tree operad receives a map from a suspended $\Ai$-operad and in turn acts on Morse homology with DG coefficients associated to a fibration over a manifold with a fibrewise monoidal structure.
\end{abstract}

\section{Introduction} 

Given a fibration $F \hookrightarrow E \overset{\pi}{\to} X$ over a manifold with a fibrewise monoidal structure, we produce a homotopy coherent algebra structure -- an $\Ai$-{algebra} structure -- on chains of the total space. In the case of the loop-loop fibration $\Omega X \hookrightarrow \ls X \overset{\ev_0}{\to} X$, this extends the Chas-Sullivan product, also known as the string topology product \cite{ChaSu99} to an $\Ai$-{algebra} structure on $C_*(\mathcal{L} X)$.

We do this in the framework of Morse homology with DG coefficients as introduced by Barraud-Damian-Humilière-Oancea in \cite{BDHO23} that extends Morse homology on a manifold $X$ to allow differential graded coefficients. In the context of a fibration over $X$, there exists a chain complex $C_*(X,C_*(F))$ that computes the singular homology of the total space $E$ given a Morse function $f\colon X\to \R$. 

Its underlying graded groups are given by $\bigoplus_{x\in \Crit(f)} C_{-|x|+*}(F_x)$ where $|x|$ denotes the index of a critical point $x$ and $F_x$ the fibre over that point. The differential is computed using a family of parametrisations of the spaces of Morse trajectories and a transitive lifting function, 
\begin{align*}
    \Phi\colon E \ {}_\pi \!\times_{\ev_0} \mathcal{P}X\to E,
\end{align*}
giving a preferred way to push from one fibre into another fibre.

Given an associative {fibrewise} monoidal structure  $\mu \colon E \ftimes{\pi}{\pi} E\to E $, we produce an $\Ai$-algebra structure on $C_*(X,C_*(\widetilde{F}))$ where $\widetilde{F}_x$ is a homotopy equivalent replacement of $F_x$. This exhibits the generalised Chas-Sullivan product as in \cite{GS07} and \cite{Rie24} as the binary operation of an $\Ai$-algebra structure. This result was proved for $\Z/2\Z$ coefficients in \cite[Theorem 6.1.14]{phdthesisrie}.

\begin{theorem}[Proposition \ref{prop : binary operation} and Theorem \ref{thm:main thm}]\label{theorem : Ai structure}
    Let $\fibration$ be a fibration with an associative fibrewise monoidal structure. A transitive lifting function induces a shifted $\Ai$-{algebra} structure on Morse homology with DG coefficients $C_*(X,C_*(\widetilde{F}))\simeq C_*(E)$. 
    
    Moreover, the binary operation is a chain-level model of the generalised Chas-Sullivan product.
\end{theorem}
As an application of this theorem, we produce an $\Ai$-algebra structure extending the original Chas-Sullivan product or string topology product \cite{ChaSu99}.

\begin{corollary}[Example \ref{exmp:loop-loop}]
    There exists an $\Ai$-algebra structure on $C_{n+*}(X ,C_*(\Omega X))$ that induces the Chas-Sullivan product $H_{n+*}(\ls X)^{\otimes 2} \to H_{n+*}(\ls X).$ 
\end{corollary}

This result is of interest beyond the framework of Morse homology with DG coefficients: it is an open question whether this operation is homotopy invariant. At the same time, our use of the framework of Morse homology with DG coefficients gives an explicit geometric description of this $\Ai$-algebra structure. In the following, we sketch the geometric origins of this construction.

\subsection{Perturbed Morse gradient trees}

The operations of the $\Ai$-{algebra} structure in Theorem \ref{theorem : Ai structure} are induced by perturbed Morse gradient trees as defined by Abouzaid in \cite{Abo11}. Abouzaid used Morse gradient trees to produce an $\Ai$-algebra structure on Morse homology for the intersection product. In particular, we extend this construction to allow for Morse homology with DG coefficients coming from a fibration over $X$ which carries a fibrewise monoidal structure $\mu\colon E\ {}_\pi\!\times_\pi E\to E$.

\begin{figure}[H]
    \centering
    \includegraphics[width=0.2\linewidth]{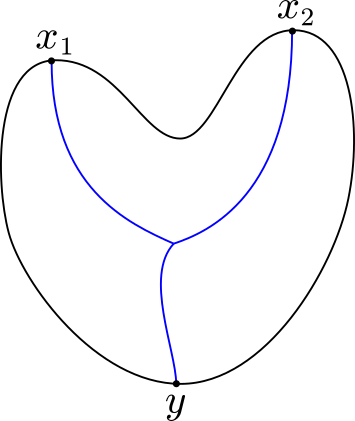}
    \caption{A perturbed Morse gradient tree in the sphere with two input legs.}
    \label{fig:perturb Morse tree}
\end{figure}

A perturbed Morse gradient tree is a metric tree embedded in the manifold $X$. It has $k$ input legs that start at critical points $x_1,\dots, x_k$ and one output leg that ends at a critical point $y$ (see Figure \ref{fig:perturb Morse tree}). The edges of the trees follow the trajectory of the gradient of a Morse function up to a perturbation $\X$ supported near the vertices of the tree. 

We denote the moduli space of such perturbed Morse gradient trees by $\Tree_k(\ux;y)\coloneqq\Tree_k(x_1,\dots,x_k;y)$. Abouzaid showed in \cite{Abo11} that $\Tree_k(\ux,y)$ is a compact, orientable manifold with boundary (see Section \ref{sec : perturbed Morse}).

A tree in $\Tree_k(\ux;y)$ defines $k$ trajectories from the critical points $x_i$ to $y$. Choosing appropriate parametrisations of these trajectories we obtain maps
\begin{align*}
    \Gamma\colon \Tree_k(\ux;y)\to \mathcal{P}_{x_1\to y}X \times \dots \times \mathcal{P}_{x_k\to y} X
\end{align*}
where $\mathcal{P}_{x_i\to y}X$ denotes the space of piecewise smooth paths from $x_i$ to $y$ (see Lemma \ref{lemme : evaluation maps trees}). Pushing along these paths using the transitive lifting function $\Phi$ and then applying the monoidal structure $\mu$, we get a map
\begin{align}\label{eq:topo operad action}
     \Tree_k(\ux;y)\times F_{x_1}\times \dots \times F_{x_k}\xrightarrow{\Gamma}\mathcal{P}_{x_1\to y}X \times \dots \times \mathcal{P}_{x_k\to y} X\times F_{x_1}\times \dots \times F_{x_k} \xrightarrow{\Phi} F_y^k\xrightarrow{\mu}F_y.
\end{align}

This map is the key geometric input. After possibly replacing $F_x$ with homotopy equivalent spaces $\widetilde{F}_x$, we get maps
\begin{align*}
    \mu_\Phi\colon \Tree_k(\ux;y)\to \mathrm{Map}(\widetilde{F}_{x_1}\times\dots \times \widetilde{F}_{x_k},\widetilde{F}_y).
\end{align*}
that induce maps on Morse homology with DG coefficients. Explicitly, the $k$-ary operation of the $\Ai$-{algebra} structure in Theorem \ref{theorem : Ai structure} is given on 
\begin{align*}
    \alpha_1\otimes \dots \otimes\alpha_k \in C_{-|x_1|+*}(\widetilde{F}_{x_1})\otimes \dots \otimes C_{-|x_k|+*}(\widetilde{F}_{x_k})\subseteq C_*(X,C_*(\widetilde{F}))^{\otimes k}
\end{align*}
by
\begin{align}\label{eq:alg operad action}
    \alpha_1\otimes \dots \otimes\alpha_k \mapsto\sum_{y} (  \mu_\Phi)_*(t_{\ux,y})( \alpha_1\times \dots \times \alpha_k)
\end{align}
for $t_{\ux,y}\in C_*(\Tree_k(\ux;y))$ a generator of the top homology relative to the boundary.

The binary operation has a particular description: given chains in $C_{-|x_1|+*}(F_{x_1})$ and $C_{-|x_2|+*}(F_{x_2})$, their product is calculated by pushing them along the trees in the moduli space $\Tree_2(x_1,x_2;y)$ for all $y\in \Crit(f)$. These moduli spaces $\Tree_2(x_1,x_2;y)$ can be described in the following way.

A tree in $\Tree_2(x_1,x_2;y)$ is given by two perturbed trajectories going down from $x_1$ and $x_2$ and one perturbed trajectory going up from $y$. In fact, we need a perturbation $\X$ such that these trajectories intersect transversely  (see Figure \ref{fig:perturb Morse tree}). This gives a description of the product as \enquote{intersecting chains on the base}. The second author showed in \cite{Rie24} that this describes the generalised Chas-Sullivan product as introduced in \cite{GS07}.

\subsection{Operadic perspective}

It remains to show that the maps \eqref{eq:alg operad action} assemble to an $\Ai$-algebra structure for an appropriate choice of generators $t_{\ux,y}\in C_*(\Tree_k(\ux;y))$. To show this, it is helpful to phrase the above construction in terms of operad actions. 

We first note that the moduli spaces $\Tree_k(\ux;y)$ together with concatenation of trees define a $\Crit(f)$-coloured operad in topological spaces. The maps \eqref{eq:topo operad action} induce a homotopy coherent action of this operad on $(F_x)_{x\in \Crit(f)}$. Using a standard result from the operadic literature (\cite{boardman2006homotopy,berger2006boardman,berger2005resolution}), we can strictify this action: after possibly replacing $F_x$ by some homotopy equivalent spaces $\widetilde{F}_x$, we get a strict operad action via the maps $\mu_\Phi$.

We recall that Morse homology is given by the complex $\bigoplus_{x\in \Crit(f)} C_{-|x|+*}(F_x)$ with a twisted differential. The operations of taking chains, shifting degrees, taking the direct sum and twisting the differential have a counterpart on the level of operads.  Applying these operations to the coloured topological operad of moduli spaces of gradient trees, we obtain a (one-coloured algebraic) operad $\mathcal{MT}^{(X,f,\X)}$ (see Definition \ref{def:Morse tree op}). This operad has an action on Morse homology with DG coefficients.

\begin{theorem}[Corollary \ref{cor:dgM acting on Morse homology}]\label{theorem:MT acting on Morse homology}
    Let $\pi\colon E\to X$ be a fibration with an associative fibrewise monoidal structure. The map \eqref{eq:topo operad action} defines an $\mathcal{MT}^{(X,f,\X)}$-action on $C_*(X,C_*(\widetilde{F}))\simeq C_*(E)$.
\end{theorem}

To prove Theorem \ref{theorem : Ai structure}, we thus need to produce an interesting map of operads $\underline{\Sigma}^n\otimes \Ai\to \mathcal{MT}^{(X,f,\X)}$ where $\underline{\Sigma}^n\otimes \Ai$ denotes the $n$-suspended $\Ai$-operad. 

This map of operads only sees the top homology of the operad $\mathcal{MT}^{(X,f,\X)}$ and in turn only the top homology of the moduli space $\Tree_k(\ux;y)$.

To produce the map of operads $\underline{\Sigma}^n\otimes \Ai\to \mathcal{MT}^{(X,f,\X)}$, we produce a $\Crit(f)$-coloured version of the $\Ai$-operad. Algebras over this operad are given by $\Ai$-algebras whose underlying graded abelian group has a $\Crit(f)$-decomposition. This fits our application because $C_*(X,C_*(\widetilde{F}))$ has a canonical decomposition of the underlying graded abelian group by $C_*(X,C_*(\widetilde{F}))\cong \bigoplus_{x\in \Crit(f)} C_{-|x|+*}(\widetilde{F}_x)$. 

Picking generators $t_{\ux,y}\in C_*(\Tree_k(\ux;y))$ of the top homology thus is equivalent to giving a map of operads from a suspended version of the $\Crit(f)$-coloured $\Ai$-operad. This description has the advantage that it gives explicit signed relations that the orientations of the spaces $\Tree_k(\ux;y)$ must satisfy. The operadic machinery then gives the desired map of operads $\underline{\Sigma}^n\otimes \Ai\to \mathcal{MT}^{(X,f,\X)}$. Composing with the operad action from Theorem \ref{theorem:MT acting on Morse homology}, this recovers the suspended $\Ai$-action as in Theorem \ref{theorem : Ai structure}.

\subsection{Structure of this paper}

In Section \ref{sec : perturbed Morse}, we recall the construction of perturbed Morse gradient trees. In Section \ref{sec:Morse hom with DG}, we set our conventions of Morse homology with DG coefficients. In Section \ref{sec:geom}, we explain the geometric ingredients that produce the binary operation and compare it to the generalised Chas-Sullivan product. In Section \ref{sec:coloured ops}, we recall the relevant background about coloured operads. In Section \ref{sec:Coloured Ainfty}, we construct the various versions of the $\Ai$-operad and explain how they relate. In Section \ref{sec:top op}, we construct the homotopy coherent action on the fibres. In Section \ref{sec:main thm}, we explain how to get an operad action of $\mathcal{MT}^{(X,f,\X)}$ on $C_*(X ,C_*(F))$ and produce the map of operads $\underline{\Sigma}^n\otimes \Ai\to \mathcal{MT}^{(X,f,\X)}$ using the orientation of the moduli spaces $\Tree_k(\ux;y)$ that we produce in Appendix \ref{sec:or trees}. In Appendix \ref{appendix: The BDHO toolset in our convention}, we compare our convention for Morse homology with DG coefficients to the conventions used in the literature.

\subsection{Conventions}
Throughout this article, $X$ denotes a smooth closed oriented manifold of dimension $n$. Moreover, $f\colon X\to \R$ denotes a Morse function on $X$ and $\xi$ denotes a pseudo-gradient vector field on $X$ such that $(f,\xi)$ is a Morse-Smale pair. We denote 
\begin{align*}
    \mathcal{P}X\coloneqq\{\gamma\colon [0,a]\to X\mid a\geq 0, \gamma\text{ is piecewise smooth}\}
\end{align*}
for the piecewise smooth Moore paths in $X$. We  denote by $\ev_0$ and $\ev_1$ the evaluation at the starting point and the endpoint of a Moore path, respectively.

For any $x,y\in X$, we denote $\mathcal{P}_{x\to y}X\coloneqq\{\gamma\in \mathcal{P}X\mid \ev_0(\gamma)=x, \ev_1(\gamma)=y\}$.

If $B,C \subseteq X$, we use the notations $$\mathcal{P}_{B \to C} X \coloneqq \displaystyle \bigcup_{b \in B, c \in C} \mathcal{P}_{b \to c} X\subseteq \mathcal{P}X, \quad \mathcal{P}_{b \to C} X \coloneqq \mathcal{P}_{\{b\} \to C} X \ \textup{and} \ \mathcal{P}_{B \to c} X \coloneqq \mathcal{P}_{B \to \{c\}} X.$$

We by $I\coloneqq [0,1]$ the unit interval

\subsection*{Acknowledgements}

We thank Mihai Damian, Guillaume Laplante-Anfossi, Alexandru Oancea, Maximilian Stegemeyer and Nathalie Wahl for helpful discussions and feedback on earlier versions. The first author was supported by the Danish National Research Foundation through the Copenhagen Centre for Geometry and Topology (DNRF151) and the Swiss National Science Foundation.

\subsection*{AI usage statement}
AI was only used for literature review and proofreading.

\tableofcontents

\newpage 

\section{Perturbed Morse gradient trees}\label{sec : perturbed Morse}

In this section, we recall the definition of perturbed Morse gradient trees. Intuitively, such a tree is given by a metric tree embedded in $X$ that follows the pseudo-gradient $\xi$ up to small perturbations. They were originally introduced by Abouzaid in \cite{Abo09,Abo11}. We mostly follow Mazuir's exposition in \cite{mazuir2022higheralgebraainftyomega}.

Let $f\colon X \to \R$ be a Morse function. We denote by $\|x\|=n-|x|$ the cohomological degree of a critical point $x\in \Crit(f)$.

If $(x_1, \dots, x_k) \in \Crit(f)^k$, we denote $\ux = (x_1, \dots, x_k)$ and $$\| \ux \| = \sum_{i=1}^k \|x_i\|.$$
Let $k\geq 1$ and $x_i,y \in \Crit(f)$ for all $i \in \{1, \dots, k\}$. The \textit{moduli space of perturbed Morse gradient trees} $\Tree_k(\ux;y)$ is built from \textit{stable metric ribbon trees}, which are closely related to associahedra, and therefore are well suited to encode $\Ai$-algebra structures.

These spaces have been shown by Abouzaid \cite[Lemma 2.13]{Abo11} to be compact, oriented manifolds with boundary of dimension $$k-2 + \|y\| - \sum_{1 \leq i \leq k} \|x_i\|.$$

We now define the moduli spaces of perturbed gradient trees. We use the definitions of Mazuir in \cite{mazuir2022higheralgebraainftyomega}.

\subsection{Stable metric ribbon trees}\label{subsec:stable metric ribbon trees}

\begin{figure}[!ht]
    \centering
    \includegraphics[scale = 0.15]{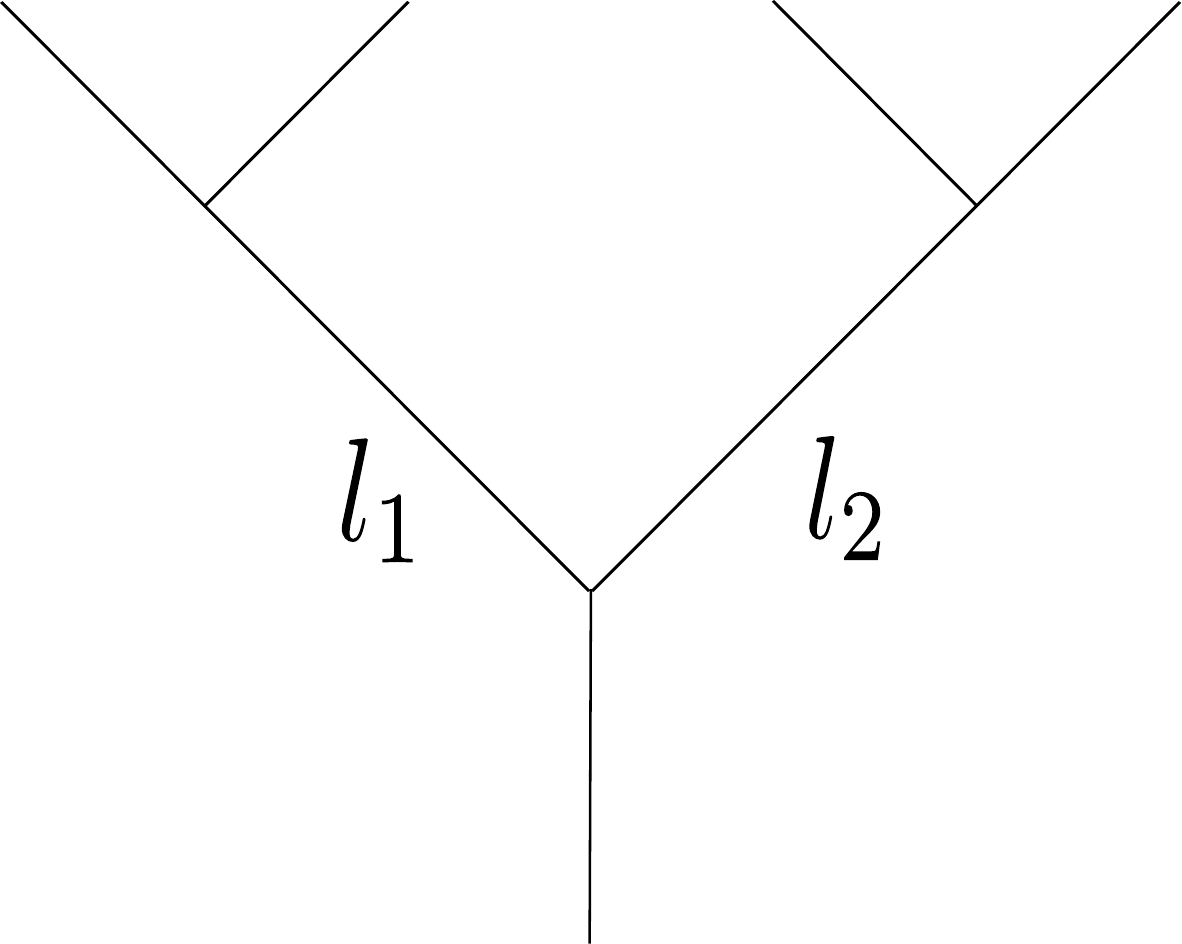} \hspace{5em}
    \includegraphics[scale = 0.15]{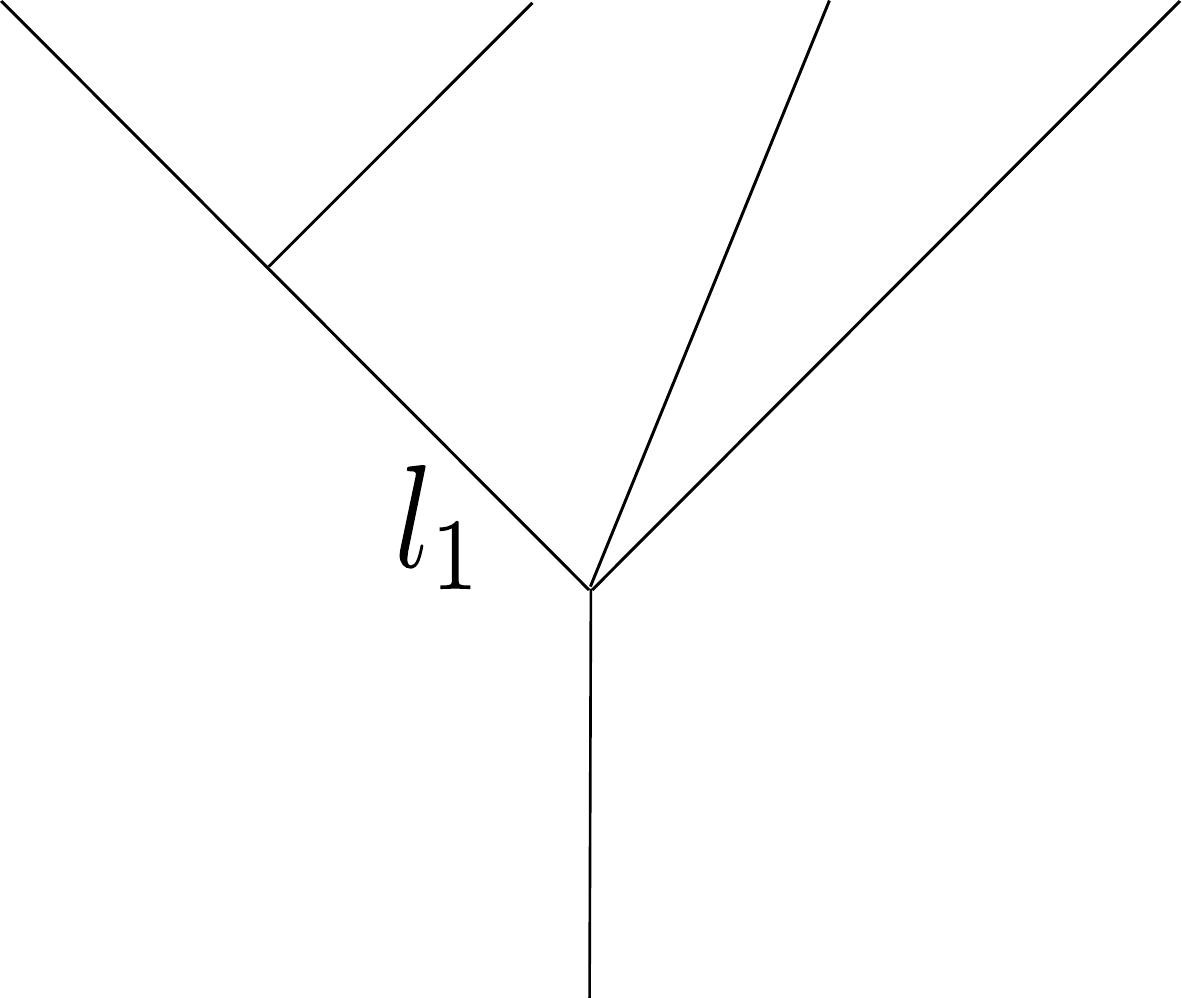}
    \hspace{5em}
    \includegraphics[scale = 0.15]{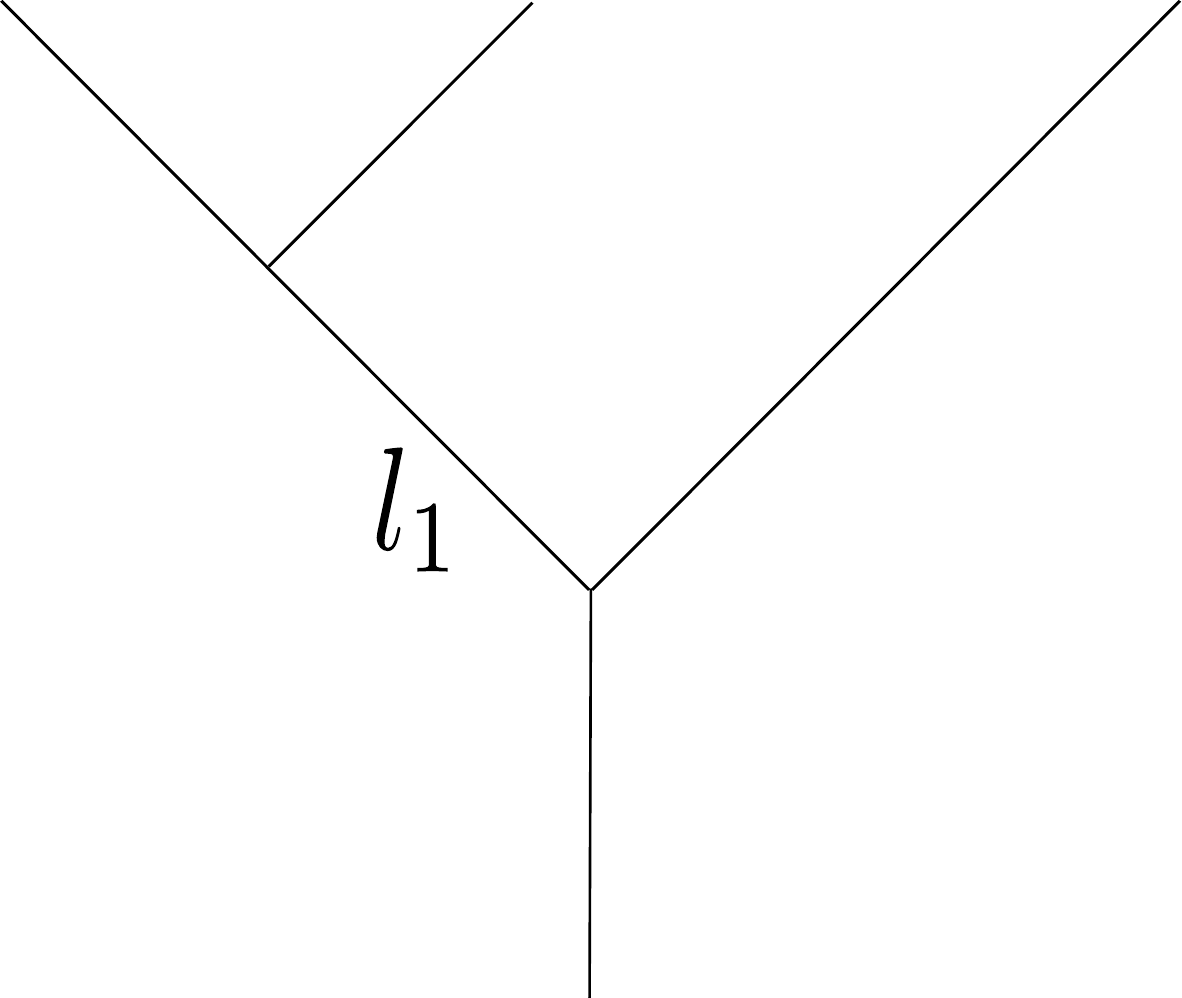}
    \caption{Examples of stable metric ribbon trees.}
    \label{fig: Examples of stable metric ribbon trees}
\end{figure}

\begin{defi}[{\cite[Definition 3.1.1]{mazuir2022higheralgebraainftyomega}}]
    A \textit{stable ribbon tree} $\tau$ is a tree whose inner vertices are at least trivalent together with a cyclic ordering on the edges at each vertex of the tree and a distinguished \enquote{\textit{outgoing}} external vertex.

    Denote $\E(\tau)$ the set of internal edges of $\tau$. If we add the data of a length $l_e \in [0, + \infty]$ for any internal edge $e \in \E(\tau)$, we call $\left(\tau, (l_e)_{e \in \E(\tau)}\right)$ a \textit{stable metric ribbon tree}. 

    The \textit{moduli space of metric trees}, $\overline{\mathcal{T}_k}$, is given by the space of all stable metric ribbon trees with $k+1$ external edges modulo the relation that a stable metric ribbon tree with internal edge $e$ of length $0$ is identified with the corresponding stable metric ribbon tree with $e$ collapsed (see Figure \ref{fig: Examples of stable metric ribbon trees}).
\end{defi}

\begin{figure}[!ht]
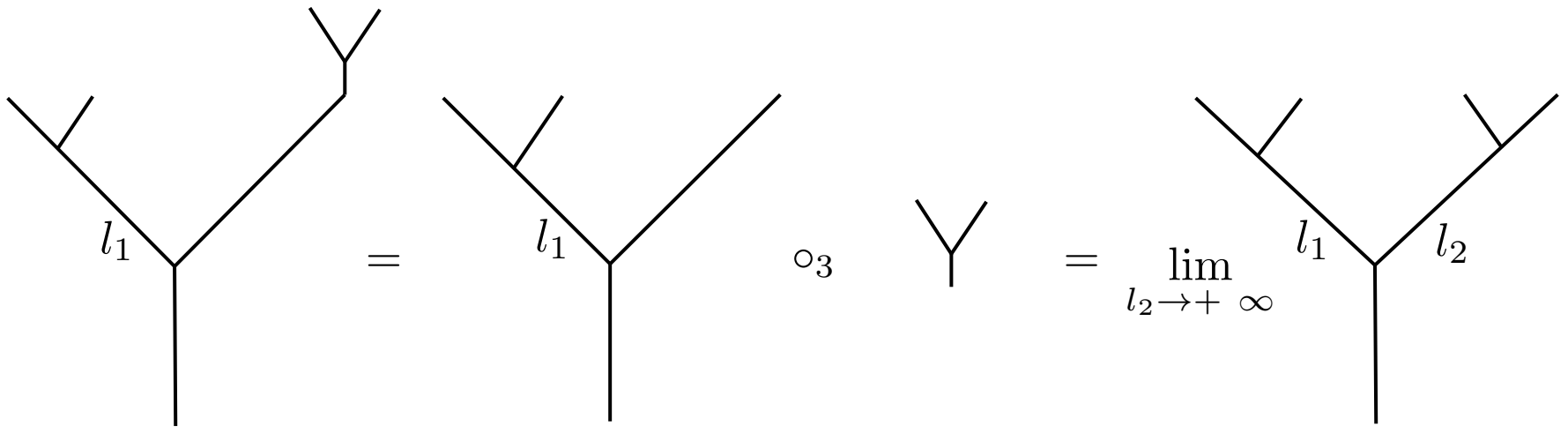

    \centering
    \drawingbrokentree
    \caption{A singular tree }
    \label{fig: singular tree}
\end{figure}

This space has the structure of a $(k-2)$-dimensional CW-complex with interior $\mathcal{T}_k$ given by trees with edges of finite length 
\begin{equation}\label{eq : gluing of trees}
    \mathcal{T}_l \circ_i \mathcal{T}_m \subset \partial \overline{\mathcal{T}}_k, 
\end{equation} 

where $l,m \geq 2$, $l+m = k+1$, $i \in \{1, \dots,l\}$ and $\circ_i$ is the operation where the tree in $\mathcal{T}_m$ connects to the $i$-th incoming edge of a tree in $\mathcal{T}_l$.
We interpret such a tree as a tree with $k$ incoming edges having an internal edge of infinite length. For example, the singular tree represented in Figure \ref{fig: singular tree} is seen as the limit when $l_2 \to + \infty$ of the first tree represented in Figure \ref{fig: Examples of stable metric ribbon trees}.

\subsection{Perturbation data}

\begin{defi}[{\cite[Definition 6.2.1]{mazuir2022higheralgebraainftyomega}}]
    Let $T\coloneqq(\tau,(l_e)_{e \in \E(\tau)}) \in \overline{\mathcal{T}}_k$. A \textit{choice of perturbation data $\X$ on $T$} consists of the following data:
    \begin{enumerate}[(i)]
        \item a vector field $$[0,l_e] \times X \overset{\mathbb{X}_e}{\to} TX,$$ that vanishes away from $[1,l_e -1]$, for every internal edge $e$ of $\tau.$
    
        \item a vector field $$[0,+ \infty) \times X \overset{\X_{e_0}}{\to} TX,$$ that vanishes away from $[0, 1]$, for the outgoing edge $e_0$ of $\tau.$
    
        \item a vector field $$(- \infty, 0] \times X \overset{\X_{e_i}}{\to} TX,$$ that vanishes away from $[-1, 0]$, for every incoming edge $e_i$, $ i \in \{1,\dots, k\}$ of $\tau$. 
    \end{enumerate}    

    A continuous choice for all $T\in \mathcal{T}_k$ is called a \textit{choice of perturbation data}.
\end{defi}

If $(T=(\tau,(l_e)),\X)$ is a stable metric ribbon tree together with a choice of perturbation data, we write $D_e$ to denote the segment $[0,l_e]$ or $(-\infty,0]$ or $[0, + \infty)$ for $e$ an edge of $\tau$.

\begin{defi}[{\cite[Definition 6.2.2]{mazuir2022higheralgebraainftyomega}}] \label{def : perturbed Morse gradient trees}
A \textit{perturbed Morse gradient tree} associated with $(T,\X)$ is the data of a smooth map $\gamma_e \colon D_e \to X$ for each edge $e$ of $T$ such that $\gamma_e$ is a trajectory of the perturbed pseudo-gradient $\xi + \X_e$ and such that the endpoints of these trajectories coincide as prescribed by the edges of the tree $T.$
\end{defi}

The vector fields on the external edges are all equal to $\xi$ away from a length 1 segment. Therefore, for every $i \in \{1, \dots,k\}$, the trajectory $\gamma_{e_i} \colon (-\infty, 0] \to X$ of the incoming edge $e_i$ converges to a critical point $x_i \in \Crit(f)$ as $t\to -\infty$. The trajectory $\gamma_{e_0} \colon [0, +\infty) \to X$ converges towards a critical point $y \in \Crit(f)$ as $t\to +\infty$.

\begin{figure}[!ht]
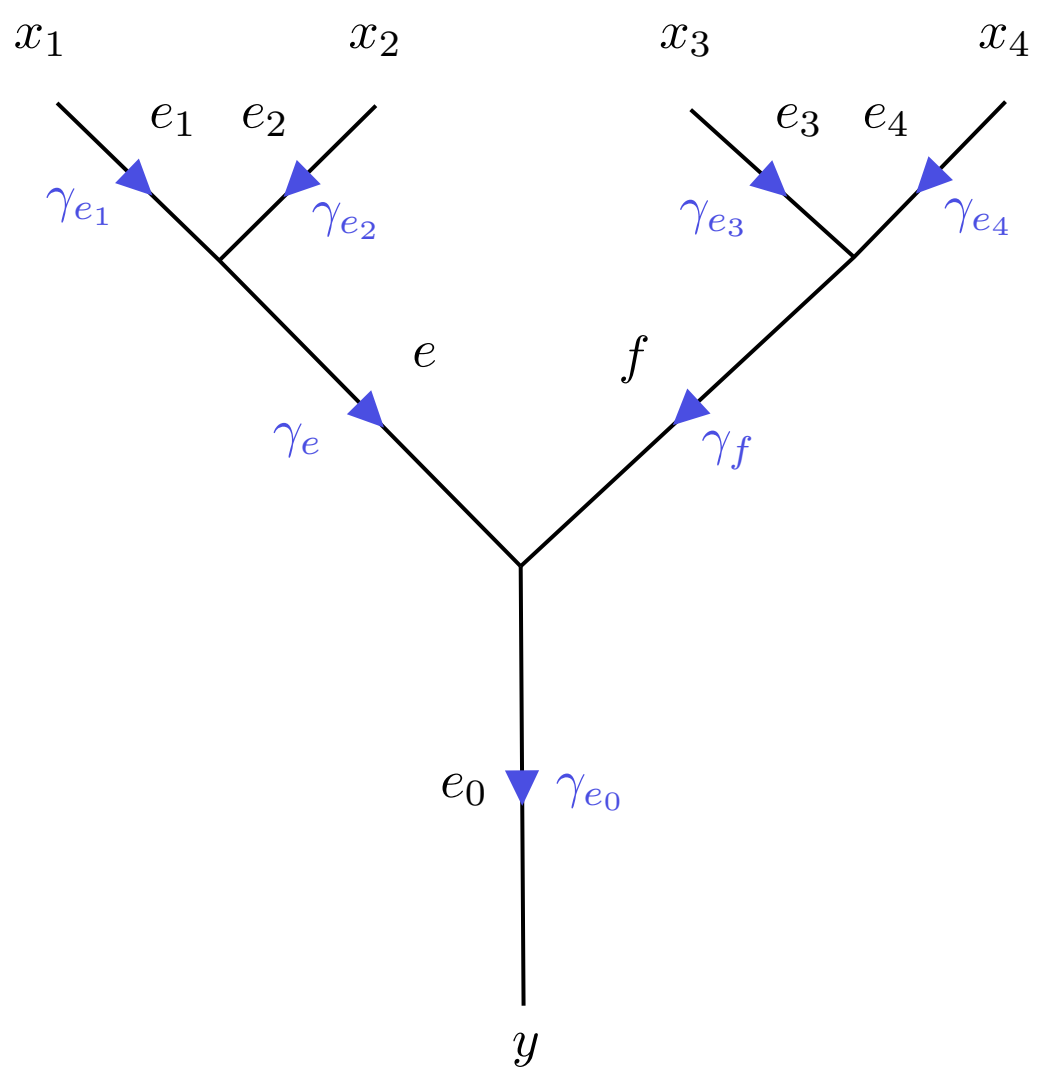

    \centering
    \drawingperturbedgradienttree
    \caption{A perturbed Morse gradient tree}
    \label{fig: A perturbed Morse gradient tree}
\end{figure}

\begin{defi}[{\cite[Definition 6.3.2]{mazuir2022higheralgebraainftyomega}}]\label{def:Tree_k(x,y)}
    Let $\X$ be a choice of perturbation data, let $y \in \Crit(f)$ and $\ux=(x_1, \dots,x_k) \in \Crit(f)^k$. Denote $\mathcal{T}^{\X}_k(\ux;y)$ the space of perturbed Morse gradient trees such that the trajectory of the outgoing edge converges to $y$ and the trajectory of the $i$-th incoming edge converges to $x_i$ for $i \in \{1, \dots,k\}$.
\end{defi}

Abouzaid and Mazuir explained that a perturbed gradient tree is uniquely determined by  the data of the time $-1$ points on its incoming edges plus the time $1$ point on its outgoing edge and used this characterisation to see  $\mathcal{T}^{\X}_k(\ux;y)$ as a fibre product. Since we heavily use this viewpoint, let us give a slightly modified version of Mazuir's construction.

Let $T \in \mathcal{T}_k$. For two points in the tree $s\in D_e$ and $t\in D_{e'}$, there exists a unique non-self-crossing path connecting $s$ to $t$: $s$ to $ s'\in \partial D_e\cap D_e$, $0$ to $l_{f_1}$ in $D_{f_1}$,..., $0$ to $l_{f_a}$ in $D_{f_a}$, $l_{f'_1}$ to $ 0$ in $D_{f_1'}$,..., $l_{f_b'}$ to $ 0$ in $t'\in  \partial D_{e'}\cap D_e$ and $l_{e'}$ to $ t$ in $D_{e'}$ (see Figure \ref{fig: A perturbed Morse gradient tree}).

\begin{ex}
In Figure \ref{fig: A perturbed Morse gradient tree}, $s$ in $D_{e_1}$ is connected to $t\in D_{e_3}$ via the path $s\to 0$ in $D_{e_1}$, $0\to l_e$ in $D_e$, $l_f\to 0$ in $D_f$ and $0\to t$ in $D_{e_3}$.
\end{ex}

\begin{defi}(\cite[\S 6.2]{mazuir2022higheralgebraainftyomega})\label{defi : Phi_X}
    Define the following diffeomorphisms $X\to X$:

\begin{itemize}
    \item for $s,t\in D_e$, let $\phi_{T,e}^{s,t}$ denote the flow along the time-dependent vector field $\xi+ X_e$ from $s$ to $t$;
    \item for $e_i$ an incoming edge, define $\phi_{e_i,\X} \coloneqq \phi^{-1, 0}_{e_i,\X}\colon X \to X$;
    \item for $e_0$ the outgoing edge, define $\phi_{e_0,\X} \coloneqq \phi^{0, 1}_{e_0,\X}\colon X \to X$;
    \item for $e$ an internal edge of $T$, denote $\phi_{e,\X} \coloneqq \phi^{0,l_e}_{e,\X}$;
    \item for $s\in D_e$ and $t\in D_{e'}$ with path connecting from $s$ to $t$ as above, define
    \begin{align*}
        \phi_{T,\X}^{s,t}:=\phi_{e',\X}^{t',t}\circ \phi_{f_b',\X}^{-1}\circ \dots \circ \phi_{f'_1,\X}^{-1}\circ \phi_{f_a,\X}\circ \dots \circ \phi_{f_1} \circ \phi_{e,\X}^{s,s'}
    \end{align*}
    the flow from $s$ to $t$.
\end{itemize}

 Then, for $i \in \{1, \dots, k\}$, define $\phi_{i,\X} \coloneqq \phi_{\X}^{s,t}\colon X \to X$ for $s\coloneqq-1\in D_{e_i}$ and $t\coloneqq1\in D_{e_0}$ the map that sends the time $-1$ points of $e_i$ to the time $1$ point of $e_0$ following the only non-self-crossing path connecting $e_i$ to $e_0$.
 
    Let $\Phi_{T,\X}\colon X^k \to X^k$ be the diffeomorphism isotopic to the identity defined by $$\Phi_{T,\X} =\left( \phi_{1,\X}, \dots, \phi_{k,\X} \right).$$

Using the decomposition by the unstable manifolds, this map induces the map $$\Phi_{\ux}\coloneqq \Phi_{\ux,\X}\colon W^u(\ux) \times \mathcal{T}_k \to X^k,$$ whose restriction to $T \in \mathcal{T}_k $ is $\Phi_{T,\X}.$
\end{defi}

\begin{prop}[{\cite[Proposition 6.3.3]{mazuir2022higheralgebraainftyomega}}] \label{prop : perturbed gradient trees as a fibre product}
    Let $\ux,y \in \Crit(f)$. $$\mathcal{T}^{\X}_k(\ux;y) \cong \Phi^{-1}_{\ux,\X}\left( \Delta_k(W^s(y))\right),$$
    where $\Delta_k \colon X \to X^k$ is the thin diagonal $\Delta_k(a) = (a, \dots, a)$.
\end{prop}

\begin{figure}[h!]
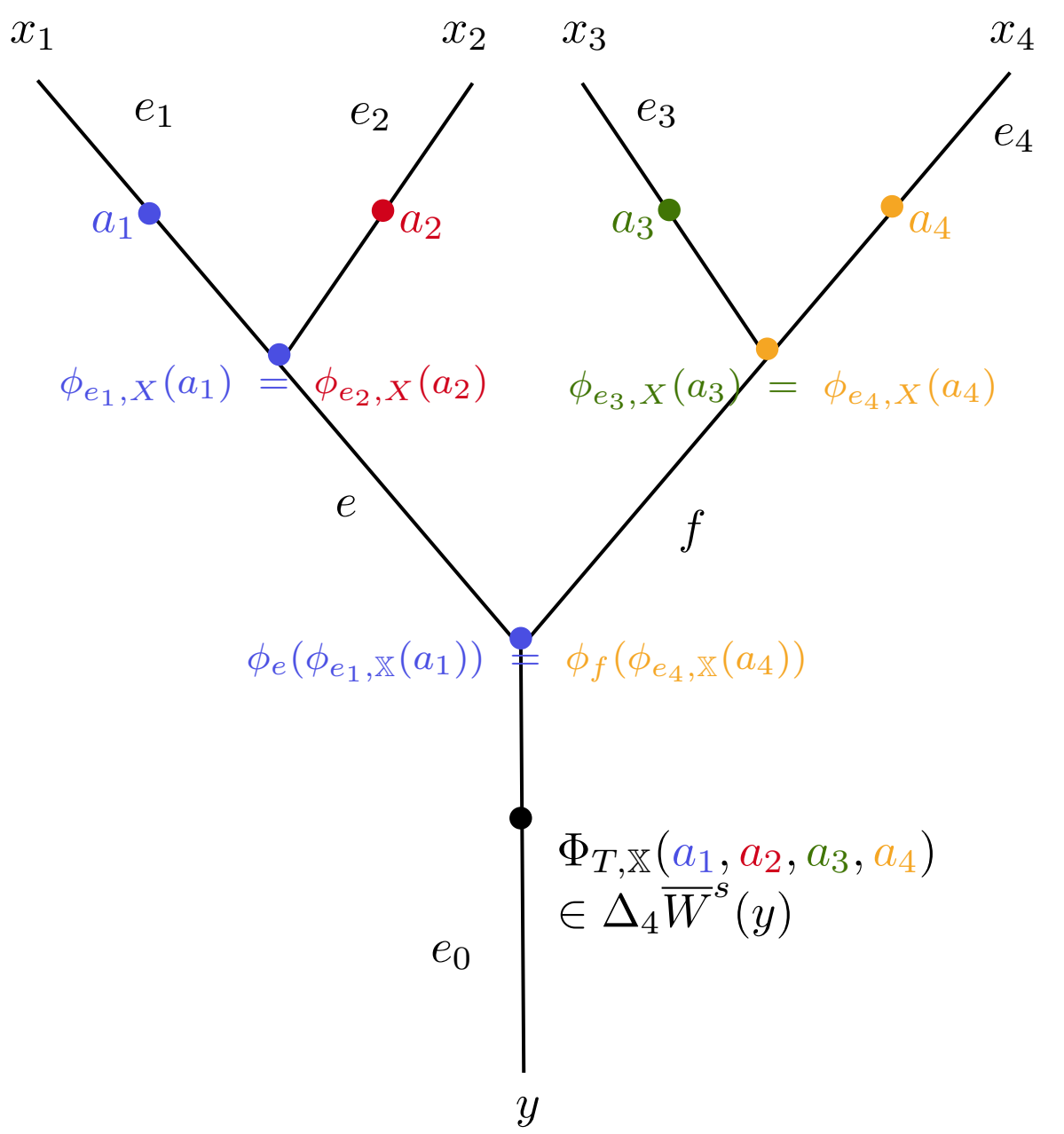

    \centering
    \flowtree
    \caption{A perturbed Morse gradient tree $(a_1,a_2,a_3,a_4) \in \mathcal{T}^{\X}_4(x_1,x_2,x_3,x_4;y)$.}
    \label{fig:flowtree}
\end{figure}

\begin{defi}(\cite[Definition 6.4.4]{mazuir2022higheralgebraainftyomega})
    An \textit{admissible choice of perturbation data} $(\X_k)_{k \geq 2}$ is a smooth family of perturbation data for elements of $\mathcal{T}_k$ for every $k \geq 2$, which is compatible with the gluing \eqref{eq : gluing of trees} and such that $\Phi_{\ux,\X}$ is transverse to $\Delta_k\left(W^s(y)\right)$ for all $\ux,y \in \Crit(f)$.
\end{defi}

Mazuir proved \cite[Theorem 6.4.5]{mazuir2022higheralgebraainftyomega} that such data exist.

The space $\mathcal{T}_k(\ux;y)$ admits a canonical compactification $\Tree_k(\ux;y)$ by adding singular gradient trees.

\begin{lemme}[{\cite[Lemma 2.13]{Abo11}, \cite[Lemma B.12]{Abo09}, \cite[\S6.4]{mazuir2022higheralgebraainftyomega}}]\label{lem:gradient trees are compact manifolds}
    For a generic choice of admissible perturbation data, the spaces $\Tree_k(\ux;y)$ are compact manifolds with boundary of dimension $$k-2 + \|y\| - \sum_{ i=1}^k \|x_i\| = k-2 + (1-k)n + \sum_{i=1}^k |x_i| - |y|.$$
        The boundary of $\Tree_k(\ux;y)$ is

    \begin{multline}\label{eq : boundary of tree}
        \partial \Tree_k(\ux;y) = \overset{\blue{(1)}}{\bigcup_w   \trajb{w,y}\times  \Tree_k(\ux;w) }  \cup \overset{\blue{(2)}}{\bigcup_{\substack{z\\ 1 \leq i \leq k}} \Tree_k((x_1, \dots,x_{i-1}, z, x_{i}, \dots, x_k);y)}\times  \trajb{x_i,z}  \\
        \cup \underset{\blue{(3)}}{\bigcup_{\substack{l+m=k+1\\ l,m \geq 2}} \bigcup_{\substack{c\\ 1 \leq i \leq m}}   \Tree_l((x_i, \dots, x_{i+m-1});c) }\circ_i\Tree_m((x_1, \dots, x_{i-1}, c, x_{i+m}, \dots, x_k);y),
    \end{multline}
    where $\trajb{x,y}$ denotes the broken Morse trajectories from $x$ to $y$.
\end{lemme}

The first two components \blue{(1)} and \blue{(2)} come from the breaking of an external edge, i.e., the trajectory $\gamma_e$ associated with an external edge breaks into a broken trajectory. The last component \blue{(3)} comes from the breaking of an internal edge.
We introduce the notation $\Tree_1(x;y) \coloneqq\overline{\mathcal{L}}(x,y)$ for $x,y \in \Crit(f)$ so that \eqref{eq : boundary of tree} reads
\begin{equation*}
   \partial \Tree_k(\ux;y) = \bigcup_{\substack{l+m=k+1\\ l,m \geq 1}} \bigcup_{\substack{c\\ 1 \leq i \leq l}}    \Tree_l((x_1, \dots, x_{i-1}, c, x_{i+m}, \dots, x_k);y)\circ_i\Tree_m((x_i, \dots, x_{i+m-1});c).
\end{equation*}

\section{Morse homology with DG coefficients and fibrations}\label{sec:Morse hom with DG}

In this section, we define a Morse homology with DG coefficients in the fibres of a fibration. This homology allows us to extend the use of Morse theory to study spaces that admit a fibration over a closed, smooth manifold. We recall the notion of fibration in \S \ref{subsubsec: Fibrations and transitive lifting functions}.

Brown proved in \cite{Br69} that if $\pi\colon E \to X$ is a fibration over a closed, connected and smooth manifold, then the homology of $E$ can be recovered as the homology of a tensor product of the singular complex of $X$ and the singular complexes of the fibres of $\pi$, twisted by a collection of \textit{twisting cochains}. Barraud, Damian, Humilière and Oancea constructed in \cite{BDHO23} a Morse homology with DG coefficients that computes the homology of $E$ as the homology of the tensor product of the Morse complex of $X$ with the cubical complex of the fibres, twisted by a \textit{Barraud-Cornea twisting cocycle}.
The construction presented in this section differs from that of Barraud, Damian, Humilière and Oancea only by a few conventions. We compare these two constructions in Appendix \ref{appendix: The BDHO toolset in our convention}.

\subsection{The Barraud-Cornea twisting cocycle}

We consider a Morse-Smale pair $(f,\xi)$ on $X$, that is to say that $\xi$ is a pseudo-gradient adapted to $f$ that satisfies the \textit{Smale} condition (see \cite[\S2]{AD14}). We define a family of chains $$\left\{ m_{x,y} \in C_{|x|-|y|-1}(\mathcal{P}_{x\to y} X), \ x,y \in \Crit(f)\right\}$$ in a  manner analogous to how Barraud, Damian, Humilière and Oancea define the Barraud-Cornea twisting cocycle in \cite{BDHO23} except that we do not perform the last step that consists in collapsing a tree $\mathcal{Y}$ connecting the critical points of $f$ to a basepoint.

We fix a choice $([\wb^u(x)])_{x \in \Crit(f)}$ of orientation for each unstable manifold $W^u(x)$. This gives an orientation on the moduli space of Morse trajectories $\overline{\mathcal{L}}(x,y)$ using the following orientation rule: it holds as a local orientation at $x$
\begin{equation}\label{eq: or L(x,y)}
    \left([ \trajb{x,y}], -\xi, [\wb^u(y)]\right) = [\wb^u(x)].
\end{equation}
See \cite[\S 5.2.1]{BDHO23} for more details.

\begin{defi}[{\cite[Definition 5.2.7]{BDHO23}}]\label{defi : representing chain system traj}
    A \textit{representing chain system} of the moduli spaces of trajectories is a collection $\{s_{x,y} \in C_{|x|-|y|-1}(\trajb{x,y}) , \ x,y \in \textup{Crit}(f)\}$ such that 

    \begin{enumerate}
        \item each $s_{x,y}$ is a cycle relative to the boundary and represents the fundamental class $[\trajb{x,y}]$;\\
        \item each $s_{x,y}$ satisfies 

        \begin{equation}\label{eq : relation sys repr}
            \partial s_{x,y} = \sum_z (-1)^{|x|-|z|} s_{x,z} \times s_{z,y}.
        \end{equation}
    \end{enumerate}
\end{defi}

Such a representing chain system always exists and is constructed inductively. It is not unique, although all representing chain systems are homologous in a sense defined in \cite[Proposition 5.2.8]{BDHO23}.

\begin{defi}\label{defi : DG Morse data}
    Let $\Xi = (f,\xi, o,s_{x,y})$ be the data of a Morse-Smale pair $(f,\xi)$, an orientation $o_x\coloneqq[\wb^u(x)]\in H_{|x|}(\wb^u(x))$ for every $x \in \Crit(f)$ and a representing chain system $\{s_{x,y}\}_{x,y \in \Crit(f)}$. 

    The data $\Xi$ is called a \textit{set of DG Morse data}. 
\end{defi}

Since the goal of this article is to build an operation on chains, we fix the set $\Xi$ of DG Morse data this article and only refer to it as a way to recall the notations and prior choices that had to be made.

\begin{defi}\label{def : parametrisation map}
  For $x,y \in \Crit(f)$, define $\Gamma_{x,y}\colon \trajb{x,y} \to \mathcal{P}_{x \to y}X$ by $$\Gamma_{x,y}(\lambda) =\left( \begin{array}{ccc}[0, f(x)-f(y)]  & \to & X,  \\
      t & \mapsto & (\lambda \cap f^{-1}(f(x)-t))
  \end{array}\right)$$ the parametrisation induced by the values of $f$.    
\end{defi}

\begin{defi}\label{def : mCR}
  For $x,y \in \Crit(f)$, denote $m_{x,y} = \Gamma_{x,y,*}(s_{x,y}) \in C_{|x|-|y|-1}(\mathcal{P}_{x\to y} X)$. The family $$\left\{m_{x,y}  \in C_{|x|-|y|-1}(\mathcal{P}_{x\to y} X), \ x,y \in \Crit(f)\right\}$$ is called the Barraud-Cornea twisting cocycle associated with $\Xi$ and it satisfies the Maurer-Cartan equation 

  \begin{equation}\label{eq : MC CR}
      \partial(m_{x,y}) = \sum_{z\in \Crit(f)} (-1)^{|x|-|z|} m_{x,z}\cdot m_{z,y},
  \end{equation}

  where the product denotes the concatenation of paths at the chain-level.
\end{defi}

\subsection{Fibrations and transitive lifting functions}\label{subsubsec: Fibrations and transitive lifting functions}

The system of coefficients considered in our construction of Morse homology with DG coefficients consists only of fibres of a (Hurewicz) fibration.

\begin{defi}
    A \textit{Hurewicz fibration} $ E \overset{\pi}{\to} X$ is a continuous map such that there exists a map $\phi \colon E \ _{\pi}\!\!\times_{\ev_0} \mathcal{P} X \to \mathcal{P}E$ such that $\pi \circ \phi = pr_2$ and $\ev_0 \circ \phi = pr_1.$
\end{defi}

\begin{defi}\label{defi: TLF}
    A \textit{transitive lifting function} is a map $$\Phi \colon E \ftimes{\pi}{\ev_0} \mathcal{P} X \to E$$ such that $\pi \circ \Phi = \ev_1 \circ pr_2$, such that for any $e\in E$,  $\Phi(e,\pi(e)) = e$ where $\pi(e)$ here denotes the constant path in $\pi(e) \in X$ and such that for every $\gamma, \delta \in \mathcal{P} X$ such that $\ev_1(\gamma)= \ev_0(\delta)$,
$$\Phi(\Phi(e,\gamma),\delta) = \Phi( e, \gamma \cdot  \delta).$$
\end{defi}

It is proved in \cite[Proposition 5.5]{DK69}  that every fibration is fibre homotopy equivalent to a fibration that admits a transitive lifting function and we therefore assume that any fibration considered is endowed with a transitive lifting function.

\begin{rem}
    The reader may be more accustomed to the definition of a Hurewicz fibration as a map that satisfies the \textit{homotopy lifting property} with regard to all spaces. The two definitions are equivalent and we only use the transitive lifting function in this article.
\end{rem}

\subsection{The Morse complex with DG coefficients}\label{subsubsec: Morse complex DG}

We consider a fibration $\fibration$ over $X$ and we denote by $\Phi\colon E \ftimes{\pi}{\ev_0}\mathcal{P}X \to E$ a transitive lifting function associated with this fibration. Let $\Xi = (f,\xi, o,s_{x,y})$ be a set of DG Morse data.

We now define a complex with coefficients in the fibres of this fibration using the twisting cocycle $m_{x,y} \in C_{|x|-|y|-1}(\mathcal{P}_{x \to y} X)$ defined in Definition \ref{def : mCR}.

\begin{defi}\label{def:Morse hom with DG coeffs}
    Define $\deffct{\Phi^{op}}{\mathcal{P}X \ftimes{\ev_0}{\pi} E}{E,}{(\gamma,a)}{\Phi(a,\gamma)}.$ Consider the complex
    $$C_*(X,C_*(F)) = \bigoplus_{x\in \Crit(f)} \Z\cdot x \otimes C_*(F_x) $$
    with the differential defined on the generators by $$\partial(x \otimes \alpha) = \sum_{y \in \Crit(f)} (-1)^{|x||y|} y \otimes \Phi^{op}_*(m_{x,y}\otimes \alpha) + (-1)^{|x|} x \otimes \partial_F\alpha.$$ The homology of this complex is denoted by $H_*(X,C_*(F))$.
\end{defi}

Since the set of DG Morse data $\Xi$ is fixed throughout this article, we often omit to mention it and denote this complex $C_*(X,C_*(F))$.

\begin{lemme}
    The linear map $\partial\colon C_*(X,C_*(F)) \to C_{-1+*}(X,C_*(F))$ defines a differential.
\end{lemme}

\begin{proof}
    Because $|m_{x,y}|+|y| - |x| =-1$, $\partial$ has degree $-1$. We now check that $\partial^2=0.$

    Let $x\in \Crit(f)$ and $\alpha \in C_*(F_x).$ We compute 
    
    \begin{align*}
        \partial^2(x \otimes \alpha) = &\ \partial \left( \sum_y (-1)^{|x||y|}y \otimes \Phi^{op}_*(m_{x,y}\otimes \alpha) + (-1)^{|x|} x \otimes \partial_F\alpha  \right) \\
        =&\  \sum_{y,z} (-1)^{|x||y| + |y||z|} z \otimes \Phi^{op}_*(m_{y,z} \otimes \Phi^{op}_*(m_{x,y}\otimes \alpha))\\
        &+ \sum_y (-1)^{|x||y| + |y|}y \otimes \partial_F\Phi^{op}_*(m_{x,y}\otimes \alpha) \\
        &+ (-1)^{|x|} \sum_y (-1)^{|x||y|}y \otimes \Phi^{op}_*(m_{x,y}\otimes \partial_F\alpha) + x\otimes \underbrace{\partial_F^2\alpha}_{=0}.
        \end{align*}

 Since $\Phi$ is transitive, for every $e\in E$ and $\gamma,\delta \in \mathcal{P}X$ such that $\ev_1(\gamma) = \ev_0(\delta)$, it holds that  $$\Phi^{op}(\delta,\Phi^{op}(\gamma,e)) = \Phi(\Phi(e,\gamma),\delta) = \Phi(e, \gamma\cdot \delta) = \Phi^{op}(\gamma\cdot\delta,e).$$ Adding the Koszul signs at the chain-level, we obtain \begin{equation*}
    \Phi^{op}_*(\beta \otimes \Phi^{op}_*(\omega \otimes \alpha)) = (-1)^{|\beta||\omega|} \Phi^{op}_*(\omega\cdot \beta \otimes \alpha).\end{equation*}

We now resume the computation.

        \begin{align*}
        \partial^2(x \otimes \alpha) =&\ \sum_{y,z}(-1)^{|x||y| + |y||z|} (-1)^{(|x|-|y|-1)(|y|-|z|-1)} z \otimes \Phi^{op}_*(m_{x,y} \cdot m_{y,z} \otimes \alpha) \\
        &+ \sum_y (-1)^{|x||y| + |y|}y \otimes \Phi^{op}_*(\partial(m_{x,y})\otimes \alpha) + \sum_y (-1)^{|x||y| + |y| + |x|+|y|-1} y \otimes \Phi^{op}_*(m_{x,y}\otimes \partial_F\alpha) \\
        &+ (-1)^{|x|} \sum_y (-1)^{|x||y|}y \otimes \Phi^{op}_*(m_{x,y}\otimes \partial_F\alpha)\\
        {=}& - \sum_{y,z}(-1)^{|x||z| + |z| + |x|+|y|} z \otimes \Phi^{op}_*(m_{x,y} \cdot m_{y,z} \otimes \alpha) \\
        &+ \sum_{y,z} (-1)^{|x||y| + |y| + |x| + |z|}y \otimes \Phi^{op}_*(m_{x,z} \cdot m_{z,y}\otimes \alpha) - \sum_y (-1)^{|x||y| + |x|} y \otimes \Phi^{op}_*(m_{x,y}\otimes \partial_F\alpha) \\
        &+ (-1)^{|x|} \sum_y (-1)^{|x||y|}y \otimes \Phi^{op}_*(m_{x,y}\otimes \partial_F\alpha)\\
        =&\ 0
    \end{align*}

    We use equation \eqref{eq : MC CR} and simplify the signs in the second equality.
\end{proof}

\section{Geometric construction}\label{sec:geom}

Let $\Xi = (f,\xi,o,s_{x,y})$ be a set of DG Morse data (see Definition \ref{defi : DG Morse data}) and let $\fibration$ be a fibration over $X$ together with a transitive lifting function $\Phi \colon E \ftimes{\pi}{\ev_0} \mathcal{P}X \to E$. For $a \in X$, we denote $F_a = \pi^{-1}(\{a\}).$

As often in Morse homology with DG coefficients, we extract from  moduli spaces of compact and oriented manifolds of interest (in this instance $\Tree_k(\ux;y)$), a chain map using a family of chains called \textit{cocycle} constructed by the following process:

\begin{itemize}
    \item We build a \textit{representing chain system} (Lemma \ref{lemme : representing chain system for trees}). This is a family of chain representatives of the fundamental classes of the moduli spaces that is adapted to their boundaries and orientations. More precisely, this family satisfies a boundary equation \eqref{eq: boundary repr system of trees} that encodes the boundaries of the moduli spaces and their orientations.

    \item We choose an evaluation map  (Lemma \ref{lemme : evaluation maps trees}) adapted to the boundaries to transform this representing chain system into a family of chains of paths called a \textit{cocycle}.
    
    \item This cocycle then acts on the system of coefficients (here $(F_{x_1} \times \dots \times F_{x_k})_{\ux \in \Crit(f)^k}$) using the transitive lifting function $\Phi$.
\end{itemize}

In this section, we apply this process to define the \textit{tree cocycle} (Definition \ref{defi : twisting cocycle for trees}). We use this cocycle to define the $\Ai$-algebra structure in Section \ref{sec:main thm}. In the present section, we prove that the binary operation associated with this $\Ai$-algebra structure induces, in homology, the Chas-Sullivan product.

\subsection{The tree cocycle}\label{subsection : tree cocycle}

For $l,m \in \Z_{\geq 1}$ and $i \leq l$, we define $$  (\mathcal{P}X)^l\times_i  (\mathcal{P}X)^m\coloneqq \left\{ (\tau_1, \dots, \tau_l, \gamma_1, \dots, \gamma_m) \in (\mathcal{P}X)^{m+l}, \ \ev_1(\gamma_{1})= \dots = \ev_1(\gamma_m) = \ev_0(\tau_i) \right\}.$$

\begin{defi}\label{def:concat of path}
    Let $k \geq 1$, $l+m = k+1$ and $i \in \{1, \dots, l\}$. We define the (order-reversed) concatenation $$\circ_i\colon (\mathcal{P}X)^l\times_i (\mathcal{P}X)^m \to (\mathcal{P}X)^k$$ by $$(\tau_1, \dots, \tau_l)\circ_i(\gamma_1, \dots, \gamma_m) = (\tau_1, \dots, \tau_{i-1},\gamma_1\cdot\tau_i, \dots, \gamma_m \cdot\tau_i, \tau_{i+1}, \dots, \tau_l).$$

    We still denote $\circ_i \colon C_*((\mathcal{P}X)^l \times_i (\mathcal{P}X)^m) \to C_*((\mathcal{P} X)^k)$ the maps induced at the chain-level. 
\end{defi}

\begin{lemme}\label{lemme : evaluation maps trees}

There exists a family of evaluation maps $$\Gamma^k_{\ux,y} \colon \Tree_k(\ux;y) \to \mathcal{P}_{x_1\to y}X \times \dots \times \mathcal{P}_{x_k \to y}X$$ for $k \geq 1$, $(\ux,y) \in \Crit(f)^{k+1}$ which is compatible with gluing, that is, it holds
    $$\Gamma^k_{\ux,y}(T\circ_iT') = \Gamma^l_{(x_1, \dots, x_{i-1}, c, x_{i+m}, \dots, x_k),y}(T) \circ_{i} \Gamma^m_{(x_i, \dots, x_{i+m-1}),c}(T')$$
 for  $T\circ_i T' \in    \Tree_l((x_1, \dots, x_{i-1}, c, x_{i+m}, \dots, x_k);y)\circ_i \Tree_m((x_i, \dots, x_{i+m-1});c) \subset \partial\Tree_{k}(\ux;y) $.
\end{lemme}

 \begin{proof}
     Let $T \in \Tree_{k}(\ux;y)$. 
     Let $i \in \{1, \dots, k\}$. There exists a unique finite sequence of edges $$(g^i_0 = e_i, g^i_1, \dots,g^i_{\ell}=e_0)$$ given by the only non-self-crossing path connecting the incoming edge $e_i$ to the outgoing edge $e_0.$
    Define $$\lambda_i = \bigcup_{j = 0}^{\ell} \textup{Im}(\gamma_{g^i_j}).$$
    Parametrize $\lambda_i$ using the values of $f$ where $\X$ vanishes and by arc length elsewhere. This defines a parametrisation map $$\Gamma_{x_i,y} \colon \Tree_{k}(\ux;y) \to \mathcal{P}_{x_i \to y} X$$ for every $i \in \{1, \dots,k\}.$

     The evaluation maps $\Gamma^k_{\ux,y} \colon \Tree_{k}(\ux;y) \to \mathcal{P}_{x_1\to y}X \times \dots \times \mathcal{P}_{x_k \to y}X$ are then defined by $$\Gamma^k_{\ux,y}(T) = (\Gamma_{x_1,y}, \dots,\Gamma_{x_k,y}).$$
    The parametrisation maps $\Gamma_{x_i,y}$ only depend on local information. This shows that they satisfy the gluing condition.

 \end{proof}

Using the same arguments as \cite[Proposition 5.2.6]{BDHO23},  we can choose a \textit{representing chain system} for the moduli spaces of perturbed gradient trees.

\begin{lemme}\label{lemme : representing chain system for trees}
There exists a family of chains $$\left\{t_{\ux,y} \in C_{k-2+ \|y\| - \|\ux\|}\left(\Tree_k(\ux;y) \right), k \geq 1,  (\ux,y) \in \Crit(f)^{k+1}\right\}$$ such that 
\begin{enumerate}
    \item For each $(\ux,y) \in \Crit(f)^{k+1}$, the chain $t_{\ux,y}$ is a cycle relative to the boundary and represents the fundamental class $\left[\Tree_k(\ux;y), \partial \Tree_k(\ux;y)\right] \in H_{k-2 + \|y\| -\|\ux\|} \left(\Tree_k(\ux;y), \partial \Tree_k(\ux;y)\right)$.
    \item For each  $(\ux,y) \in \Crit(f)^{k+1}$, \begin{equation}\label{eq: boundary repr system of trees}
    \partial t_{\ux,y} =  \sum_{\substack{l+m = k+1\\ l,m \geq 1}} \sum_{\substack{c \\ 1 \leq i \leq l}} \pm   t_{(x_1, \dots, x_{i-1},c,x_{i+m}, \dots, x_k),y} \circ_i t_{(x_i, \dots, x_{i+m-1}),c}.
    \end{equation}
\end{enumerate}
\end{lemme}
\begin{flushright}
    $\blacksquare$
\end{flushright}

The signs are given by the difference in signs between the orientations of the boundary components of $\Tree_k(\ux;y)$ as boundary and as a product (see \eqref{eq : boundary of tree}). These are therefore dependent on the choice of orientation of $\Tree_k(\ux;y)$ and are dealt with in Appendix \ref{sec:or trees}.

\begin{defi}\label{defi : twisting cocycle for trees}
    Let $k \geq 1$ and $(\ux,y) \in \Crit(f)^{k+1}$. We define the \textit{tree cocycle} $$m^k_{\ux,y} = \Gamma^k_{\ux,y,*}(t_{\ux,y}) \in C_{k-2 + \|y\| - \|\ux\|}(\mathcal{P}_{x_1\to y}X \times \dots \times \mathcal{P}_{x_k \to y}X).$$  This family satisfies the equation 
    \begin{equation*}
    \partial m^k_{\ux,y} =  \sum_{\substack{l+h = k+1\\ l,h \geq 1}} \sum_{\substack{c \\ 1 \leq i \leq l}} \pm m^l_{(x_1, \dots, x_{i-1},c,x_{i+h}, \dots, x_k),y} \circ_{i} m^h_{(x_i, \dots, x_{i+h-1}),c} .
    \end{equation*}
    The signs in this equation are the same as in \eqref{eq: boundary repr system of trees}.
\end{defi}

\subsection{The binary operation}\label{subsec: binary operation}

Consider a set of DG Morse data $\Xi$ with Morse function $f$ and $\Delta\colon X \to X \times X$, $\Delta(a) = (a, a)$.
Let $F \hookrightarrow E \overset{\pi}{\to} X$ be a fibration endowed with a morphism of fibrations $\mu\colon E \ftimes{\pi}{\pi} E \to E$. In this subsection, we focus on the binary operation given by trees with two incoming edges. 

The moduli space of binary Morse trees can be identified with the moduli spaces that compute $\Delta_!$. The second author identified in \cite{Rie24} the composition of a Künneth map, $\Delta_!$ and $\mu_*$ as a generalisation of the Chas-Sullivan product as defined in \cite{GS07}. It is thus the following proposition that connects the operations that we construct to existing results.

We prove in Theorem \ref{thm : Functoriality} that $\mu$ induces a chain map $\widetilde{\mu}\colon C_*(X, \Delta^*C_*(F \times F)) \to C_*(X,C_*(F))$ that corresponds, via the Fibration Theorem (see Theorem \ref{thm : fibration theorem for our complex}), to the map $\mu_*\colon C_*(E \ftimes{\pi}{\pi} E) \to C_*(E)$ induced by $\mu$. This reproves the results from \cite{Rie24} in our conventions for Morse homology with DG coefficients.

The map $\Phi_2\colon (E \times E) \ftimes{(\pi,\pi)}{\ev_0} \mathcal{P}(X \times X) \to E \times E$, $\Phi_2((e,e'),(\gamma,\tau)) =(\Phi(e,\gamma), \Phi(e',\tau)) $ is a transitive lifting function associated to the product fibration $F\times F \hookrightarrow E \times E \overset{(\pi,\pi)}{\to} X \times X$.

\begin{prop}\label{prop : binary operation}

The map defined on the generators by $$\deffct{\CS^2}{C_{n+*}(X,C_*(F))^{\otimes 2}}{C_{n+*}(X ,C_*(F))}{(x_1 \otimes \alpha_1) \otimes (x_2 \otimes \alpha_2)}{\sum_{y \in \Crit(f)} (-1)^{|\alpha_1||x_2|} \widetilde{\mu}\left( y \otimes \Phi^{op}_{2,*}\left( m^2_{(x_1,x_2),y} \otimes (\alpha_1 \times  \alpha_2)\right)\right)}$$ is a chain map that induces the associative product $$\CS\colon H_{n+*}(X ,C_*(F))^{\otimes 2} \to H_{n+*}(X ,C_*(F))$$ of \cite[Theorem 7.1]{Rie24}.
\end{prop}

\begin{proof}

    The fact that $\CS^2$ is a chain map is a corollary of Theorem \ref{thm:main thm}. 

    We thus need to prove that $\CS^2$ induces the Chas-Sullivan product $$\CS = \widetilde{\mu} \circ \Delta_! \circ K\colon H_{n+*}(X ,C_*(F))^{\otimes 2} \to H_{n+*}(X ,C_*(F))$$ where $K$ denotes the Künneth map. We give a precise definition in Subsection \ref{appendix subsec: CS product}, where we prove that it coincides with the Chas-Sullivan product defined by the second author in \cite{Rie24}.

    By definition, \begin{align*}
        \CS^2 =  \widetilde{\mu} \circ T_{!} \circ K\colon H_{n+*}(X ,C_*(F))^{\otimes 2} \to H_{n+*}(X ,C_*(F)),
    \end{align*}
   where $$\deffct{T_!}{C_*(X \times X, m^K ,C_*(F \times F))}{C_*(X, \Delta^*C_*(F \times F))}{(x_1,x_2) \otimes \beta}{ \sum_y y \otimes \Phi_{2,*}^{op}( m^2_{(x_1,x_2),y} \otimes \beta)}.$$

    Therefore, it suffices to compare $$T_!\colon H_*(X \times X ,C_*(F \times F)) \to H_{-n+*}(X, \Delta^*C_*(F \times F))$$ and $$\Delta_!\colon H_*(X \times X ,C_*(F \times F)) \to H_{-n+*}(X, \Delta^*C_*(F \times F)).$$

    Both operations act via parametrisations of moduli spaces. In the rest of the proof, we identify the respective moduli spaces and compare the relevant parametrisations. While they are not the same parametrisations, they are homotopic. We show that being homotopic is enough to induce the same map on homology.

We prove that $\Delta_! = (\Phi^{-1}_{\X} \circ \Delta)_!  = T_!$, where $\Phi_{\X} \coloneqq \Phi_{T,\X}$ is the diffeomorphism isotopic to the identity defined in Definition \ref{defi : Phi_X} for $T$ the tripod, which is the only element in $\mathcal{T}_2$.

Since $\Phi^{-1}_{\X}$ is isotopic to the identity, the equality $\Delta_! = (\Phi^{-1}_{\X} \circ \Delta)_!$ is a direct consequence of the homotopy property of shriek maps \cite[Theorem 8.1.1]{BDHO23}. See also \cite[Proposition 4.4.1]{PPDG} for a construction of a chain homotopy between the shriek maps of two homotopic maps using morphisms of fibrations.

We now remark that \begin{align*}
    &\mathcal{T}^{\X}_2(x_1,x_2;y) = W^u(x_1) \times W^u(x_2) \cap \Phi^{-1}_{\X}\left( \Delta W^s(y)\right)\\
    &\cong \Delta^{-1}\Phi_{\X}\left( W^u(x_1) \times W^u(x_2) \right) \cap  W^s(y) = \mathcal{M}^{(\Phi^{-1}_{\X} \circ \Delta)_!}(x_1,x_2;y)
\end{align*} through the diffeomorphism $$\Phi_{\X}^{-1}\circ \Delta\colon \mathcal{M}^{(\Phi^{-1}_{\X} \circ \Delta)_!}(x_1,x_2;y) \to \mathcal{T}^{\X}_2(x_1,x_2;y),$$ where $\mathcal{M}^{(\Phi^{-1}_{\X} \circ \Delta)_!}(x_1,x_2;y)$, for $x_1,x_2,y \in \Crit(f)$, are the moduli spaces used to define the shriek map $(\Phi^{-1}_{\X} \circ \Delta)_!$. Since the tripod $T$ does not have internal edges, it follows that this identification extends to the compactifications $$\Phi^{-1}_{\X} \circ \Delta\colon \overline{\mathcal{M}}^{(\Phi^{-1}_{\X} \circ \Delta)_!}(x_1,x_2;y) \to \Tree_2(x_1,x_2;y).$$

We now compare the evaluation maps. The evaluation maps $\Gamma^2_{(x_1,x_2),y} \colon \Tree_2(x_1,x_2;y) \to \mathcal{P}_{x_1\to y} X \times \mathcal{P}_{x_2 \to y} X$ are defined in Lemma \ref{lemme : evaluation maps trees}. Let $(\lambda_{a_1},\lambda_{a_2},\lambda_b) \in \Tree_2(x_1,x_2;y).$ We see $\lambda_{a_1} \in \wb^u(x_1)$ as a half-infinite trajectory of $f$ going from $x_1$ to $a_1 \in X$ and $\lambda_{a_2} \in \wb^u(x_2)$ as a half-infinite trajectory of $f$ going from $x_2$ to $a_2 \in X$. Assuming $(\lambda_{a_1},\lambda_{a_2},\lambda_b) \in \Tree_2(x_1,x_2;y)$ implies that $\Phi_\X(a_1,a_2)=(b,b)$ for some $b \in \wb^s(y)$ and that $\lambda_b$ is a half-infinite trajectory going from $b$  to $y$.

Then,  $$\Gamma^2_{(x_1,x_2),y} (\lambda_{a_1},\lambda_{a_2},\lambda_b) = \left( \Gamma_{x_1,y}(\lambda_{a_1},\lambda_b), \Gamma_{x_2,y}(\lambda_{a_2},\lambda_b) \right),$$

where
$$ \Gamma_{x_1,y}(\lambda_{a_1},\lambda_b) = \Gamma_{x_1}(\lambda_{a_1}) \cdot (t \mapsto \phi^t_1(a_1)) \cdot (t \mapsto \phi_0^t(\phi^1_1(a_1)))\cdot\Gamma_{y}(\lambda_b)$$ and  $$\Gamma_{x_2,y}(\lambda_{a_2},\lambda_b)=\Gamma_{x_2}(\lambda_{a_2}) \cdot (t \mapsto \phi^t_2(a_2)) \cdot (t \mapsto \phi_0^t(\phi^1_2(a_2)))\cdot\Gamma_{y}(\lambda_b)$$ are the parametrisations by $f$ where the perturbation data vanish and by arc length otherwise.

The evaluation maps $\Gamma^{(\Phi_\X^{-1} \circ \Delta)_!}_{(x_1,x_2),y} \colon \overline{\mathcal{M}}^{(\Phi_\X^{-1} \circ \Delta)_!}(x_1,x_2;y) \to \mathcal{P}_{x_1 \to y} X \times \mathcal{P}_{x_2 \to y} X $ are defined in \S \ref{subsec:Shriek map} as follows. Let $(\lambda_{\Phi_{\X}^{-1}(b,b)}, \lambda_b) \in \overline{\mathcal{M}}^{(\Phi_\X^{-1} \circ \Delta)_!}(x_1,x_2;y)$. We see $\lambda_b \in \wb^s(y)$ as a half-infinite trajectory going from $b \in X$ to $y$ and $\lambda_{\Phi_{\X}^{-1}(b,b)}=(\lambda_{a_1},\lambda_{a_2}) \in \wb^u(x_1,x_2) = \wb^u(x_1)  \times \wb^u(x_2)$ as a half-infinite trajectory in $X\times X$ going from $(x_1,x_2)$ to $(a_1,a_2) \coloneqq \Phi_{\X}^{-1}(b,b)$.

Then, $$\Gamma^{(\Phi_\X^{-1} \circ \Delta)_!}_{(x_1,x_2),y}(\lambda_{a_1},\lambda_{a_2}, \lambda_b) = \Gamma_{x_1,x_2}^K(\lambda_{a_1,a_2}) \cdot \Phi_\X^{-1} \circ \Delta\circ \Gamma_y(\lambda_b),$$ where $\Gamma^K_{x_1,x_2} = \Gamma_{x_1} \times \Gamma_{x_2}\colon \wb^u(x_1,x_2)\to \mathcal{P}_{x_1 \to X} X \times \mathcal{P}_{x_2 \to X} X$ is the parametrisation map on a Cartesian product that parametrizes both trajectories by the values of $f$ independently. This is the evaluation map used to define the Künneth twisting cocycle $m^K$, see \cite[\S 6]{Rie24}.

It follows that
$$\Gamma^{(\Phi_\X^{-1} \circ \Delta)_!}_{(x_1,x_2),y}(\lambda_{a_1},\lambda_{a_2}, \lambda_b)= \left(\Gamma_{x_1}(\lambda_{a_1})\cdot \phi^{-1}_1\phi^{-1}_0\circ \Gamma_y(\lambda_b), \Gamma_{x_2}(\lambda_{a_2})\cdot \phi^{-1}_1\phi^{-1}_0\circ \Gamma_y(\lambda_b) \right).$$

To conclude the proof, we follow the idea of \cite[Lemma 6.3.7]{BDHO23} and build a homotopy of evaluation maps $$\Gamma^{(s)}_{(x_1,x_2),y} \colon \Tree_2(x_1,x_2;y) \to \mathcal{P}_{x_1 \to y} X \times \mathcal{P}_{x_2 \to y} X$$ between $\Gamma^2_{(x_1,x_2),y}$ and $\Gamma^{(\Phi_\X^{-1} \circ \Delta)_!}_{(x_1,x_2),y}$ that satisfies the boundary relations

\begin{itemize}
    \item If $(T,\lambda) \in \Tree_2(x'_1,x_2;y)\times \trajb{x_1,x'_1} \subset \partial \Tree_2(x_1,x_2;y)$, $$\Gamma^{(s)}_{(x_1,x_2),y}(T,\lambda) =  \Gamma^{(s)}_{(x'_1,x_2),y}(T)\circ_1 \Gamma_{x_1,x'_1}(\lambda).$$
    \item If $(T,\lambda) \in  \Tree_2(x_1,x'_2;y)\times \trajb{x_2,x'_2} \subset \partial \Tree_2(x_1,x_2;y)$, $$\Gamma^{(s)}_{(x_1,x_2),y}(T,\lambda) = \Gamma^{(s)}_{(x_1,x'_2),y}(T)\circ_2 \Gamma_{x_2,x'_2}(\lambda).$$
    \item If $(\lambda,T) \in \trajb{x_1,x'_1} \times   \Tree_2(x_1,x_2;z) \subset \partial \Tree_2(x_1,x_2;y)$, $$\Gamma^{(s)}_{(x_1,x_2),y}(\lambda,T) =  \Gamma_{z,y}(\lambda)\circ_1 \Gamma^{(s)}_{(x_1,x_2),z}(T).$$
\end{itemize}

Such a homotopy induces a chain homotopy $\Psi_!\colon C_{*}(X \times X, m^K ,C_*(F \times F)) \to C_{1+*}(X , \Delta^*C_*(F \times F))$ between $(\Phi^{-1}_\X \circ \Delta)_!$ and $T_!$. That is to say $$\partial\Psi_! + \Psi_!\partial = (\Phi^{-1}_\X \circ \Delta)_! - T_!.$$

This then proves that $\Delta_! = (\Phi^{-1}_\X \circ \Delta)_! = T_!\colon H_{*}(X \times X ,C_*(F \times F)) \to H_{*}(X, \Delta^*C_*(F \times F)).$

For $s \in [0,1]$ and $i \in \{1,2\}$, denote $$\psi^s_i= \left\{\begin{array}{cc}
    \phi^{2s}_i & \textup{if} \ s \leq 1/2   \\
     \phi^{2s-1}_0 \circ \phi^1_i& \textup{if} \  s\geq 1/2 
\end{array}. \right.$$

Let $s \in [0,1]$. Define $$\Gamma^{(s)}_{(x_1,x_2),y} (\lambda_{a_1},\lambda_{a_2},\lambda_b) = \left( \Gamma^{(s)}_{x_1,y}(\lambda_{a_1},\lambda_b), \Gamma^{(s)}_{x_2,y}(\lambda_{a_2},\lambda_b) \right),$$

where
$$ \Gamma^{(s)}_{x_1,y}(\lambda_{a_1},\lambda_b) = \Gamma_{x_1}(\lambda_{a_1}) \cdot (t \mapsto \psi_1^{ts}(a_1))\cdot\psi_1^{s}\circ \phi^{-1}_1 \circ \phi^{-1}_0\circ\Gamma_{y}(\lambda_b)$$ and  $$\Gamma^{(s)}_{x_2,y}(\lambda_{a_2},\lambda_b)=\Gamma_{x_2}(\lambda_{a_2}) \cdot (t \mapsto \psi_2^{ts}(a_2))\cdot\psi_2^{s}\circ \phi^{-1}_2 \circ \phi^{-1}_0\circ\Gamma_{y}(\lambda_b).$$

It is clear from the definition of $\Gamma^{(s)}_{(x_1,x_2),y}$ that $\Gamma^{(0)}_{(x_1,x_2),y} = \Gamma^{(\Phi^{-1}_\X \circ \Delta)_!}_{(x_1,x_2),y}$, that $\Gamma^{(1)}_{(x_1,x_2),y} = \Gamma^2_{(x_1,x_2),y}$ and that it satisfies the boundary relations.

\end{proof}

\section{Coloured operads}\label{sec:coloured ops}

In this section, we recall the facts from coloured operad theory that are relevant for this article. In the following, we sketch the relevant constructions and fix our notation in the process.

\subsection{Notation}

With the exception of Section \ref{sec:top op}, we work with operads in chain complexes. We thus introduce the notation in this setting. However, everything in this subsection is also true for coloured topological operads after replacing the symbol $\otimes$ with $\times$.
\begin{notation}
    Let $C$ denote a set.
    \begin{itemize}
        \item  For a $C$-coloured operad $\mathcal{O}$ and $(\ux,y)\in C^{k+1}$, we denote $\mathcal{O}(\ux;y)$ the space of operations with input colours $\ux$ and output colours $y$;
        \item for tuples $\ux=(x_1,\dots, x_l)\in C^l$ and $\ux'=(x_1',\dots, x_m')\in C^m$, we denote 
        \begin{align*}
            \ux\circ_i \ux'=(x_1,\dots, x_{i-1},x_1',\dots,x_m',x_{i},\dots, x_l)\in C^{l+m-1};
        \end{align*}
        \item we denote the partial composition by
        \begin{align*}
            \circ_i\colon \mathcal{O}(\ux;y)\otimes \mathcal{O}(\ux';x_{i})\to \mathcal{O}(\ux\circ_i\ux';y).
        \end{align*}
    \end{itemize}

    Together with a permutation action on the spaces $\mathcal{O}(\ux;y)$, this defines a $C$-coloured operad if the composition and permutation action satisfy certain associativity relations (see e.g.\ \cite[\S 16.2]{yau2016colored}).
\end{notation}

The partial compositions are chain maps, i.e., they satisfy the Leibniz rule $\partial(\mu\circ_i\nu)=\partial(\mu)\circ_i \nu+(-1)^{|\mu|}\mu\circ_i \partial(\nu)$.

\begin{ex}
    For a collection of chain complexes $(A_x)_{x\in C}$, the \textit{coloured endomorphism operad} $\End_{(A_x)_{x\in C}}$ is given by 
    \begin{align*}
        \End_{(A_x)_{x\in C}}(x_1,\dots,x_k; x)=\Hom(A_{x_1}\otimes \dots \otimes A_{x_k},A_x)
    \end{align*}
    and composition is given by composition of morphisms.
\end{ex}

\begin{defi}
    Let $\mathcal{O}$ be a $C$-coloured operad. An \textit{$\mathcal{O}$-algebra structure} on a collection of chain complexes $(A_x)_{x\in C}$ is given by a morphism of $C$-coloured operads $\mathcal{O}\to \End_{(A_x)_{x\in C}}$.
\end{defi}

\subsection{Signs and suspension}\label{subsec:signs and susp}
 
Working in a graded context, there are some signs in formulas that change depending on certain choices. In this subsection, we fix our conventions.

We work with the convention that maps act from the left, that is, for chain complexes $A$ and $B$, the canonical evaluation map is given by
\begin{align*}
    \Hom(A,B)\otimes A&\to B,\\
    f\otimes a&\mapsto  f(a).
\end{align*}

\begin{notation}
    We denote the unit of the symmetric monoidal category of chain complexes by $\mathbbm{1}$. It is given by the graded abelian group
    \begin{align*}
        \mathbbm{1}_*=\begin{cases}
            \Z &\text{if } *=0;\\
            0&\text{else;}
        \end{cases}
    \end{align*}
    with a fixed generator $1$ in degree zero.
\end{notation}

\begin{rem}\label{rem:natural decomp of Hom}
    As a consequence of our convention, the natural isomorphism for decomposition of $\Hom$-spaces is the following
    \begin{align*}
        \Hom(A,B)\cong B\otimes \Hom(A,\mathbbm{1}).
    \end{align*}
    Indeed, this decomposition is natural. This can be seen by noting commutativity of the following diagram
    \begin{center}
        \begin{tikzcd}
            \Hom(A,B)\otimes A \ar[r, "\cong"] \ar[d, "\mathrm{ev}"] & B\otimes \Hom(A,\mathbbm{1})\otimes A\ar[d, "\mathrm{ev}"]\\
            B\ar[r, "\cong"] & B\otimes \mathbbm{1}.
        \end{tikzcd}
    \end{center}
\end{rem}

\begin{defi}\label{def:susp op}
    The \textit{suspension chain complex} $\Sigma$ is given by the graded abelian group
    \begin{align*}
        \Sigma_*\coloneqq\begin{cases}
            \Z &\text{if } *=1;\\
            0&\text{else;}
        \end{cases}
    \end{align*}
    with a fixed generator $s$ in degree one.    
\end{defi}

\begin{notation}
    For $k\geq 0$, we denote $\Sigma^k\coloneqq\Sigma^{\otimes k}$ with fixed generator $s^k\coloneqq s^{\otimes k}$ in degree $k$. For $k<0$, we denote $\Sigma^k\coloneqq\Hom(\Sigma^{|k|},\mathbbm{1})$ with fixed generator $s^k\coloneqq(s^{|k|}\mapsto 1)$.
\end{notation}

\begin{defi}
    For a chain complex $A=(A_*,d_A)$ and $k\in \Z$, the \textit{$k$-th suspension of $A$} is given by $A_{-k+*}\coloneqq\Sigma^kA\coloneqq\Sigma^k\otimes A$ with differential $\partial_{\Sigma^kA}\coloneqq\id \otimes d_A$. For $a\in A$, we denote $s^ka\coloneqq s^k\otimes a\in\Sigma^k\otimes A$.
\end{defi}

\subsubsection{Suspension operad}\label{subsubsec:susp operad}
We now explain how to encode degree shifts on the level of operads.
\begin{defi}\label{def:coloured susp op}
    Let $g\colon C\to \Z$ denote a map. The \textit{$g$-graded suspension operad} $\underline{\Sigma}^{g}$ is given by the endomorphism operad of the collection $(\underline{\Sigma}^{g(x)})_{x\in C}$.
\end{defi}

\begin{lemme}\label{lem:suspending operads}
    Let $C$ be a set with a grading $g\colon C\to \Z$. The evaluation maps\footnote{This map together with the coevaluation map \eqref{eq:coev maps} exhibit $\Sigma^{-g(x)}$ as the dual object of $\Sigma^{g(x)}$. Hence our naming convention.}
    \begin{align}\label{eq:map inducing equivalence}
        \begin{split}
            \underline{\Sigma}^{-g(x)}\otimes \underline{\Sigma}^{g(x)}&\to \mathbbm{1},\\
            s^{-g(x)}\otimes s^{g(x)}&\mapsto 1
        \end{split}
    \end{align}
    exhibit tensoring with $\underline{\Sigma}^{g}$ as an automorphism of the category of $C$-coloured operads.
\end{lemme}
\begin{proof}
    The grading $-g\colon C\to \Z$ defines a suspension operad $\underline{\Sigma}^{-g}$. The collection of maps \eqref{eq:map inducing equivalence} defines a map of operads $\underline{\Sigma}^{-g}\otimes \underline{\Sigma}^{g}\to \underline{\Sigma}^0$. Tensoring with $\underline{\Sigma}^0$ is canonically identified with the identity functor on $C$-coloured operads. We thus obtain a counit transformation $\underline{\Sigma}^{-g}\otimes \underline{\Sigma}^{g}\otimes -\to\id $. The corresponding unit transformation $\id\to \underline{\Sigma}^{-g}\otimes \underline{\Sigma}^{g}\otimes -$ is induced by the maps
    \begin{align*}
        \mathbbm{1}&\to \Sigma^{-g(x)}\otimes \Sigma^{g(x)},\\
        1&\mapsto s^{-g(x)}\otimes s^{g(x)}.
    \end{align*}
    This exhibits tensoring with $\underline{\Sigma}^{-g}$ as a left adjoint to tensoring with $\underline{\Sigma}^g$. Moreover, the above maps are equivalences and thus tensoring with $\underline{\Sigma}^{-g}$ defines an inverse of $\underline{\Sigma}^g$. 
\end{proof}

\begin{rem}\label{rem:ev and coev}
    The above statement is not symmetric under changing the roles of $g$ and $-g$. In the statement with reversed roles, the coevaluation maps
    \begin{align}\label{eq:coev maps}
        \begin{split}
            \Sigma^{g(x)}\otimes \Sigma^{-g(x)}&\to \mathbbm{1},\\
            s^{g(x)}\otimes s^{-g(x)}&\mapsto 1
        \end{split}
    \end{align}
    induce the equivalence. The evaluation maps \eqref{eq:map inducing equivalence} and the coevaluation maps \eqref{eq:coev maps} composed with the Koszul twist differ by the sign $(-1)^{g(x)}$.

    The upshot is that, when identifying composition of shifts, there is a choice. Throughout this article, we fix the convention that we use the evaluation maps to induce the isomorphism.
\end{rem}

\begin{ex}\label{exmp:suspending algebras}
    For a grading $g\colon C\to \Z$, a $C$-coloured operad $\mathcal{O}$ and a $C$-collection $(A_x)_{x\in C}$, it is equivalent to give an $\mathcal{O}$-algebra structure on $(A_x)_{x\in C}$ and to give a $\underline{\Sigma}^g\otimes \mathcal{O}$-algebra structure on $(\underline{\Sigma}^{g(x)} A_x)_{x\in C}$: an $\mathcal{O}$-algebra structure on $(A_x)_{x\in C}$ is given by a map of operads $\mathcal{O}\to \End_{(A_x)_{x\in C}}$. After tensoring with the suspension operad $\underline{\Sigma}^g$, we obtain a map of operads
    \begin{align*}
        \underline{\Sigma}^g\otimes \mathcal{O}\to \underline{\Sigma}^g\otimes \End_{(A_x)_{x\in C}}=\End_{(\Sigma^{g(x)})_{x\in C}}\otimes \End_{(A_x)_{x\in C}}\cong \End_{(\Sigma^{g(x)}A_x)_{x\in C}},
    \end{align*}
    i.e.\ a $\underline{\Sigma}^{g}\otimes \mathcal{O}$-algebra structure on $(\Sigma^{g(x)}A_x)_{x\in C}$. By Lemma \ref{lem:suspending operads}, this assigment is an equivalence.
\end{ex}

\subsection{Decolouring}\label{subsec:dec}
In this subsection, we review a canonical way to produce a one-coloured operad out of a $C$-coloured operad. 

\begin{defi}\label{def:decolouring}
    Let $\mathcal{O}$ be a $C$-coloured operad. The \textit{decolouring} of $\mathcal{O}$ is the one-coloured operad given by 
    \begin{align*}
        \dec(\mathcal{O})(k)\coloneqq\bigoplus_{\substack{\ux\in C^k,\\ y\in C}}\mathcal{O}(\ux;y) 
    \end{align*}
    and composition is given by the extension by zero, i.e., by the following composition
    \begin{align*}
        \dec(\mathcal{O})(k)\otimes  \dec(\mathcal{O})(l)&\cong \bigoplus_{\substack{\ux\in C^k y\in C,\\
        \ux'\in C^ly'\in C}}\mathcal{O}(\ux;y)\otimes \mathcal{O}(\ux';y')\\ 
        &\to\bigoplus_{\substack{\ux\in C^k y\in C,\\
        \ux'\in C^l}}\mathcal{O}(\ux;y)\otimes \mathcal{O}(\ux';x_i)\\
        &\overset{\circ_i}{\to}\bigoplus_{\substack{\ux\in C^k y\in C,\\
        \ux'\in C^l}}\mathcal{O}(\ux\circ_i\ux';y)\to \dec(\mathcal{O})(k+l-1).
    \end{align*}
\end{defi}

\begin{rem}
    The decolouring of a $C$-coloured operad can be extended to a functor from $C$-coloured operads to one-coloured operads.
\end{rem}

\begin{ex}\label{ex:decolouring algebras}
    Let $C$ be a finite set and $(A_x)_{x\in C}$ be a collection of chain complexes. The decolouring of $\End_{(A_x)_{x\in C}}$ is isomorphic to the endomorphism operad of $\bigoplus_{x\in C} A_x$. The isomorphism is given by the natural isomorphism:
    \begin{align*}
        \bigoplus_{(x_1,\dots,x_k,x)\in C^{k+1}}\Hom(A_{x_1}\otimes \dots \otimes A_{x_k},A_x)\cong \Hom\left(\bigoplus_{x_1\in C^k}A_{x_1}\otimes \dots \otimes\bigoplus_{x_k\in C} A_{x_k},\bigoplus_{x\in C}A_x\right), 
    \end{align*}
    which exists because $C$ is finite.
\end{ex}

\section{\texorpdfstring{$\Ai$}{Ainfty}-operads}\label{sec:Coloured Ainfty}
In this section, we define a coloured version of the $\Ai$-operad. It encodes $\Ai$-algebras where the underlying graded abelian group has a decomposition. First, we recall various characterisations of the one-coloured $\Ai$-operad.

\subsection{One-coloured \texorpdfstring{$\Ai$}{Ainfty}-operads}
In the literature, there are various operads that are called the $\Ai$-operad. In this subsection, we recall a few of the definitions and give them names to differentiate them. 

Stasheff's original definition is given in terms of the associahedron \cite{Sta70}. A related construction is via the moduli space of metric trees $\mathcal{T}_k$ with $k$ incoming leaves (see Subsection \ref{subsec:stable metric ribbon trees}). This moduli space has a canonical cell structure where the interior of the top cell is given by all trees with finite length interior edges and the boundary is given by all trees with at least one interior edge of infinite length. The moduli spaces of trees and concatenation make this into an operad whose composition maps are cellular. We refer to \cite[\S I.4]{boardman2006homotopy} for a more detailed explanation and an argument why this CW space is equivalent to Stasheff's definition of the associahedron.

This explains the following definition of a version of the $\Ai$-operad.

\begin{defi}
    The \textit{cellular $\Ai$-operad} is given by $\Ai^{cell}(k)\coloneqq C_*^{cell}(\mathcal{T}_k)$ and composition is given by gluing of metric trees.
\end{defi}

While $\mathcal{T}_k$ has a canonical cell structure, the cells are not canonically oriented. Several orientation conventions appear in the literature. We adopt the convention of \cite[\S 9.2.4]{loday2012algebraic}.

\begin{defi}
    The \textit{oriented $\Ai$-operad}, $\Ai^{or}$, is given by the operad generated by elements $\{m_k\in (\Ai^{or})_{k-2}(k)\mid k\geq 2\}$ modulo the relations that for all $k\geq 2$, it holds
    \begin{align*}
        \partial(m_k)=\sum_{\substack{l+m=k+1 \\
        l,m \geq 2}} \sum_{i=1}^{l} (-1)^{i + m(l-i)} m_l \circ_{i}m_m.
    \end{align*}
\end{defi}

\begin{lemme}[{\cite[Proposition 9.2.4]{loday2012algebraic}}]\label{lem:or Ai to cell Ai}
    There exist orientations $[\mathcal{T}_k]\in C_{k-2}^{cell}(\mathcal{T}_k)$ of the top cell of $\mathcal{T}_k$ such that
    \begin{align*}
        \Ai^{or}&\to \Ai^{cell},\\
        m_k&\mapsto [\mathcal{T}_k]
    \end{align*}
    is an isomorphism of operads.
\end{lemme}

For our purposes, it is useful to define an operad that encodes $\Ai$-algebras where the differential $\partial$ of the underlying chain complex can be written as $\partial=\partial_1+\partial_2$ where $\partial_1$  is a differential. This situation occurs for the complex with DG coefficients, whose differential naturally splits into the fibre differential and the Morse contribution.

\begin{defi}
    The \textit{augmented oriented $\Ai$-operad}, $\Ai^+$, is the operad freely generated by elements $\{m_k\in (\Ai^+)_{k-2}(k)\mid k\geq 1\}$ subject to the relation that, for all $k\geq 1$, it holds
    \begin{align}\label{eq:Ainfty rel with extra diff}
        \partial(m_k)=\sum_{l+m=k+1}\sum_{i=1}^{l} (-1)^{i+m(l-i)} m_l\circ_{i} m_m.
    \end{align}
\end{defi}

We follow the naming convention from \cite{laplante2024universal} for this operad. It corresponds to $\hat{\mathcal{A}}_\infty$ in the notation of \cite[\S 5.2]{chuang2013combinatorics}, where it is shown that $\hat{\mathcal{A}}_\infty$ is an acyclic operad \cite[Theorem 5.7 (2)]{chuang2013combinatorics}. This operad thus does not carry any homological information. However, it allows a twisted differential which is of interest to us.

To twist the differential, we use a Maurer-Cartan element. By \eqref{eq:Ainfty rel with extra diff} for $k=1$, the element $m_1\in \Ai^+(1)_{-1}$ satisfies the associative Maurer-Cartan equation
\begin{align}\label{eq:MC}
    \partial(m_1)+m_1\circ_1m_1=0.
\end{align}
Such a Maurer-Cartan element can be used to twist the differential:

\begin{lemme}[e.g. {\cite[Proposition 5.5]{dotsenko2023maurer}}]\label{lem:MC gives twisted diff}
    Let $\mathcal{O}$ be an operad with an element $m_1\in \mathcal{O}(1)_{-1}$ that satisfies the associative Maurer-Cartan equation \eqref{eq:MC}. The map 
    \begin{align*}
        \partial_{m_1}\colon \mathcal{O}(k)&\to \mathcal{O}(k),\\
        \mu&\mapsto \partial(\mu)+m_1\circ_1\mu+(-1)^{|\mu|+1}\sum_{i=1}^k \mu \circ_{i} m_1
    \end{align*}
    defines a differential.
\end{lemme}

\begin{lemme}\label{lemme:Ainfty or to Ainfty split}
    The assignment 
    \begin{align*}
        (\Ai^{or},\partial)&\to (\Ai^+,\partial_{m_1}),\\
        m_k&\mapsto m_k
    \end{align*}
    defines a quasi-isomorphism of operads.
\end{lemme}
\begin{proof}
    Rewriting \eqref{eq:Ainfty rel with extra diff}, we obtain 
    \begin{align*}
        \partial_{m_1}(m_k)=\sum_{\substack{l+m=k+1 \\
        l,m \geq 2}} \sum_{i=1}^{l} (-1)^{i + m(l-i)} m_l \circ_{i}m_m.
    \end{align*}
    In the notation of \cite[\S 5.2]{chuang2013combinatorics}, the operad $(\Ai^+,\partial_{m_1})$ corresponds to $\Ai[m]$ and the above map gives the canonical inclusion $\Ai\to \Ai[m]$ which is a quasi-isomorphism by \cite[Theorem 5.7 (1)]{chuang2013combinatorics}.
\end{proof}

\begin{ex}\label{ex:Ai split alg struct}
    If we have an $\Ai^+$-algebra structure $\varphi \colon \Ai^+\to \End_{(A,\partial)}$ on a chain complex $(A,\partial)$, we get that $\partial_{\varphi(m_1)}\coloneqq \partial+\varphi(m_1)$ defines a twisted differential on $A$. Moreover, the collection $\{\varphi(m_k)\mid k\geq 2\}$ defines an $\Ai^{or}$-algebra structure on $(A,\partial_{\varphi(m_1)})$ (see \cite[Example 5.6]{dotsenko2023maurer}).
\end{ex}

\begin{rem}
    There is another version of the $\Ai$-operad that appears in the literature (e.g.\ \cite[\S 9.2.4]{loday2012algebraic}) but is not used in our article. For completeness, we briefly discuss how it compares to the above versions of the $\Ai$-operad.
    
    Like the augmented $\Ai$-operad, $\Ai^+$, this operad is generated by operations $m_k$ in arity $k$ and degree $k-2$ for all $k\geq 1$. But instead of requiring \eqref{eq:Ainfty rel with extra diff}, these elements satisfy
    \begin{align*}
        \sum_{l+m=k+1}\sum_{i=1}^{l} (-1)^{i+m(l-i)} m_l\circ_{i} m_m=0.
    \end{align*}
    This operad encodes $\Ai^{or}$-algebras together with the differential of the algebra.
\end{rem}

\subsection{The coloured \texorpdfstring{$\Ai$}{Ainfty}-operad}
In this subsection, we define a $C$-coloured operad whose algebras are $\Ai$-algebras with a $C$-decomposition of the underlying graded abelian group. 

\begin{notation}
    For a tuple $\ux\in C^k$, we denote $\ell(\ux)\coloneqq k$.
\end{notation}

\begin{defi}\label{defi : Ai coloured operad}
    Let $C$ be a set. The \textit{$C$-coloured $\mathcal{A}_\infty$-operad} $\mathcal{A}^C_\infty$ is generated by elements $m_{\underline{x},y}\in \mathcal{A}_\infty^C(\underline{x};y)_{\ell(\ux)-2}$ for all tuples $\underline{x}, y$, such that 
    \begin{align}\label{eq:coloured Ainfty rel}
        \partial (m_{\underline{x},y})=\sum_{\underline{x}=\underline{w}\circ_{i}\underline{z}}(-1)^{i+\ell(\underline{z})(\ell(\underline{w})-i)}m_{\underline{w},y}\circ_{i} m_{\underline{z},w_{i}}. 
    \end{align}
\end{defi}

\begin{prop}\label{prop:from one coloured Ainfty to coloured}
    Let $C$ be a finite set. There exists a morphism of operads 
    \begin{align*}
        \Ai^+&\to \dec(\Ai^C),\\
        m_k&\mapsto \sum_{\substack{(\ux,y)\in C^k\times C}} m_{\ux,y},
    \end{align*}
    where $\dec$ denotes the decolouring functor (see Definition \ref{def:decolouring}).
\end{prop}
\begin{proof}
    Let $\varphi$ denote the map in the statement. It suffices to check that $\varphi(m_k)$ satisfy \eqref{eq:Ainfty rel with extra diff}, the defining relation of $\Ai^+$, for all $k\geq 1$:
    \begin{align*}
        \partial\varphi(m_k)&= \sum_{\substack{(\ux,y)\in C^k\times C}}\partial( m_{\ux,y})\\
        &\!\overset{\eqref{eq:coloured Ainfty rel}}{=} \sum_{\substack{(\ux,y)\in C^k\times C}}\sum_{\underline{x}=\underline{w}\circ_{i}\underline{z}}(-1)^{i+\ell(\underline{z})(\ell(\underline{w})-i)}m_{\underline{w},y}\circ_{i} m_{\underline{z},w_{i}}\\
        &= \sum_{l+m=k+1}\sum_{i=1}^{l} (-1)^{i+m(l-i)} \sum_{\substack{(\uw;y)\in C^l\times C,\\
        (\uz;v)\in C^m\times C,\\
        v=w_{i}}} m_{\uw,y}\circ_{i} m_{\uz,v}\\
        &= \sum_{l+m=k+1}\sum_{i=1}^{l} (-1)^{i+m(l-i)} \sum_{\substack{(\uw;y)\in C^l\times C,\\
        (\uz;v)\in C^m\times C}} m_{\uw,y}\circ_{i} m_{\uz,v}\\
        &=\sum_{l+m=k+1}\sum_{i=1}^{l} (-1)^{i+m(l-i)}\varphi( m_l)\circ_{i} \varphi(m_m).
    \end{align*}
\end{proof}

\begin{cor}\label{cor:Ai to dec coloured Ai}
    There exists a morphism of operads
    \begin{align*}
        \phi\colon \Ai^{or}&\to (\dec(\Ai^C),\partial_{(\sum m_{x,y})}),\\
        m_k&\mapsto \sum_{\substack{(\ux,y)\in C^{k+1}}}m_{\ux,y}.
    \end{align*}
\end{cor}
\begin{proof}
    Let $\varphi$ be the map from Proposition \ref{prop:from one coloured Ainfty to coloured}. The twisted differential $\partial_{m_1}$ on $\Ai^+$ induces a twisted differential $\partial_{\varphi(m_1)}$ on $\dec(\Ai^C)$. Since $\varphi(m_1)=\sum_{x, y} m_{x,y}$, the twisted differential on the decolouring is precisely $\partial_{\sum m_{x,y}}$, proving the claim.
\end{proof}

We conclude this section by giving an example that explains how decolouring a coloured $\Ai$-algebra produces an $\Ai$-algebra.

\begin{ex}\label{ex:Ai structure via coloured Ai}
    Let $C$ be a finite set and let $(A_x,\partial)_{x\in C}$ be a $C$-collection of chain complexes. We fix an $\Ai^C$-algebra, $\varphi\colon \Ai^C\to \End_{(A_x,\partial)}$. This induces a morphism of operads $\Ai^+\to \dec(\Ai^C)\to \dec(\End_{(A_x,\partial)_{x\in C}})\cong \End_{\bigoplus_{x\in C} A_x}$ by Example \ref{ex:decolouring algebras}. We thus get an $\Ai^+$-algebra structure on $\bigoplus_{x\in C} A_x$ with the differential $\sum_{x\in C} \partial$. 

    By Example \ref{ex:Ai split alg struct}, we get an $\Ai$-algebra structure on $\bigoplus_{x\in C} A_x$ with the twisted differential given on $a\in A_x$ by
    \begin{align*}
        \partial(a)+\sum_{x, y} \varphi(m_{x,y})(a).
    \end{align*}
    This differential has the same form as the differential on Morse homology with DG coefficients. Up to shift, the goal is thus to find an $\Ai^{\Crit(f)}$-action on $(C_*(F_x))_{x\in \Crit(f)}$.
\end{ex}

\section{Topological operads}\label{sec:top op}

In this section, we construct an action of a topological operad given a fibration with a fibrewise monoidal structure. The action is constructed by first flowing according to a transitive lifting function and then applying the monoidal operation. In Subsection \ref{subsec:Path operad}, we give this action in the special case of lift-preserving monoidal fibrations (see Definition \ref{defi: monoidal morphism of fibration}). In general, this action is not strict but only homotopy coherent. In Subsection \ref{subsec:W-constr}, we recall how to encode such a homotopy coherent action and how to rectify this to a strict action. We conclude this section by giving and rectifying such a homotopy coherent action in Subsection \ref{subsec:homotopy path alg}.

\subsection{Path operads}\label{subsec:Path operad}
In this subsection, we define a path operad and an operad of gradient trees in a manifold with a Morse function. We explain how these operads act on the fibres of a monoidal fibration with a compatible transitive lifting function.

\begin{defi}\label{def:top path op}
    The \textit{topological path-operad of $X$ associated to a Morse function $f\colon X\to \R$}, denoted by $\mathcal{PO}(X,f)$, is the topologically enriched $\Crit(f)$-coloured operad given on $\underline{x}=(x_1,\dots,x_k)\in \Crit(f)^k$ and $y\in \Crit(f)$ by
    \begin{align*}
        \mathcal{PO}(X,f)(\underline{x};y)=\mathcal{P}_{x_1\to y}X\times \dots \times \mathcal{P}_{x_k\to y}X
    \end{align*}
    and composition given by the diagonal and then composition of paths (see Definition \ref{def:concat of path}).
\end{defi}

The goal is to construct an action of this operad on the fibres of a fibration over $X$. To construct this, we need a fibrewise monoidal structure. As a first step, we assume that the monoidal structure is compatible with the transitive lifting function (see Definition \ref{defi: TLF}). This assumption can be dropped using some standard operad/model category theory (see Subsection \ref{subsec:homotopy path alg}).

\begin{defi}\label{defi: monoidal morphism of fibration}
    A \textit{monoidal fibration} is a fibration $\pi\colon E\to X$ endowed with an associative fibrewise monoidal map $\mu\colon E \ftimes{\pi}{\pi} E \to E$.
    
    Given a monoidal fibration $(E,\mu)$ with a transitive lifting function $\Phi\colon E \ {}_\pi\!\times_{ev_0} \mathcal{P}X\to E$, the triple $(E,\Phi,\mu)$ is a \textit{lift-preserving monoidal fibration} if it holds
    \begin{equation}\label{eq: monoidal}
        \mu(\Phi(e,\gamma),\Phi(e',\gamma)) = \Phi(\mu(e,e'),\gamma )
    \end{equation}  
    for all $(e,e',\gamma) \in E \ftimes{\pi}{\pi} E\ftimes{\pi}{\ev_0}\mathcal{P}X$.
\end{defi}

\begin{rem}\label{rem:paths acting on fibres}
    A lift-preserving monoidal fibration $(E,\Phi,\mu)$ inherits an action of the operad $\mathcal{PO}(X,f)$ on $(F_x)_{x\in \Crit(f)}$ given by
    \begin{align*}
        \mathcal{PO}(X,f)(\ux;y)&\to \End_{(F_x)_{x\in \Crit(f)}},\\
        (\lambda_1,\dots, \lambda_k)&\mapsto [(\alpha_1,\dots,\alpha_k)\mapsto \mu(\Phi(\alpha_1,\lambda_1),\dots,\Phi(\alpha_k,\lambda_k))].
    \end{align*}

    Indeed, a direct computation shows that the associativity of $\mu$ and the fact that $\mu$ and $\Phi$ commute as in \eqref{eq: monoidal} imply that the above maps define an operad action.
\end{rem}

\begin{ex}\label{ex : assumption holds for loop loop}
    We consider the loop-loop fibration $\ev_0\colon \mathcal{L}X\to X$ together with $\mu$ the concatenation of loops. It is a lift-preserving monoidal fibration if we take $\Phi$ a transitive lifting function coming from a transitive diffeomorphism lift which were defined by the first author in \cite{Cli25}. 
\end{ex}

\begin{defi}\label{def:top gradient tree op}
    The \textit{topological gradient tree operad associated to $(X,f)$} is the topologically enriched operad $\mathcal{GTO}_{top}^{(X,f,\X)}$ with space of operations $\mathcal{GTO}_{top}^{(X,f,\X)}(\ux;y)\coloneqq\Tree_k(\ux;y)$ as in Definition \ref{def:Tree_k(x,y)} and composition given by composition of trees.
\end{defi}

\begin{rem}\label{rem:evaluating trees}
    We recall that in Lemma \ref{lemme : evaluation maps trees}, we defined evaluation maps $\Gamma_{\ux,y}^k\colon\Tree_k(\ux;y)\to \mathcal{P}_{x_1\to y}X\times \dots \times \mathcal{P}_{x_k\to y}X$ that are compatible with gluing. In the language of operads, these evaluation maps assemble to a map of operads $\Gamma\colon\mathcal{GTO}_{top}^{(X,f,\X)}\to \mathcal{PO}(X,f)$.
\end{rem}

Assembling Remarks \ref{rem:paths acting on fibres} and \ref{rem:evaluating trees}, we obtain the following construction.

\begin{cor}\label{cor:GTO on fibre lift}
    A lift-preserving monoidal fibration $(E,\Phi,\mu)$ inherits an action of the operad $\mathcal{GTO}_{top}^{(X,f,\X)}$ on the collection $(F_x)_{x\in \Crit(f)}$.
\end{cor}

\subsection{W-construction}\label{subsec:W-constr}
In this subsection, we recall the W-construction or Boardman-Vogt resolution as introduced by \cite{boardman2006homotopy}. We only sketch the construction here. For a more detailed construction for coloured operads, we refer to the article by Berger and Moerdijk \cite[\S 3]{berger2006boardman}.

\begin{defi}\label{defi:W-const}
    For a $C$-coloured topological operad $\mathcal{O}$, the \textit{W-construction} $W(\mathcal{O})$ is given by the following $C$-coloured operad. 
    
    An operation in $W(\mathcal{O})(x_1,\dots,x_k;y)$ is represented by a labelled planar tree (see Figure \ref{fig:W-constr}) where:
    \begin{itemize}
        \item all edges are labelled by a colour such that the incoming edges are labelled by $x_1,\dots,x_k$ and the outgoing edge is labelled by $y$;
        \item all vertices are labelled by an operation of $\mathcal{O}$ compatible with the edge labelling;
        \item all internal edges are labelled by a length $t_e\in [0,1]$.
    \end{itemize}

    \begin{figure}[H]
        \centering
        \includegraphics[width=0.12\linewidth]{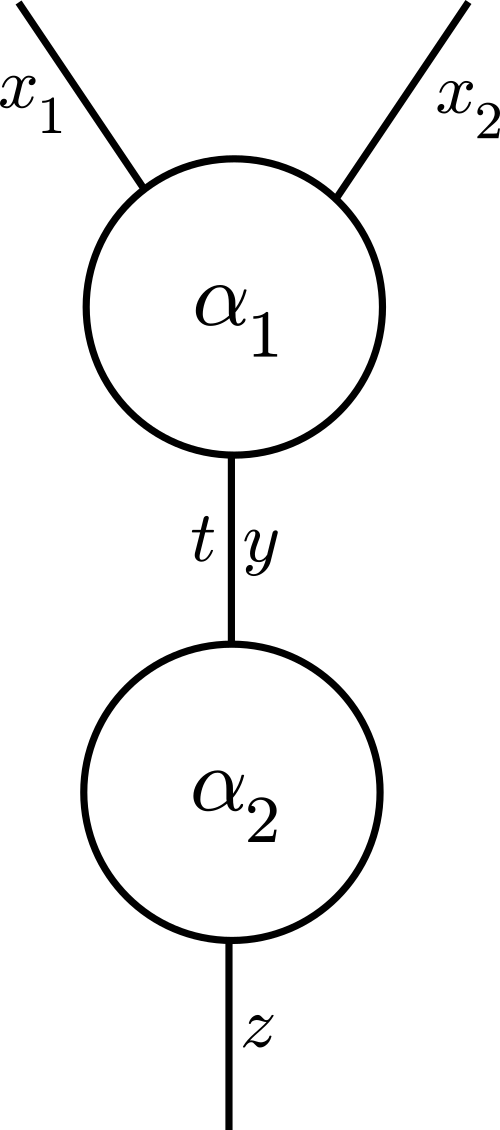}
        \caption{An element of $W(\mathcal{O})(x_1,x_2;z)$ labelled by $\alpha_1\in \mathcal{O}(x_1,x_2;y)$ and $\alpha_2\in \mathcal{O}(y;z)$ and the sole internal edge has length $t$.}
        \label{fig:W-constr}
    \end{figure}
    
    This space is subject to the relations:
    \begin{enumerate}[(i)]
        \item if the tree contains an edge $e$ of length $t_e=0$, the element is identified with the tree where $e$ is contracted and the vertex labellings are composed;

        \item if the tree contains a vertex $v$ labelled by the identity element with input edge of length $s$ and output edge of length $t$, the element is identified with the tree where the vertex is deleted and the edge is labelled by $s\star t\coloneqq 1-(1-s)(1-t)$.\footnote{We use the original \enquote{reversed multiplication} as in \cite{boardman2006homotopy} instead of the maximum operation as in \cite{berger2005resolution}. Both operations work as is explained in \cite{berger2006boardman}. The reversed multiplication simplifies our explicit definition of an action of a $W$-operad (see Definition \ref{defi:homotopy path algebra}).}
    \end{enumerate}
    Composition is given by concatenation of trees and labelling the newly created internal edge with $1$.
\end{defi}

The idea is that $W(\mathcal{O})$ encodes homotopy $\mathcal{O}$-algebras. By composing the vertex labellings, we obtain a map of operads $\epsilon \colon W(\mathcal{O})\to \mathcal{O}$. This induces a functor $\epsilon^*$ from $\mathcal{O}$-algebras to $W(\mathcal{O})$-algebras.

Under suitable conditions, the homotopy theory of $W(\mathcal{O})$-algebras is equivalent to the homotopy theory of $\mathcal{O}$-algebras. This can be made precise and explicit in broad generality (see \cite{berger2005resolution, berger2006boardman}). However in this article, we only need the following rectification lemma:

\begin{lemme}\label{lem:W-constr}
    Let $\mathcal{O}$ be a $C$-coloured operad where each operation space is a CW-complex with a free permutation action. For a $W(\mathcal{O})$-algebra structure on $(A_x)_{x\in C}$, there exists an $\mathcal{O}$-algebra structure on $(\widetilde{A}_x)_{x\in C}$ and a map of $W(\mathcal{O})$-algebras $(A_x)_{x\in C}\to \epsilon^*(\widetilde{A}_x)_{x\in C}$ such that each map $A_x\to \widetilde{A}_x$ is a weak equivalence.
\end{lemme}
\begin{proof}
    This follows from the more general statement \cite[Corollary 5.1]{berger2005resolution} which states that under certain conditions $\epsilon^*$ is the right adjoint of a Quillen equivalence. 

    We thus only need to ensure that the assumptions are satisfied: 
    \begin{itemize}
        \item $\mathcal{O}$ and $W(\mathcal{O})$ are {admissible};

        \item the underlying model category is left proper;

        \item $\mathcal{O}$ is $\Sigma$-cofibrant.
    \end{itemize}
    By \cite[Theorem 2.1]{berger2005resolution} any operad in compactly generated spaces is automatically admissible. Moreover the category of compactly generated spaces is left proper. It remains to check that $\mathcal{O}$ is $\Sigma$-cofibrant. This is ensured by the assumption that the permutation action is free on the operation spaces and that they are CW-complexes. 
\end{proof}

\subsection{Homotopy path algebras}\label{subsec:homotopy path alg}

In this subsection, we explain how to drop the assumption that the monoidal fibration is lift-preserving and still obtain an action of $\mathcal{GTO}_{top}^{(X,f,\X)}$ on the fibres $(F_x)_{x\in \Crit(f)}$. The key for that result is the W-construction of the operads $\mathcal{PO}(X,f)$ and $\mathcal{GTO}_{top}^{(X,f,\X)}$.

Our strategy is the following:
\begin{enumerate}[(1)]
    \item we define an operad action of $W(\mathcal{PO}(X,f))$ on $(F_x)_{x\in \Crit(f)}$ (see Definition \ref{defi:homotopy path algebra});
    \item by Remark \ref{rem:evaluating trees}, this induces an operad action of $W(\mathcal{GTO}_{top}^{(X,f,\X)})$ on $(F_x)_{x\in \Crit(f)}$;
    \item the operad $\mathcal{GTO}_{top}^{(X,f,\X)}$ satisfies the assumption of the rectification lemma (Lemma \ref{lem:W-constr}) and we thus obtain a homotopy equivalent replacement $(\widetilde{F}_x)_{x\in \Crit(f)}$ of the collection $(F_x)_{x\in \Crit(f)}$ with a strict $\mathcal{GTO}_{top}^{(X,f,\X)}$-action (see Corollary \ref{cor:top action by ass}).
\end{enumerate}

\begin{defi}[Homotopy path algebra structure]\label{defi:homotopy path algebra}    
    Let $(E,\Phi,\mu)$ be a fibration with transitive lifting function and an associative fibrewise monoidal structure $\mu\colon E\ {}_\pi \!\times_\pi E\to E$.
    
    We construct in the following an operad action 
    \begin{align*}
        \varphi\colon W(\mathcal{PO}(X,f))\to \End_{(F_x)_{x\in \Crit(f)}}.
    \end{align*}
    We define the action $\varphi([\eta])$ of $[\eta]\in W(\mathcal{PO}(X,f))(\ux,y)$ on $\underline{\alpha}\in F_{\ux}$ via induction on the number of inputs plus the internal edges of the metric tree $\eta$.

    The operad involved has no operations in arity $0$. Hence the base case is given by $\eta$ being a single vertex labelled by $\lambda\in \mathcal{P}_{x\to y} X$. In this case, we define
    \begin{align*}
        \mu(\Phi)([\eta])({\alpha})=\Phi(\alpha,\lambda).
    \end{align*}
    
    For the induction step, we let $v$ be the last vertex of $\eta$. We may assume that all the inputs of $v$ are internal edges by possibly adding a vertex labelled with the constant path.
    
    We assume that $v$ is labelled by $\lambda_1,\dots,\lambda_l$ with input edges $e_1,\dots,e_l$. We denote by $t_i$ the length of edge $e_i$ and $\eta_i$ for the maximal subtree of $\eta$ that has $e_i$ as output (see Figure \ref{fig:induction_step}).

    \begin{figure}[H]
        \centering
        \includegraphics[width=0.25\linewidth]{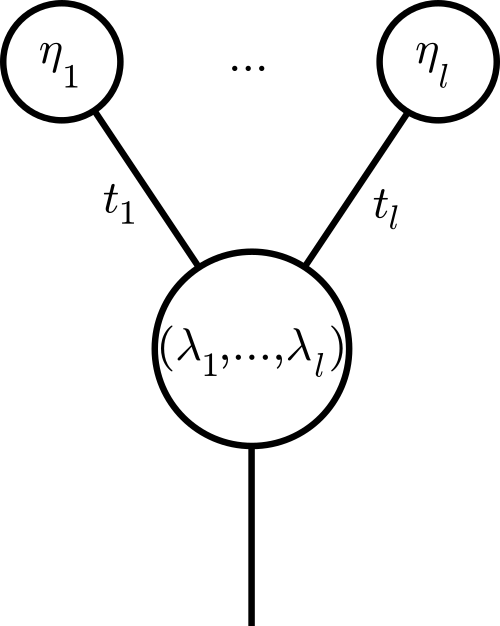}
        \caption{The tree $\eta$ with the last vertex $v$ labelled by $\lambda_1,\dots,\lambda_l$.}
        \label{fig:induction_step}
    \end{figure}
    
    The permutation action on $\mathcal{PO}(X,f)$ and thus on $W(\mathcal{PO}(X,f))$ is free. We thus may assume that the inputs of $\eta_i$ are $j_{i-1}+1,\dots, j_i$ for some $0=j_0< \dots < j_l=k$. We denote $\underline{\alpha}_i\coloneqq (\alpha_{j_{i-1}+1},\dots,\alpha_{j_i})$. We thus have $\underline{\alpha}=(\underline{\alpha}_1,\dots,\underline{\alpha}_l)$.

    For a Moore path $\lambda\colon [0,a]\to X$ and $s\leq t\in [0,1]$, we denote by $\lambda|_{[s,t]}$ the restriction of $\lambda$ to $[sa,ta]$. 
    
    We can now define 
    \begin{align}\label{eq:induction step}
        \varphi([\eta])(\underline{\alpha})\coloneqq \mu\left(\widetilde{\varphi}(t_1,\lambda_1,\eta_1)(\underline{\alpha}_1),\dots, \widetilde{\varphi}(t_l,\lambda_l,\eta_l)(\underline{\alpha}_l)\right) 
    \end{align}
    where
    \begin{align*}
        \widetilde{\varphi}(t_i,\lambda_i,\eta_i)(\underline{\alpha}_i)\coloneqq \Phi(\varphi([\eta_i\cdot \lambda_i|_{[0,1-t_i]}])(\underline{\alpha}_i),\lambda_i|_{[1-t_i,1]})
    \end{align*}
    and $\eta_i\cdot \lambda_i|_{[0,1-t_i]}$ denotes the tree $\eta_i$ where the label of the last vertex is changed from $\underline{\nu}$ to $\lambda_i|_{[0,1-t_i]}\circ_1\underline{\nu}$ (see Figure \ref{fig:etai_lambdai}).
    \begin{figure}[H]
        \centering
        \includegraphics[width=0.5\linewidth]{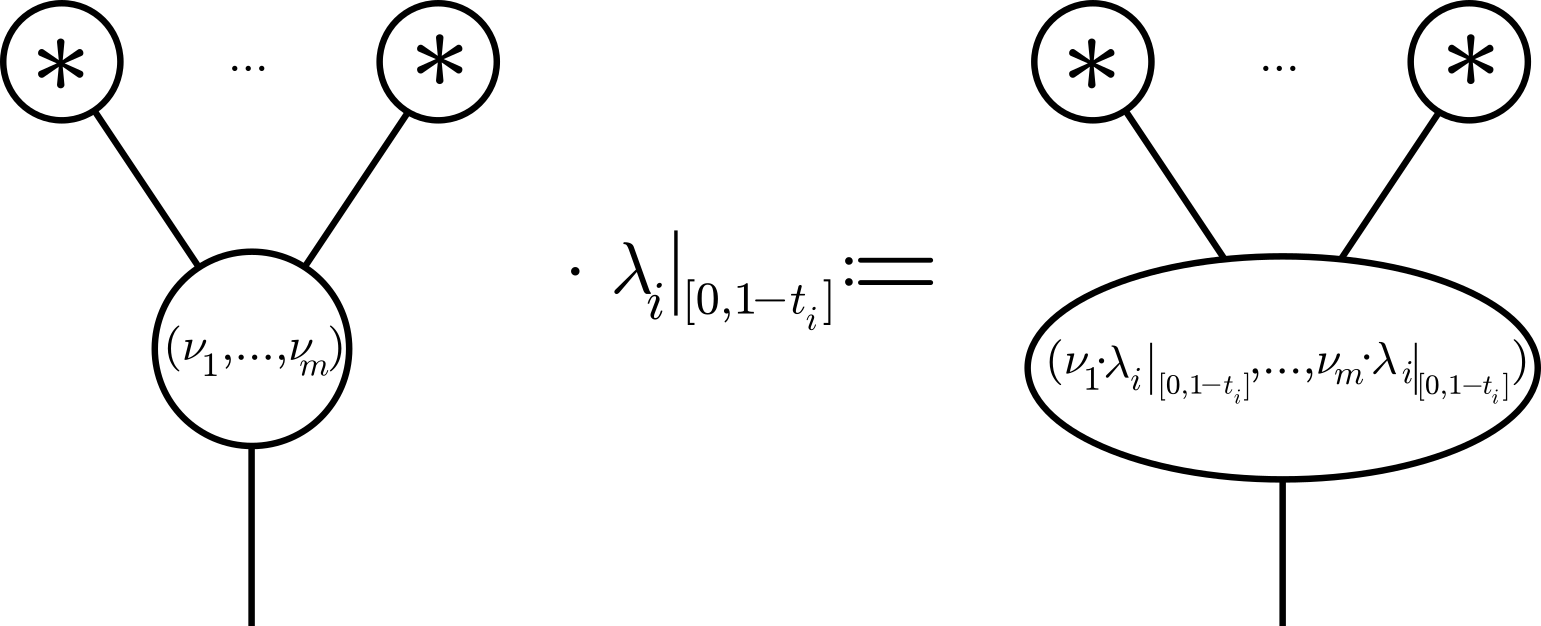}
        \caption{The tree $\eta_i\cdot \lambda_i|_{[0,1-t_i]}$}
        \label{fig:etai_lambdai}
    \end{figure}
\end{defi}

\begin{lemme}\label{lem:homotopy path algebra is well-def}
    Definition \ref{defi:homotopy path algebra} gives a well-defined operad action 
    \begin{align*}
        \varphi \colon W(\mathcal{PO}(X,f))\to \End_{(F_x)_{x\in \Crit(f)} }.
    \end{align*}
\end{lemme}
\begin{proof}
    We need to check that $\varphi$ does not depend on the labelled tree $\eta$ representing the element in $W(\mathcal{PO}(X,f))$.

    We recall that the trees in $W(\mathcal{PO}(X,f))$ are subject to the following two relations:
    \begin{itemize}
        \item contraction of edges of length $0$;
        \item deletion of vertices labelled by the constant path.
    \end{itemize}
    We thus need to show that the induction step of Definition \ref{defi:homotopy path algebra} is invariant under these relations if $e_i$ has length $0$, the last vertex $v$ of $\eta$ is labelled by the constant path or the last vertex $v_i$ of $\eta_i$ is labelled by the constant path.

    We first assume that $t_i=0$. We compute
    \begin{align*}
        \widetilde{\varphi}(t_i,\lambda_i,\eta_i)(\underline{\alpha}_i)=\varphi([\eta_i\cdot \lambda_i])(\underline{\alpha}_i).
    \end{align*}
    The last vertex of $\eta_i\cdot \lambda_i$ is labelled by $\nu_1\cdot \lambda_i,\dots,\nu_m\cdot \lambda_i$ where $\nu_1,\dots,\nu_m$ is the label of the last vertex of $\eta_i$. 

    Associativity of $\mu$ then gives that $\varphi([\eta])$ gives the same operation to the tree where $e_i$ is contracted and $v$ is labelled by $\lambda_1,\dots, \nu_1\cdot \lambda_i,\dots,\nu_m\cdot \lambda_i,\dots,\lambda_l$ which is precisely the composition $\underline{\lambda}\circ_i \underline{\nu}$. This shows that the definition is invariant under edge collapse.

    We next assume that $v$ is labelled by the constant path $c_y$. We have in particular that $v$ has only one input, that is $l=1$. Equation \eqref{eq:induction step} thus reads as
    \begin{align*}
        \varphi([\eta])(\underline{\alpha})=\widetilde{\varphi}(t_1,c_y,\eta_1)=\varphi([\eta_1])(\underline{\alpha})
    \end{align*}
    which is what we wanted to show because $\eta_1$ is the tree we obtain by deleting $v$.

    Finally, we consider the case where the last vertex $v_i$ of $\eta_i$ is labelled by a constant path $c_z$.

    If $v_i$ is the only vertex of $\eta_i$, there is nothing to prove: indeed if we delete $v_i$, $e_i$ is an edge with no tree attached. Hence in the induction step, we add back the unit to apply formula \eqref{eq:induction step}. 

    If $v_i$ is not the only vertex of $\eta_i$, it is preceded by an edge $e'$ of length $s$ and a tree $\eta_i'$ (see Figure \ref{fig:etai'}).

    \begin{figure}[H]
        \centering
        \includegraphics[width=0.2\linewidth]{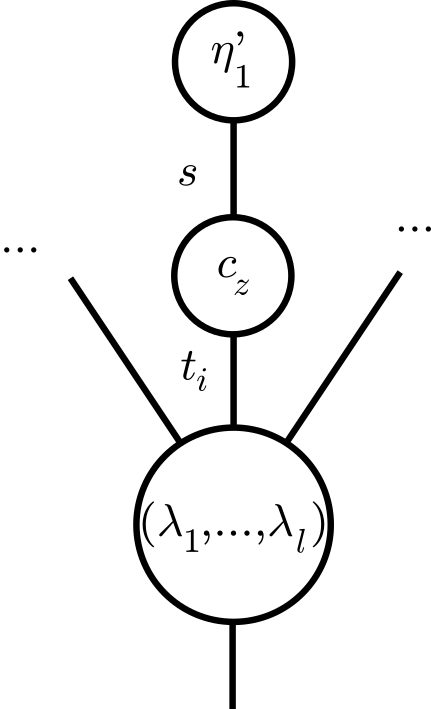}
        \caption{The set up where $v_i$ is labelled by $c_z$ and is preceded by a tree $\eta_i'$.}
        \label{fig:etai'}
    \end{figure}

    The tree $\eta_i\cdot \lambda_i|_{[0,1-t_i]}$ has thus a last vertex labelled by $\lambda_i|_{[0,1-t_i]}$. Applying formula \eqref{eq:induction step} to $\eta_i\cdot \lambda_i|_{[0,1-t_i]}$, we find
    \begin{align*}
        \widetilde{\varphi}(s,\lambda_i|_{[0,1-t_i]},\eta_i')(\underline{\alpha}_i)&=\Phi(\varphi([\eta_i'\cdot (\lambda_i|_{[0,1-t_i]})|_{[0,1-s]}])(\underline{\alpha}_i),(\lambda_i|_{[0,1-t_i]})|_{[1-s,1]})\\
        &=\Phi(\varphi([\eta_i'\cdot \lambda_i|_{[0,(1-t_i)(1-s)]}])(\underline{\alpha}_i),\lambda_i|_{[(1-s)(1-t_i),1-t_i]})\\
        &=\Phi(\varphi([\eta_i'\cdot \lambda_i|_{[0,1-t_i\star s]}])(\underline{\alpha}_i),\lambda_i|_{[1-t_i\star s,1-t_i]}),
    \end{align*}
    where $s\star t_i$ denotes the {reversed multiplication} $1-(1-s)(1-t_i)$ as in Definition \ref{defi:W-const}.
    
    Plugging that back in, we find
    \begin{align*}
        \widetilde{\varphi}(t_i,\lambda_i,\eta_i)(\underline{\alpha}_i)&=\Phi(\Phi(\varphi([\eta_i'\cdot \lambda_i|_{[0,1-t_i\star s]}])(\underline{\alpha}),\lambda_i|_{[1-t_i\star s,1-t_i]}),\lambda_i|_{[1-t_i,1]})\\
        &=\Phi(\varphi([\eta_i'\cdot \lambda_i|_{[0,1-t_i\star s]}])(\underline{\alpha}),\lambda_i|_{[1-t_i\star s,1]})\\
        &=\widetilde{\varphi}(s\star t_i,\lambda_i,\eta_i')(\underline{\alpha}_i).
    \end{align*}
    If we delete $v_i$ beforehand, the newly created edge is labelled by $s\star t_i$ by definition. Therefore the above term is the $i$-th term in Formula \eqref{eq:induction step} for the tree where we deleted $v_i$. 

    This shows that Definition \ref{defi:homotopy path algebra} is invariant under the choice of representative $\eta$.

    It remains to check that $\varphi$ defines an operad action, that is, it commutes with operadic composition. We recall that composition in $W(\mathcal{PO}(X,f))$ is defined by concatenation of trees and labelling the new edge by $1$.

    We assume that in the induction step $t_i=1$ and check that 
    \begin{align}\label{eq:compostion formula}
        \varphi([\eta])=\varphi([\eta']\circ_i [\eta_i])
    \end{align}
    where $\eta'$ is the subtree of $\eta$ given by deleting $\eta_i$.

    We note
    \begin{align*}
        \widetilde{\varphi}(t_i,\lambda_i,\eta_i)(\underline{\alpha}_i)=\Phi(\varphi([\eta_i])(\underline{\alpha}_i),\lambda_i).
    \end{align*}
    Plugging this back into \eqref{eq:induction step}, we find the desired composition formula \eqref{eq:compostion formula}.
\end{proof}

We directly obtain the following result:
\begin{cor}\label{cor:top action by ass}
     Let $(E,\mu)$ be a monoidal fibration. A transitive lifting function $\Phi$ induces an operad action of $W(\mathcal{GTO}_{top}^{(X,f,\X)})$ on $(F_x)_{x\in \Crit(f)}$. 
     
     Moreover, there exist homotopy equivalent spaces $\rho_x\colon F_x\simeq\widetilde{F}_x$ that inherit an operad action
     \begin{align*}
         \mu_\Phi\colon \mathcal{GTO}_{top}^{(X,f,\X)}\to \End_{(\widetilde{F}_x)_{x\in \Crit(f)}}
     \end{align*}
     such that $\rho\colon (F_x)_{x\in \Crit(f)}\to (\widetilde{F}_x)_{x\in \Crit(f)}$ induces a map of $W(\mathcal{GTO}_{top}^{(X,f,\X)})$-algebras.
\end{cor}
\begin{proof}
    By Definition \ref{defi:homotopy path algebra} and Lemma \ref{lem:homotopy path algebra is well-def}, we inherit a $W(\mathcal{PO}(X,f))$-action on $(F_x)_{x\in \Crit(f)}$.

    Using the map of operads $\Gamma\colon \mathcal{GTO}_{top}^{(X,f,\X)}\to \mathcal{PO}(X,f)$ from Remark \ref{rem:evaluating trees} and functoriality of the W-construction, we get a map of operads
    \begin{align*}
        W(\mathcal{GTO}_{top}^{(X,f,\X)})\to W(\mathcal{PO}(X,f)).
    \end{align*}
    Pulling back along this map of operads, we get the structure of a $W(\mathcal{GTO}_{top}^{(X,f,\X)})$-algebra on $(F_x)_{x\in \Crit(f)}$. This proves the first part of the statement.

    To get the second part of the statement, we apply the rectification lemma (see Lemma \ref{lem:W-constr}). We check that $\mathcal{GTO}_{top}^{(X,f,\X)}$ fulfils the assumption of Lemma \ref{lem:W-constr}.

    Indeed, the operation spaces $\mathcal{GTO}_{top}^{(X,f,\X)}(\ux;y)=\Tree_k(\ux;y)$ have a free permutation action because there is an underlying ribbon tree to an element in $\Tree_k(\ux;y)$. Moreover, the operation spaces are manifolds (Lemma \ref{lem:gradient trees are compact manifolds}) and thus CW-complexes.\footnote{This is not true for the operad $\mathcal{PO}(X,f)$ as it is not a CW-complex. This is the reason why we take the detour via the operad $\mathcal{GTO}_{top}^{(X,f,\X)}$ before applying Lemma \ref{lem:W-constr}.} 

    Therefore, we may apply Lemma \ref{lem:W-constr} and obtain the strict replacement $(\widetilde{F}_x)_{x\in \Crit(f)}$.
\end{proof}

\section{From operads to Morse homology}\label{sec:main thm}

In this section, we phrase the geometric construction from Section \ref{sec:geom} in terms of the operad action from Section \ref{sec:top op}. We define a Morse tree operad (see Definition \ref{def:Morse tree op}). We construct an action of this operad on Morse homology with DG coefficients and show that it encodes the generalised Chas-Sullivan product. In fact, we obtain this operation as the binary operation of a shifted $\Ai$-algebra structure (see Theorem \ref{thm:main thm}). We conclude this section by reviewing some examples and comparing them to the literature in Subsection \ref{subsec:exmp}.

\subsection{The Morse tree operad}
In this subsection, we define the Morse tree operad and show that it acts on Morse homology with DG coefficients given a monoidal fibration.

Applying the cubical chains functor to the topological gradient tree operad $\mathcal{GTO}_{top}^{(X,f,\X)}$ (see Definition \ref{def:top gradient tree op}), we obtain an algebraic operad:

\begin{defi}\label{def:alg gradient tree op}
    The \textit{algebraic gradient tree operad associated to $(X,f)$} is the operad given by $\mathcal{GTO}_{alg}^{(X,f,\X)}\coloneqq C_*(\mathcal{GTO}_{top}^{(X,f,\X)})$.
\end{defi}

We recall that $\underline{\Sigma}^{|\cdot |}$ denotes the coloured suspension operad (see Definition \ref{def:coloured susp op}) that encodes suspension of algebras. Secondly, we recall that $\dec$ denotes the decolouring functor that turns a coloured operad into a one-coloured operad (see Subsection \ref{subsec:dec}). 

\begin{cor}\label{cor:dec+shifted operad on fibres}
    For $(E,\mu)$ a monoidal fibration, the operad $\dec(\underline{\Sigma}^{|\cdot|}\otimes \mathcal{GTO}_{alg}^{(X,f,\X)})$ acts on 
    \begin{align*}
        \bigoplus_{x\in \Crit(f)} \Sigma^{|x|}\otimes C_{*}(\widetilde{F}_x)
    \end{align*}
    with the untwisted differential where $\widetilde{F}_x$ denotes the homotopy equivalent fibre from Corollary \ref{cor:top action by ass}.
\end{cor}
\begin{proof}
    This follows from Corollary \ref{cor:top action by ass} combined with Examples \ref{exmp:suspending algebras} and \ref{ex:decolouring algebras}.
\end{proof}

\begin{notation}\label{nota:sxy}
    For $(\ux,y)\in \Crit(f)^{k+1}$, we denote by $\sigma_{\ux,y}$ the element in $\underline{\Sigma}^{|\cdot |}(\ux;y)=\Hom(\Sigma^{|x_1|}\otimes \dots \otimes \Sigma^{|x_k|},\Sigma^{|y|})$ that sends $s^{|x_1|}\otimes \dots \otimes s^{|x_k|}$ to $s^{|y|}$. In the unary case, we have $\sigma_{z,y}\circ_1\sigma_{x,z}=\sigma_{x,y}$.
\end{notation}

\begin{lemme}\label{lem:Barrad-Cornea cocycle is MC}
    Let $\{s_{x,y}\in C_*(\trajb{x,y})=\mathcal{GTO}_{alg}^{(X,f,\X)}(x;y)\}$ be the representing chain system as defined in Definition \ref{defi : representing chain system traj}. The Maurer-Cartan element $\mathfrak{m}=\sum_{x, y} \sigma_{x, y}\otimes s_{x,y}\in \dec(\underline{\Sigma}^{|\cdot |}\otimes \mathcal{GTO}_{alg}^{(X,f,\X)})(1)_{-1}$ satisfies the associative Maurer-Cartan equation:
    \begin{align*}
        \partial(\mathfrak{m})+\mathfrak{m}\circ_1\mathfrak{m}=0.
    \end{align*}
    In particular, $\partial_{\mathfrak{m}}$ defines a twisted differential on $\dec(\underline{\Sigma}^{|\cdot |}\otimes \mathcal{GTO}_{alg}^{(X,f,\X)})$.
\end{lemme}

\begin{proof}
    We note that $|s_{x,y}|=|x|-|y|-1$ and $|\sigma_{x,y}|=|y|-|x|$. Therefore $\mathfrak{m}$ indeed defines a degree minus one element.

    We recall that $\partial(s_{x,y})=\sum_{x>z>y} (-1)^{|x|+|z|}s_{x,z}\times s_{z,y}=\sum_{x>z>y} (-1)^{|x|+|z|+|s_{x,z}||s_{z,y}|} s_{z,y}\circ_1 s_{x,z}$ (see \eqref{eq : relation sys repr}). We compute:
    \begin{align*}
        \partial(\mathfrak{m})&=\partial\left(\sum_{x>y} \sigma_{x,y}\otimes s_{x,y}\right)=\sum_{x>y} (-1)^{|x|+|y|}\sigma_{x,y}\otimes \partial(s_{x,y})\\
        &=\sum_{x>z>y} (-1)^{|y|+|z|+|s_{x,z}||s_{z,y}|} (\sigma_{z,y}\circ_1 \sigma_{x,z})\otimes ( s_{z,y}\circ_1 s_{x,z})\\
        &=\sum_{x>z>y} (-1)^{|y|+|z|+(|s_{x,z}|+|\sigma_{x,z}|)|s_{z,y}|} (\sigma_{z,y}\otimes s_{z,y})\circ_1 ( \sigma_{x,z}\otimes s_{x,z})\\
        &=-\sum_{x>z>y} (\sigma_{z,y}\otimes s_{z,y})\circ_1 ( \sigma_{x,z}\otimes s_{x,z})=-\mathfrak{m}\circ_1\mathfrak{m}.\\
    \end{align*}
    This shows that $\mathfrak{m}$ satisfies the associative Maurer-Cartan equation and thus $\partial_{\mathfrak{m}}$
    defines a differential by Lemma \ref{lem:MC gives twisted diff}.
\end{proof}

\begin{defi}\label{def:Morse tree op}
    The \textit{Morse tree operad associated} to a manifold $X$ with Morse function $f$ and perturbation $\X$, $\mathcal{MT}^{(X,f,\X)}$ is defined as
    \begin{align*}
        (\dec(\underline{\Sigma}^{|\cdot |}\otimes \mathcal{GTO}_{alg}^{(X,f,\X)}),\partial_{\mathfrak{m}}).
    \end{align*}
\end{defi}

\begin{cor}\label{cor:dgM acting on Morse homology}
    Let $(E,\mu)$ be a monoidal fibration. A transitive lifting function induces an $\mathcal{MT}^{(X,f,\X)}$-action on Morse homology with DG coefficients $C_*(X,C_*(\widetilde{F}))$. Moreover, the map
    \begin{align*}
        \rho\colon C_*(X,C_*(F))&\to C_*(X,C_*(\widetilde{F})),\\
        x\otimes \alpha&\mapsto x\otimes \rho_x(\alpha)
    \end{align*}
    defines a quasi-isomorphism.
\end{cor}
\begin{proof}
    We denote by $\phi$ the action of $\dec(\underline{\Sigma}^{|\cdot |}\otimes \mathcal{GTO}_{alg}^{(X,f,\X)})$ from Corollary \ref{cor:dec+shifted operad on fibres}. By Lemma \ref{lem:Barrad-Cornea cocycle is MC}, we get an $\mathcal{MT}^{(X,f,\X)}$-action on $\bigoplus_{x\in \Crit(f)} \Sigma^{|x|}\otimes C_*(\widetilde{F}_x)$ with the differential $\partial_{\phi(\mathfrak{m})}$. 
    
    We first show that $(\bigoplus_{x\in \Crit(f)} \Sigma^{|x|}\otimes C_*(F_x),\partial_{\phi(\mathfrak{m})})$ is isomorphic to $C_*(X,C_*(\widetilde{F}))$, Morse homology with DG coefficients in $C_*(\widetilde{F})$. We claim that 
    \begin{align*}
        \left(\bigoplus_{x\in \Crit(f)} \Sigma^{|x|}\otimes C_*(\widetilde{F}_x),\partial_{\phi(\mathfrak{m})}\right)&\to C_*(X,C_*(\widetilde{F})),\\
        s^{|x|}\otimes \alpha&\mapsto x\otimes \alpha
    \end{align*}
    is an isomorphism of chain complexes. As it is just renaming the element $s^{|x|}$ to $x$, it is an isomorphism of graded abelian groups. It thus remains to check that it is a chain map. We compute for $s^{|x|}\otimes \alpha\in  \Sigma^{|x|}\otimes C_*(\widetilde{F}_x)$:
    \begin{align*}
        \partial_{\phi(\mathfrak{m})}(s^{|x|}\otimes \alpha)&=(-1)^{|x|} s^{|x|}\otimes \partial(\alpha)+\phi(\mathfrak{m})(s^{|x|}\otimes \alpha)\\
        &=(-1)^{|x|} s^{|x|}\otimes \partial(\alpha)+\sum_{x', y}(\sigma_{x',y}\otimes \Phi^{op}(m_{x',y}\times -))(s^{|x|}\otimes \alpha)\\
        &=(-1)^{|x|} s^{|x|}\otimes \partial(\alpha)+\sum_{x, y}(-1)^{|s_{x,y}||x|}\sigma_{x,y}(s^{|x|})\otimes \Phi^{op}(m_{x,y}\times \alpha).
    \end{align*}
    Because $(-1)^{|s_{x,y}||x|}=(-1)^{|x||y|}$, we find the same differential as for $C_*(X,C_*(\widetilde{F}))$ (see Definition \ref{def:Morse hom with DG coeffs}).

    It remains to show that the map $\rho$ is a quasi-isomorphism. We recall that the homotopy equivalences $\rho_x\colon F_x\to \widetilde{F}_x$ commute with the $W(\mathcal{GTO}_{top}^{(X,f,\X)})$-action. 
    
    The differential on both sides is twisted by multiplication with $\sum_{x>y} m_{x,y}$ and $m_{x,y}$ is in $C_*(\mathcal{GTO}_{top}^{(X,f,\X)}(x,y))\subseteq C_*(W(\mathcal{GTO}_{top}^{(X,f,\X)}))$. The fact that the $\rho_x$ commute with the $W(\mathcal{GTO}_{top}^{(X,f,\X)})$-action thus implies that the map $\rho$ is a chain map. It thus induces a map on the first page of the associated spectral sequence
    \begin{align*}
        C_*(X,H_*(F))\to C_*(X,H_*(\widetilde{F}))
    \end{align*}
    which is an isomorphism. This shows that $\rho$ is indeed a quasi-isomorphism.
\end{proof}

\subsection{Main theorem}
In this subsection, we explain how to extend the binary operation from Section \ref{sec:geom} to an $\Ai$-algebra structure.

We recall that we have a map of operads $\Ai^{or}\to (\dec(\Ai^{\Crit(f)}),{\partial}_{m_1})$ (see Corollary \ref{cor:Ai to dec coloured Ai}) and an operad action of $\mathcal{MT}^{(X,f,\X)}$ on $C_*(X,C_*(F))$ (see Corollary \ref{cor:dgM acting on Morse homology}). Moreover, the operad $\mathcal{MT}^{(X,f,\X)}$ is constructed from the $\Crit(f)$-coloured operad $\underline{\Sigma}^{|\cdot|}\otimes\mathcal{GTO}_{alg}^{(X,f,\X)}$. The puzzle piece that lets us thus connect the two worlds is a map of operads from a (shifted version of the) $\Crit(f)$-coloured $\Ai$-operad to $\underline{\Sigma}^{|\cdot|}\otimes\mathcal{GTO}_{alg}^{(X,f,\X)}$. 

\begin{notation}
    We denote by $\underline{n}\colon \Crit(f)\to \Z$ the constant grading $x\mapsto n$. Accordingly, we denote by $\underline{\Sigma}^{\underline{n}}$ the $\Crit(f)$-coloured suspension operad (see Definition \ref{def:coloured susp op}). 
    
    This is to contrast it with the one-coloured suspension operad $\underline{\Sigma}^n$.
\end{notation}

\begin{thm}\label{thm:orientations}
    There exist generators $\sigma_{\ux,y}^n\in\underline{\Sigma}^{\underline{n}}(\ux;y)$ and chains $t_{\underline{x},y}\in C_{\|y\|-\|\underline{x}\|+k-2}(\Tree_k(\underline{x};y))$ such that
    \begin{align*}
        \psi^X\colon  \underline{\Sigma}^{\underline{n}}\otimes \mathcal{A}_\infty^{\Crit(f)}&\to \underline{\Sigma}^{|\cdot |}\otimes \mathcal{GTO}_{alg}^{(X,f,\X)},\\
        \sigma^n_{\ux,y}\otimes m_{\underline{x},y}&\mapsto  \sigma_{\ux,y}\otimes t_{\underline{x},y}
    \end{align*}
    is a morphism of operads. 

    Moreover $\psi_X$ sends $\id_{\Sigma^n} \otimes m_{x,y}\in \underline{\Sigma}^{\underline{n}}\otimes \Ai^{\Crit(f)}(x;y)$ to $\sigma_{x,y}\otimes s_{x,y}$, for $(s_{x,y})_{x,y\in \Crit(f)}$ the representing chain system as in Definition \ref{defi : representing chain system traj}.
\end{thm}

We postpone proving this theorem to Subsection \ref{subsec:proof of orient thm}. We first explain why this lets us prove our main theorem:

\begin{thm}\label{thm:main thm}
    Let $(E,\mu)$ be a monoidal fibration. A transitive lifting function induces a shifted $\Ai$-{algebra} structure on $C_*(X,C_*(\widetilde{F}))$. Concretely, there exists an operad action
    \begin{align*}
        \CS\colon\underline{\Sigma}^n\otimes \Ai^{or}\to \End_{C_*(X,C_*(\widetilde{F}))}
    \end{align*}
    such that the top degree binary operation is a chain-level model of the generalised Chas-Sullivan product on $C_*(X,C_*(\widetilde{F}))\simeq C_*(E)$.

    The top degree $k$-ary operation is given on $ (x_1\otimes \alpha_1)\otimes \dots \otimes (x_k\otimes \alpha_k)\in C_*(X,C_*(\widetilde{F}))^{\otimes k}$ by
    \begin{align*}
        \sum_{y} y\otimes (\mu_\Phi)_*(t_{\ux,y})(\alpha_1\times \dots\times \alpha_k)
    \end{align*}
    for some chains $t_{\ux,y}\in C_*(\Tree_k(\ux;y))$ and $\mu_\Phi$ the operad action from Corollary \ref{cor:top action by ass}.
\end{thm}
\begin{proof}
    Applying the decolouring functor to $\psi^X$ from Theorem \ref{thm:orientations}, we obtain a map of operads
    \begin{align*}
        \dec(\underline{\Sigma}^{\underline{n}}\otimes \Ai^{\Crit(f)})\to \dec(\underline{\Sigma}^{|\cdot |}\otimes \mathcal{GTO}_{alg}^{(X,f,\X)}).
    \end{align*}
    By Corollary \ref{cor:Ai to dec coloured Ai}, we get a morphism of operads
    \begin{align}\label{eq:Sigma n Ai to dec}
        \underline{\Sigma}^n\otimes \Ai^{or}\to \underline{\Sigma}^n\otimes (\dec(\Ai^{\Crit(f)}),\partial_{(\sum m_{x,y})}) \cong (\dec(\underline{\Sigma}^{\underline{n}}\otimes \Ai^{\Crit(f)}),\partial_{(\sum \id_{\Sigma^n}\otimes m_{x,y})}).
    \end{align}
    The isomorphism comes from the fact that suspension and decolouring commute. This is due to distributivity of the tensor product (defining the suspension) and direct sum (defining the decolouring).

    We note that $\dec(\psi_X)$ sends $\sum_{x, y} \id_{\Sigma^n}\otimes m_{x,y}$ to the element $\mathfrak{m}=\sum_{x, y} \sigma_{x,y}\otimes s_{x,y}$ as in Lemma \ref{lem:Barrad-Cornea cocycle is MC}, the Maurer-Cartan element defining $\mathcal{MT}^{(X,f,\X)}$. Thus $\psi_X$ induces a map of operads
    \begin{align*}
        (\dec(\underline{\Sigma}^{\underline{n}}\otimes \Ai^{\Crit(f)}),\partial_{(\sum \id_{\Sigma^n}\otimes m_{x,y})})\to \mathcal{MT}^{(X,f,\X)}.
    \end{align*}
    Precomposing this map with \eqref{eq:Sigma n Ai to dec} and postcomposing with the operad action from Corollary \ref{cor:dgM acting on Morse homology}, we obtain the desired operad action:
    \begin{align*}
         \underline{\Sigma}^n\otimes \Ai^{or}\to  (\dec(\underline{\Sigma}^{\underline{n}}\otimes \Ai^{\Crit(f)}),\partial_{(\sum \id_{\Sigma^n}\otimes m_{x,y})})\to \mathcal{MT}^{(X,f,\X)}\to \End_{C_*(X,C_*(\widetilde{F}))}.
    \end{align*}

    Let $\sigma_k^n\in \underline{\Sigma}^n(k)$ be a generator. The element $\sigma_k^n\otimes m_k$ generates the top degree of $\underline{\Sigma}^n\otimes \Ai^{or}(k)$. Under the map \eqref{eq:Sigma n Ai to dec}, this element is sent to 
    \begin{align*}
        \sum_{\ux, y}\sigma_k^n\otimes m_{\ux,y}\in (\dec(\underline{\Sigma}^{\underline{n}}\otimes \Ai^{\Crit(f)}),\partial_{(\sum \id_{\Sigma^n}\otimes m_{x,y})})
    \end{align*}
    and then under the map from Theorem \ref{thm:orientations} to
    \begin{align*}
        \sum_{\ux, y}\sigma_{\ux,y}\otimes t_{\ux,y}\in \mathcal{MT}^{(X,f,\X)}.
    \end{align*}

    By Corollary \ref{cor:dgM acting on Morse homology}, this action is given up to sign by
    \begin{align*}
        x_1\otimes \alpha_1\otimes \dots \otimes x_k\otimes \alpha_k&\mapsto \sum_{ y} \pm \sigma_{\ux,y}(x_1\otimes \dots\otimes x_k)\otimes (\mu_\Phi)_*(t_{\ux,y})(\alpha_1\otimes \dots\otimes \alpha_k)\\
        &= \sum_{ y} \pm y\otimes (\mu_\Phi)_*((t_{\ux,y})(\alpha_1\times \dots\times \alpha_k)).
    \end{align*}
    After possibly substituting $t_{\ux,y}$ with $-t_{\ux,y}$, this shows the second part of the theorem.
    
    It remains to check that the binary operation models the generalised Chas-Sullivan product.

    We recall that
    \begin{align*}
        \rho\colon C_*(X,C_*(F))&\to C_*(X,C_*(\widetilde{F})),\\
        x\otimes \alpha&\mapsto x\otimes \rho_x(\alpha)
    \end{align*}
    is a quasi-isomorphism and the homotopy equivalences $\rho_x\colon F_x\to \widetilde{F}_x$ commute with the $W(\mathcal{GTO}_{top}^{(X,f,\X)})$-action. 

    We fix times $(t_1,\dots,t_k)\in I^k$ and critical points $x_1,x_2,y,z_1,\dots,z_k\in \Crit(f)$. For perturbed Morse gradient trees $u_{\ux,z_1}\in \Tree_2(\ux;z_1)$, $v_{z_i,z_{i+1}}\in \Tree_1(z_i;z_{i+1})$ and $w_{z_k,y}\in \Tree_1(z_k;y)$, we denote 
    \begin{align*}
        [\eta_{\ux,y}(t_1,\dots,t_k,u_{\ux,z_1},v_{z_1,z_{2}},\dots,v_{z_{k-1},z_{k}},w_{z_k,y})]\in W(\mathcal{GTO}_{top}^{(X,f,\X)})(\ux;y)
    \end{align*}
    for the element that is represented by the following tree: it has vertices $v_0,\dots,v_{k}$ where only $v_0$ has two inputs and is labelled by $u_{\ux,z_1}$ and all the others have one input. For $i=1,\dots,k-1$, the vertex $v_i$ is labelled by $v_{z_i,z_{i+1}}$ and $v_k$ is labelled by $w_{z_k,y}$. The edge $e_i$ connects $v_{i-1}$ to $v_i$, is labelled by $z_i$ and has length $1-t_i$ (see Figure \ref{fig:eta_ux_y}).

    \begin{figure}[!ht]
        \centering
        \includegraphics[width=0.15\linewidth]{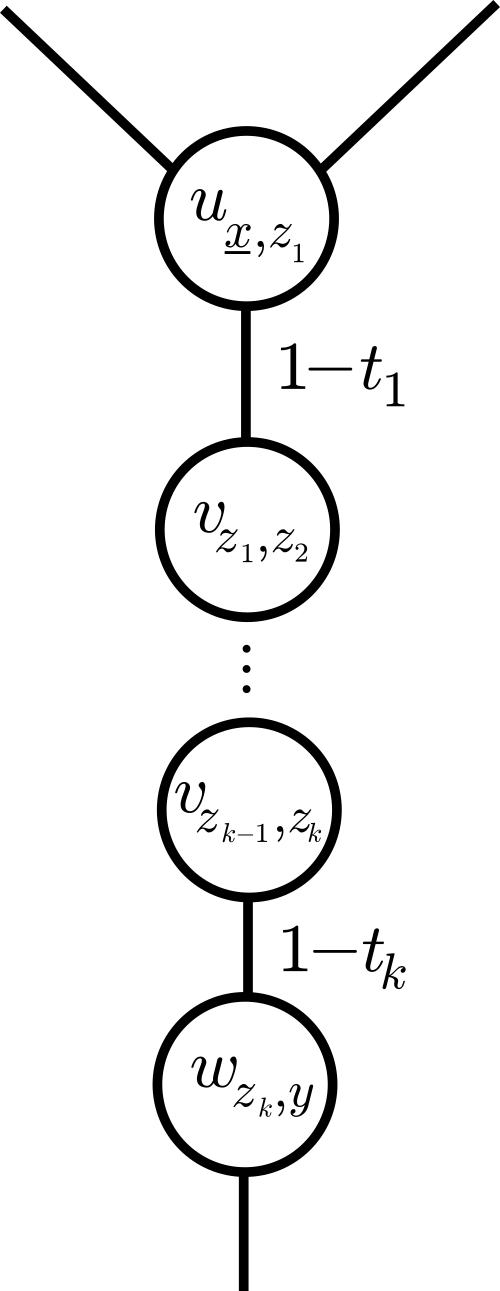}
        \caption{The labelled tree $\eta_{\ux,y}(t_1,\dots,t_k,u_{\ux,z_1},v_{z_1,z_{2}},\dots,v_{z_{k-1},z_{k}},w_{z_k,y})$.}
        \label{fig:eta_ux_y}
    \end{figure}

    We define the following chain
    \begin{align*}
        \theta_{\ux,z_1,\dots,z_k,y}\coloneqq \eta_{\ux,y}(I^{\times k}\times t_{\ux,z_1}\times t_{z_1,z_2}\times\dots\times t_{z_k,y}).
    \end{align*}

    We recall that, by Proposition \ref{prop : binary operation}, the map 
    \begin{align*}
        \mathrm{CS}_{DG}^2((x_1\otimes \alpha_1)\otimes (x_2\otimes \alpha_2))=\sum_y \pm \widetilde{\mu}(y\otimes \Phi^{op}_{2,*}(m_{(x_1,x_2),y}\otimes (\alpha_1\times \alpha_2)))
    \end{align*}
    computes the generalised Chas-Sullivan product. After unravelling the definitions, we find that this map can be rewritten as
    \begin{align*}
        \mathrm{CS}_{DG}^2((x_1\otimes \alpha_1)\otimes (x_2\otimes \alpha_2))=\sum_{y,z_1,\dots,z_k}\pm y\otimes (\mu_\Phi)_*(\theta_{\ux,z_1,\dots,z_k,y})(\alpha_1\times\alpha_2).
    \end{align*}
    We now use that $\rho_x$ commute with the action $\mu_\Phi$ and find
    \begin{align*}
        \rho(\mathrm{CS}_{DG}^2((x_1\otimes \alpha_1)\otimes (x_2\otimes \alpha_2)))&=\sum_{y,z_1,\dots,z_k}\pm y\otimes \rho_y((\mu_\Phi)_*(\theta_{\ux,z_1,\dots,z_k,y})(\alpha_1\times\alpha_2))\\
        &=\sum_{y,z_1,\dots,z_k}\pm y\otimes (\mu_\Phi)_*(\theta_{\ux,z_1,\dots,z_k,y})(\rho_{x_1}(\alpha_1)\times\rho_{x_2}(\alpha_2))).
    \end{align*}
    We recall that the action on $\widetilde{F}$ is strict and therefore the action of 
    \begin{align*}
        \eta_{\ux,y}(t_1,\dots,t_k,u_{\ux,z_1},v_{z_1,z_{2}},\dots,v_{z_{k-1},z_{k}},w_{z_k,y})
    \end{align*}
    does not depend on $t_1,\dots,t_k$. In particular, $(\mu_\Phi)_*(\theta_{\ux,z_1,\dots,z_k,y})(\rho_{x_1}(\alpha_1)\times\rho_{x_2}(\alpha_2)))$ is a degenerate chain if $k\geq 1$. We therefore find
    \begin{align*}
        \rho(\mathrm{CS}_{DG}^2((x_1\otimes \alpha_1)\otimes (x_2\otimes \alpha_2)))&=\sum_{y}\pm y\otimes (\mu_\Phi)_*(\theta_{\ux,y})(\rho_{x_1}(\alpha_1)\times\rho_{x_2}(\alpha_2)))\\
        &=\sum_y\pm y\otimes (\mu_\Phi)_*(t_{\ux,y})(\rho_{x_1}(\alpha_1)\times\rho_{x_2}(\alpha_2))
    \end{align*}
    which is precisely the binary operation evaluated on $\rho(x_1\otimes \alpha_1)\otimes \rho(x_2\otimes \alpha_2)$. This shows that the binary operation computes the generalised Chas-Sullivan product and thus concludes the proof.
\end{proof}

\subsection{Proof of Theorem \ref{thm:orientations}}\label{subsec:proof of orient thm}

In this subsection, we prove Theorem \ref{thm:orientations} assuming the technical Lemma \ref{lem:existence of orientations}, which is proved in Appendix \ref{sec:or trees}.

\begin{notation}
    We use the following notation:
    \begin{itemize}
        \item for an orientable manifold $A$ with possibly empty boundary of dimension $d$, we denote by $\Or(A)\coloneqq H_d(A,\partial A)$;
        \item we denote by $[A] \in \Or(A)$ a choice of generator;
        \item for a point $a\in A$, we denote by $\Or_a(A)$ the top exterior power of the tangent space $T_aA$ and $[A]_a$ for the image of $[A]$ under the projection $\Or(A)\to \Or_a(A)$. By orientability, we find $\Or(A)\cong \bigoplus_{[a]\in \pi_0(A)} \Or_a(A)$;
        
        \item similar to $\sigma_{x,y}$ as in Notation \ref{nota:sxy}, we denote by $\tau_{\ux,y}\in \underline{\Sigma}^{\|\cdot\|}(\ux;y)=\Hom(\Sigma^{\|x_1\|}\otimes \dots \otimes \Sigma^{\|x_k\|},\Sigma^{\|y\|})$  the element that sends $s^{\|x_1\|}\otimes \dots \otimes s^{\|x_k\|}$ to $s^{\|y\|}$.
    \end{itemize}
\end{notation}

By the identity that $|x|+\|x\|=n$, we can fix the degree-preserving isomorphisms
\begin{align}\label{eq:Sigma |x| to Sigma n}
    \begin{split}
        \Sigma^{|x|}\otimes \Sigma^{\|x\|}&\to \Sigma^n,\\
        s^{|x|}\otimes s^{\|x\|}&\mapsto s^n.
    \end{split}
\end{align}

\begin{lemme}\label{lem:sxy and tau xy}
    The maps as in \eqref{eq:Sigma |x| to Sigma n} induce a map of operads
    \begin{align*}
        \underline{\Sigma}^{|\cdot |}\otimes \underline{\Sigma}^{\|\cdot \|}\cong \underline{\Sigma}^{\underline{n}}
    \end{align*}
    that, in the unary case, sends $\sigma_{x,y}\otimes \tau_{x,y}$ to $(-1)^{|x|+|x||y|}\id_{\Sigma^n}$.
\end{lemme}
\begin{proof}
    We fix $(\ux,y)\in \Crit(f)^{k+1}$. An element $\sigma \otimes \tau \in \underline{\Sigma}^{|\cdot |}\otimes \underline{\Sigma}^{\|\cdot \|}(\ux;y)$ is sent to the following map in $\underline{\Sigma}^{\underline{n}}(\ux;y)=\Hom((\Sigma^n)^{\otimes k},\Sigma^n)$:
    \begin{align*}
        (\Sigma^n)^{\otimes k}&\overset{\eqref{eq:Sigma |x| to Sigma n}}{\cong }\Sigma^{|x_1|}\otimes \Sigma^{\|x_1\|}\otimes \dots \otimes \Sigma^{|x_k|}\otimes \Sigma^{\|x_k\|}\cong \Sigma^{|x_1|}\otimes\dots \otimes \Sigma^{|x_k|}\otimes \Sigma^{\|x_1\|} \otimes\dots \otimes \Sigma^{\|x_k\|} \\
        &\overset{\sigma\otimes \tau}{\to } \Sigma^{|y|}\otimes \Sigma^{\|y\|}\overset{\eqref{eq:Sigma |x| to Sigma n}}{\cong }\Sigma^n.
    \end{align*}
    This is natural and thus defines a map of operads.
    
    In the particular case $k=1$, we find:
    \begin{align*}
        s^n & \overset{\eqref{eq:Sigma |x| to Sigma n}}{\mapsto} s^{|x|}\otimes s^{\|x\|}\\
        & \!\!\!\!\!\!\!\!\, \overset{\sigma_{x,y}\otimes \tau_{x,y}}{\mapsto} (-1)^{|\tau_{x,y}||x|} \sigma_{x,y}(s^{|x|})\otimes \tau_{x,y}(s^{\|x\|})=(-1)^{|x|+|x||y|} s^{|y|}\otimes s^{\|y\|}\\
        &\overset{\eqref{eq:Sigma |x| to Sigma n}}{\mapsto} (-1)^{|x|+|x||y|} s^n=(-1)^{|x|+|x||y|} \id_{\Sigma^n}(s^n).
    \end{align*}
\end{proof}

\begin{rem}
    The generators $\sigma_{\ux,y}$ and $\tau_{\ux,y}$ are not canonical as illustrated by the signs in Lemma \ref{lem:sxy and tau xy}. In fact, it is possible to prove everything we prove with other generators of $\underline{\Sigma}^{|\cdot|}(\ux;y)$ and $\underline{\Sigma}^{\|\cdot \|}(\ux;y)$. They can be thought of as variables, but in the following, it is sometimes useful to have concrete generators fixed to facilitate some computations.
\end{rem}

By Lemma \ref{lem:suspending operads}, tensoring with $\underline{\Sigma}^{|\cdot|}$ is an equivalence. It therefore is equivalent to produce a map $\underline{\Sigma}^{\|\cdot \|}\otimes \Ai^{\Crit(f)}\to \mathcal{GTO}_{alg}^{(X,f,\X)}$ and to produce a map 
\begin{align*}
    \Sigma^n\otimes \Ai^{\Crit(f)}\cong\Sigma^{|\cdot|}\otimes \Sigma^{\|\cdot \|}\otimes \Ai^{\Crit(f)}\to \Sigma^{|\cdot |}\otimes \mathcal{GTO}_{alg}^{(X,f,\X)}.
\end{align*}
The former is easier as only one suspension operad is involved. 

The operad $\underline{\Sigma}^{\|\cdot \|}\otimes \Ai^{\Crit(f)}$ is generated by the elements $\tau_{\ux,y}\otimes m_{\ux,y}$. They satisfy the relations of the form
\begin{align}\label{eq:eps0 defi}
    \partial(\tau_{\ux,y}\otimes m_{\ux,y})=\sum_{\substack{\uw\circ_{i}\uz=\ux}} \varepsilon_0(\uw;\uz;i) (\tau_{\uw,y}\otimes m_{\uw,y})\circ _{i}(\tau_{\uz,w_{i}}\otimes m_{\uz,w_{i}})
\end{align}
for certain signs $\varepsilon_0(\uw;\uz;i)\in \{\pm 1\}$. The signs $\varepsilon_0(\uw;\uz;i)$ are directly computable (see Lemma \ref{lem:eps0}).

Giving a map $\underline{\Sigma}^{\|\cdot \|}\otimes \Ai^{\Crit(f)}\to \mathcal{GTO}_{alg}^{(X,f,\X)}$ thus comes down to producing elements $t_{\ux,y}\in \mathcal{GTO}_{alg}^{(X,f,\X)}$ that satisfy the relations
\begin{align*}
    \partial(t_{\ux,y})=\sum_{\uw\circ_i \uz=\ux} \varepsilon_0(\uw;\uz;i) t_{\uw,y}\circ_it_{\uz,w_i}.
\end{align*}

We recall that $\mathcal{GTO}_{alg}^{(X,f,\X)}$ is given by chains on the moduli space of Morse trees. We compute the degree $|\tau_{\ux,y}\otimes m_{\ux,y}|= \ell(\ux)-2-\|\ux\|+\|y\|$ which is precisely the top (relative) homology of 
\begin{align*}
    \mathcal{GTO}_{alg}^{(X,f,\X)}(\ux;y)=C_*(\Tree_k(\ux;y))
\end{align*}
by Lemma \ref{lem:gradient trees are compact manifolds}.

Thus the following lemma is the key to produce the map $\underline{\Sigma}^{\|\cdot \|}\otimes \Ai^{\Crit(f)}\to \mathcal{GTO}_{alg}^{(X,f,\X)}$:

\begin{lemme}[{Lemma \ref{lem:existence of orientations appendix}}]\label{lem:existence of orientations}
    There exist orientations $[\Tree_k(\ux;y)]\in \Or(\Tree_k(\ux;y))$ such that, for all points $p\in \mathcal{T}_l^{\mathbb{X}}(\uw;y)\times \mathcal{T}^{\mathbb{X}}_m(\uz;w_{i})\subseteq \partial \Tree_k(\ux;y)$, it holds
    \begin{align}\label{eq:orient with n}
        [\Tree_k(\ux;y)]_p=\varepsilon_0(\uw;\uz;i)[\vec{n}]\times [\Tree_l(\uw;y)]_p\times [\Tree_m(\uz;w_{i})]_p\in \Or_p(\Tree_{k}(\ux;y))
    \end{align}
    where $\vec{n}\in T_p\Tree_{k}(\ux;y)$ is a vector pointing out of the manifold and $\varepsilon_0(\uw;\uz;i)$ are as in \eqref{eq:eps0 defi}.
    
    Moreover for $k=1$, we have $[\Tree_1(x;y)]=(-1)^{|x|+|x||y|}[\overline{\mathcal{L}}(x,y)]$ for $[\overline{\mathcal{L}}(x,y)]$ as in \eqref{eq: or L(x,y)}.
\end{lemme}

\begin{rem}\label{rem:interpreting bound orient}
    The relevance of the vector $\vec{n}$ is the following: for a $d$-manifold $M$, a point $p$ in its boundary and a vector $\vec{n}$ pointing out, the map 
    \begin{align*}
        H_{d-1}(\partial M)\cong \Or_p(\partial M)&\to \Or_p(M)\cong H_d(M),\\
        \alpha&\mapsto \vec{n}\times \alpha
    \end{align*}
    is an inverse to the connecting isomorphism given by taking the differential. Hence, we may interpret \eqref{eq:orient with n} as ensuring that the map
    \begin{align*}
        \Or(\Tree_k(\ux;y))&\cong H_{\dim \Tree_k(\ux;y)} (\Tree_k(\ux;y),\partial \Tree_k(\ux;y))\\
        &\overset{\partial}{\to} H_{\dim \Tree_k(\ux;y)-1}(\partial \Tree_k(\ux;y))\to\Or_p(\mathcal{T}_l^{\mathbb{X}}(\uw;y)\times \mathcal{T}^{\mathbb{X}}_m(\uz;w_{i}))
    \end{align*}
    sends $[\Tree_k(\ux;y)]$ to $\varepsilon_0(\uw;\uz;i) [\Tree_l(\uw;y)]\times [\Tree_m(\uz;w_{i})]$.
\end{rem}

\begin{proof}[Proof of Theorem \ref{thm:orientations}]
    To prove the theorem, we need to construct a map of operads
    \begin{align*}
        \psi_X\colon \underline{\Sigma}^{\underline{n}}\otimes \Ai^{\Crit(f)}&\to \underline{\Sigma}^{|\cdot |}\otimes \mathcal{GTO}_{alg}^{(X,f,\X)}
    \end{align*}
    that sends $\id_{\Sigma^n}\otimes m_{x,y}$ to $\sigma_{x,y}\otimes s_{x,y}$ where $s_{x,y}$ is the representing chain system as in Definition \ref{defi : representing chain system traj}.

    We first define a map 
    \begin{align*}
        \underline{\Sigma}^{\|\cdot \|}\otimes \Ai^{\Crit(f)}\to \mathcal{GTO}_{alg}^{(X,f,\X)}.
    \end{align*}
    
    If $\Tree_k(\ux;y)$ is empty, we define $t_{\ux,y}=0$ and, if $k=1$, we define $t_{x,y}\coloneqq(-1)^{|x|+|x||y|}s_{x,y}$.  Inductively on the dimension of $\Tree_k(\ux;y)$, we define chains $t_{\ux,y}\in C_*(\Tree_k(\ux;y))$ representing $[\Tree_k(\ux;y)]\in H_*(\Tree_k(\ux;y),\partial\Tree_k(\ux;y))$ such that
    \begin{align*}
        \partial(t_{\ux,y})=\sum_{\uw\circ_i \uz=\ux} \varepsilon_0(\uw;\uz;i) t_{\uw,y}\circ_it_{\uz,w_i}
    \end{align*}
    where $\varepsilon_0(\uw;\uz;i)$ are the signs appearing in 
    \begin{align*}
        \partial(\tau_{\ux,y}\otimes m_{\ux,y})=\sum_{{\uw\circ_{i}\uz=\ux}} \varepsilon_0(\uw;\uz;i) (\tau_{\uw,y}\otimes m_{\uw,y})\circ _{i}(\tau_{\uz,w_{i}}\otimes m_{\uz,w_{i}}).
    \end{align*}
    This is possible by Lemma \ref{lem:existence of orientations} and Remark \ref{rem:interpreting bound orient}.

    We claim that the map
    \begin{align*}
        \psi'_X\colon \underline{\Sigma}^{\|\cdot \|}\otimes \Ai^{\Crit(f)}&\to \mathcal{GTO}_{alg}^{(X,f,\X)},\\
        \tau_{\ux,y}\otimes m_{\ux,y}&\mapsto t_{\ux,y}
    \end{align*}
    defines a map of operads.

    Because $\underline{\Sigma}^{\|\cdot \|}\otimes \Ai^{\Crit(f)}$ is freely generated by the elements $\tau_{\ux,y}\otimes m_{\ux,y}$ as an operad enriched in graded abelian groups, $\psi'_X$ defines a unique map of operads enriched in graded abelian groups.

    It remains to check that $\psi_X'$ is a chain map. In fact, it suffices to check that $\psi_X'$ commutes with the differential $\partial$ on generators $\tau_{\ux,y}\otimes m_{\ux,y}$. We thus compute:
    \begin{align*}
        \partial(\psi_X'(\tau_{\ux,y}\otimes m_{\ux,y}))=\partial(t_{x,y})&=\sum_{\uw\circ_i \uz=\ux} \varepsilon_0(\uw;\uz;i) t_{\uw,y}\circ_it_{\uz,w_i}\\
        &=\sum_{{\uw\circ_{i}\uz=\ux}} \varepsilon_0(\uw;\uz;i) \psi_X'(\tau_{\uw,y}\otimes m_{\uw,y})\circ_{i} \psi_X'(\tau_{\uz,w_{i}}\otimes m_{\uz,w_{i}})\\
        &=\psi_X'\left(\sum_{\substack{\uw\circ_{i}\uz=\ux}} \varepsilon_0(\uw;\uz;i) (\tau_{\uw,y}\otimes m_{\uw,y})\circ _{i}(\tau_{\uz,w_{i}}\otimes m_{\uz,w_{i}})\right)\\
        &=\psi'_X(\partial(\tau_{\ux,y}\otimes m_{\ux,y})).
    \end{align*}

    This shows that we have a map of operads $\psi_X'\colon  \underline{\Sigma}^{\|\cdot \|}\otimes \Ai^{\Crit(f)}\to \mathcal{GTO}_{alg}^{(X,f,\X)}$. We then define $\psi_X$ as the composition 
    \begin{align*}
        \underline{\Sigma}^{n}\otimes \Ai^{\Crit(f)}\overset{\text{Lemma }\ref{lem:sxy and tau xy}}{\cong } \underline{\Sigma}^{|\cdot |}\otimes \underline{\Sigma}^{\|\cdot \|}\otimes \Ai^{\Crit(f)}\overset{\underline{\Sigma}^{|\cdot |} \psi_X'}{\to} \underline{\Sigma}^{|\cdot |}\otimes \mathcal{GTO}_{alg}^{(X,f,\X)}.
    \end{align*}
    We find by Lemmas \ref{lem:sxy and tau xy} and \ref{lem:existence of orientations}:
    \begin{align*}
        \psi_X(\id_{\Sigma^n}\otimes m_{x,y})\overset{\ref{lem:sxy and tau xy}}{=}(-1)^{|x|+|x||y|} \sigma_{x,y} \otimes \psi'_X(\tau_{x,y}\otimes m_{x,y})=(-1)^{|x|+|x||y|}\sigma_{x,y}\otimes t_{x,y}\overset{\ref{lem:existence of orientations}}{=}\sigma_{x,y}\otimes s_{x,y}.
    \end{align*}
\end{proof}

\subsection{Examples}\label{subsec:exmp}
In this subsection, we discuss applications of Theorem \ref{thm:main thm} and we compare them to existing results in the literature.

In particular, we focus on the binary operation induced on homology by Theorem \ref{thm:main thm} which gives an operation of degree $-n$ that is graded associative in the sense of \cite{cieliebak2022topological,clivio2025graded}:
\begin{lemme}\label{lem:graded ass}
    For a generator $\sigma_2^n\in \underline{\Sigma}^n(2)$, the element $\sigma_2^n\otimes m_2\in \underline{\Sigma}^n\otimes \Ai^{or}$ induces a binary operation
    \begin{align*}
        \operatorname{CS}^2\colon H_*(E)\otimes H_*(E)\to H_{-n+*}(E)
    \end{align*}
    under the map from Theorem \ref{thm:main thm}. This operation satisfies
    \begin{align*}
        \operatorname{CS}^2\circ (\operatorname{CS}^2\otimes \id)=(-1)^n\operatorname{CS}^2\circ (\id \otimes \operatorname{CS}^2).
    \end{align*}
\end{lemme}
\begin{proof}
    Let $\sigma_2^n\in \underline{\Sigma}^n(2)$ be a generator.  We define $\operatorname{CS}^2\coloneqq H_*(\operatorname{CS}(\sigma^n_2\otimes m_2))$. It acts on $H_*(C_*(X,C_*(F)))\cong H_*(E)$.
    
    Evaluating both sides on $s^n\otimes s^n\otimes s^n$, we find 
    \begin{align*}
        \sigma_3^n\coloneqq\sigma^n_2\circ_1\sigma^n_2=(-1)^n\sigma_2^n\circ_2\sigma_2^n.
    \end{align*}
    We thus compute in $\underline{\Sigma}^n\otimes \Ai^{or}$:
    \begin{align*}
        \partial(\sigma_3^n\otimes m_3)&= \sigma^n_3\otimes \partial(m_3)= \sigma^n_3\otimes(m_2\circ_2 m_2-m_2\circ_1 m_2)\\
        &=(-1)^n(\sigma^n_2\circ_2\sigma^n_2)\otimes (m_2\circ_2m_2)-(\sigma_2\circ_1\sigma_2^n)\otimes (m_2\circ_1 m_2)\\
        &=(-1)^n(\sigma_2^n\otimes m_2)\circ_2(\sigma^n_2\otimes m_2) - (\sigma_2^n\otimes m_2)\circ_1(\sigma^n_2\otimes m_2).
    \end{align*}
    In homology, this gives $(\sigma_2^n\otimes m_2)\circ_1(\sigma^n_2\otimes m_2) =(-1)^n (\sigma_2^n\otimes m_2)\circ_2(\sigma^n_2\otimes m_2)$. Therefore also its image, $\operatorname{CS}^2$, satisfies
    \begin{align*}
        \operatorname{CS}^2\circ_1 \operatorname{CS}^2=(-1)^n\operatorname{CS}^2\circ_2  \operatorname{CS}^2
    \end{align*}
    which shows the desired relation.
\end{proof}

\begin{ex}\label{exmp:triv}
    For the trivial fibration $\id\colon X\to X$, Morse homology with DG coefficients computes singular homology of $X$. The product that is computed is the intersection product (see \cite[\S 2.1]{Abo11}).

    In particular by Lemma \ref{lem:graded ass}, it satisfies graded associativity. After choosing an appropriate generator $\sigma_2^n$, it corresponds to the dual of the Thom coproduct $\nu_{Th}$ in \cite[\S 6.1]{clivio2025graded}. It also corresponds to the dual of the coproduct $\boldsymbol{\nu}$ in \cite[\S 7]{cieliebak2022topological} after choosing an appropriate generator $\sigma^n_2$. We discuss below in \S \ref{subsubsec:choice of generators} which choices of $\sigma_2^n$ correspond to which convention.
\end{ex}

\begin{ex}\label{exmp:loop-loop}
    By Example \ref{ex : assumption holds for loop loop}, there exists a transitive lifting function $\Phi$ such that $(\ls X, \Phi, \textup{concat})$ is a lift-preserving monoidal fibration. We can thus apply Theorem \ref{thm:main thm} without even replacing the fibres. This computes the map $\CS^2$ defined in Proposition \ref{prop : binary operation}. The second author showed in \cite[Theorem 7.1]{Rie24} that this map computes the Chas-Sullivan product as in \cite{ChaSu99}. 

    In particular, our sign convention provides us with a graded associative product.  This product agrees with the definition of the loop product in \cite[Definition 1.2]{hingston2022product} up to choosing the correct generator $\sigma_2^n$ (see \S \ref{subsubsec:choice of generators}). The signs from their \enquote{algebraic loop product} \cite[Definition 1.3]{hingston2022product} can be recovered by applying $\underline{\Sigma}^{-n}$ and thus getting an associative operation on $\underline{\Sigma}^{-n}\otimes H_*(\mathcal{L}X)= H_{n+*}(\mathcal{L}X)$.
    
\end{ex}

\begin{rem}
    The above example shows the operad action $\mathcal{MT}^{(X,f,\X)}$ on $C_*(X,C_*(\Omega X))$ recovers the original Chas-Sullivan product or string topology product \cite{ChaSu99}. It is an open problem whether this action \enquote{knows} more than the string topology product which only depends on the homotopy type of $X$ \cite{cohen2008homotopy}.

    We suggest questions of the following type for future research: does the operad action $\mathcal{MT}^{(X,f,\X)}$ on $C_*(X,C_*(\Omega X))$ depend on the perturbation $\X$, the Morse function $f$, the smooth structure of $X$, the homeomorphism type of $X$ or the simple homotopy type of $X$?
\end{rem}

\subsubsection{Choice of generator}\label{subsubsec:choice of generators}

\begin{notation}
    For any $x_1,x_2,y\in \Crit(f)$, we denote by $\kappa_{\ux,y}\in \underline{\Sigma}^{\underline{n}}(x_1,x_2;y)$ the map that sends $s^n\otimes s^n$ to $s^n$. We denote by $\kappa_n$ the uncoloured version, i.e., the element in $\underline{\Sigma}^n(2)=\Hom(\Sigma^n\otimes \Sigma^n,\Sigma^n)$ that sends $s^n\otimes s^n$ to $s^n$.
\end{notation}

\begin{lemme}\label{lem:sigma tau binary}
    The identification $\underline{\Sigma}^{|\cdot |}\otimes \underline{\Sigma}^{\|\cdot \|}\cong \underline{\Sigma}^{\underline{n}}$ from Lemma \ref{lem:sxy and tau xy} sends $\sigma_{\ux,y}\otimes \tau_{\ux,y}\in \underline{\Sigma}^{|\cdot |}\otimes \underline{\Sigma}^{\|\cdot \|}(x_1,x_2;y)$ to
    \begin{align*}
        (-1)^{\|x_1\||x_2|+(\|x_1\|+\|x_2\|)(1+\|y\|)}\kappa_{\ux,y}.
    \end{align*}
\end{lemme}
\begin{proof}
    We recall that the map from Lemma \ref{lem:sxy and tau xy} are induced by the identification
    \begin{align*}
        \Sigma^{|x|}\otimes \Sigma^{\|x\|}&\to \Sigma^{{n}},\\
        s^{|x|}\otimes s^{\|x\|}&\mapsto s^n.
    \end{align*}
    Under this map, we have:
    \begin{align*}
        s^n\otimes s^n&\mapsto s^{|x_1|}\otimes s^{\|x_1\|}\otimes s^{|x_2|}\otimes s^{\|x_2\|}\\
        &\mapsto (-1)^{\|x_1\||x_2|} s^{|x_1|}\otimes s^{|x_2|}\otimes s^{\|x_1\|}\otimes s^{\|x_2\|}.
    \end{align*}
    Applying $\sigma_{\ux,y}\otimes \tau_{\ux,y}$, we find:
    \begin{align*}
        (-1)^{\|x_1\||x_2|} s^{|x_1|}\otimes s^{|x_2|}\otimes s^{\|x_1\|}\otimes s^{\|x_2\|}&\mapsto (-1)^{\|x_1\||x_2|+(|x_1|+|x_2|)|\tau_{\ux,y}|} s^{|y|}\otimes s^{\|y\|}\\
        &\mapsto (-1)^{\|x_1\||x_2|+(\|x_1\|+\|x_2\|)(\|x_1\|+\|x_2\|+\|y\|)} s^n\\
        &=(-1)^{\|x_1\||x_2|+(\|x_1\|+\|x_2\|)(\|y\|+1)} s^n.
    \end{align*}
\end{proof}

In Example \ref{exmp:triv}, we argue that, for the trivial fibration, the binary operation computes the intersection product in a sign convention that is graded associative as in \cite[\S 7]{cieliebak2022topological} and \cite[\S 6.1]{clivio2025graded}. However these two sign conventions are slightly different. In \cite{cieliebak2022topological}, it holds
\begin{align*}
    [X]\cdot [*]=(-1)^n[*]
\end{align*}
and, in \cite{clivio2025graded}, it holds
\begin{align*}
    [X]\cdot [*]=(-1)^{\frac{n(n+1)}{2}}[*]
\end{align*}
for $[*]\in H_0(X)$ the element given by the point.

\begin{prop}
    For $\operatorname{CS}^2=H_*(\operatorname{CS}(\kappa_n\otimes m_2))$, it holds $\operatorname{CS}^2([X]\otimes [*])=(-1)^n[*]$.
\end{prop}
\begin{proof}
    We pick a Morse function $f\colon X\to \R$ with a unique critical point of maximal index $n$, $x_+$, a unique critical point of minimal index $0$, $x_-$ (and arbitrarily many critical points of other indices). It thus holds $\|x_+\|=0$ and $\|x_-\|=n$. 

    We pick orientations $[\wb^u(x_+)]\coloneqq[X]$ and $[\wb^u(x_-)]\coloneqq[*]$. We have that $\operatorname{CS}^2([X]\otimes [*])=[\operatorname{CS}(\kappa_n\otimes m_2)(x_+\otimes x_-)]$. Indeed $[X]\otimes [*]$ is represented by $(x_+ + \sum_{y\neq x_+}\lambda_y y)\otimes x_-$ in Morse homology. But for degree reasons, the terms $y\otimes x_-$ get sent to zero by $\operatorname{CS}(\kappa_n\otimes m_2)$ because $|y|<n$. 
    
    Moreover, $\operatorname{CS}(\kappa_n\otimes m_2)(x_+\otimes x_-)$ is a multiple of $x_-$ for degree reasons. In fact, $X$ intersects the point transversely and $X\cap *=*$. We therefore have $\operatorname{CS}(\kappa_n\otimes m_2)(x_+\otimes x_-)=\pm x_-$.

    We set $\ux\coloneqq(x_+,x_-)$ and $y\coloneqq x_-$. By the above orientation choices, we have $[\wb^s(y)]=[X]$. Let $[\nu_{\ux,y}]\in \Or(T_yX^2/T_y \Delta(W^s(y)))$ be as in Lemma \ref{lem:explicit orient k>1}. This corresponds to the orientation $[X]\times [*]$ because the defining equation of $[\nu_{\ux,y}]$ holds:
    \begin{align*}
        [X]\times [*]\times \Delta_*([X])=[X]\times [X]\in \Or_y(X).
    \end{align*}
    We next compute $\eta(\ux;y)$ as in Lemma \ref{lem:explicit orient k>1}. The formula gives:
    \begin{align*}
        \eta(\ux;y)&=\sum_{j}(k-j)n\|x_j\|+\sum_{j>j'}\|x_j\|\|x_{j'}\|+\|y\|\\
        &=n\|x_1\|+\|x_2\|\|x_1\|+\|y\|=n.
    \end{align*}
    This shows that 
    \begin{align*}
        [X]\times [*]\times [\Tree_2(x_1,x_2;y)]_p=(-1)^n([\wb^u(\ux)]\times [\overline{\mathcal{T}_k}])_p=(-1)^n([X]\times [*]\times [*])_p.
    \end{align*}
    The upshot is that $[\Tree_2(x_1,x_2;y)]=(-1)^n[*]$ on each connected component of $\Tree_2(x_1,x_2;y)$. We can thus choose $t_{\ux,y}\coloneqq(-1)^n\sum_i t_i$ where $t_i$ are trees in the connected components of $\Tree_2(x_1,x_2;y)$. The map 
    \begin{align*}
        \psi_X'\colon \underline{\Sigma}^{\|\cdot \|}\otimes \Ai^{\Crit(f)}\to \mathcal{GTO}_{alg}^{(X,f,\X)}
    \end{align*}
    thus sends $\tau_{\ux,y}\otimes m_{\ux,y}$ to $t_{\ux,y}$. 

    By Lemma \ref{lem:sigma tau binary}, the induced map $\psi_X\colon \underline{\Sigma}^{\underline{n}}\otimes \Ai^{\Crit(f)}\to \underline{\Sigma}^{|\cdot |}\otimes \mathcal{GTO}_{alg}^{(X,f,\X)}$ thus sends $\kappa_{\ux,y}\otimes m_{\ux,y}$ to 
    \begin{align*}
        (-1)^{\|x_+\||x_-|+(\|x_+\|+\|x_-\|)(1+\|x_-\|)}\sigma_{\ux,y}\otimes t_{\ux,y}=\sigma_{\ux,y}\otimes t_{\ux,y}.
    \end{align*}
    We conclude by computing:
    \begin{align*}
        \operatorname{CS}(\kappa_n\otimes m_2)(x_+\otimes x_-)&= \psi_X(\kappa_{\ux,y}\otimes m_{\ux,y})(s^{|x_+|}\otimes [*]\otimes s^{|x_-|}\otimes [*])\\
        &=(\sigma_{\ux,y}\otimes (\Phi_\mu)_*(t_{\ux,y}))(s^{|x_+|}\otimes s^{|x_-|}\otimes [*]\times [*])\\
        &=\sigma_{\ux,y}(s^{|x_+|}\otimes s^{|x_-|})\otimes (\mu_\Phi)_*(t_{\ux,y})([*]\times [*])\\
        &=(-1)^n x_-\otimes \sum_i  (\mu_\Phi)_*(t_i)([*]\times [*])=(-1)^n\sum_i x_-\otimes [*].
    \end{align*}
    By the previous discussion, this sum equals $\pm x_-\otimes [*]$. This shows that this sum is just one term. In fact, this shows that $\Tree_2(\ux;y)$ is connected. We moreover find the desired relation:
    \begin{align*}
        \operatorname{CS}^2([X]\otimes [*])=(-1)^n[x_-\otimes [*]]=(-1)^n[*]\in H_0(X).
    \end{align*}
\end{proof}

\begin{ex}
    The previous proposition shows that the map $H_*(\operatorname{CS}(\sigma^n_2\otimes m_2))$ produces the sign convention from \cite{cieliebak2022topological} if we pick $\sigma_2^n\coloneqq\kappa_n$. To recover the sign convention from \cite{clivio2025graded}, we must pick the $\sigma^n_2\coloneqq(-1)^{\frac{n(n+1)}{2}}\kappa_n=(\kappa_1)^{\otimes n}$. 
\end{ex}

\begin{ex}
    In Example \ref{exmp:loop-loop}, we showed that for the loop-loop fibration, we obtain a graded associative version of the Chas-Sullivan product. The loop product as defined in \cite[Definition 1.2]{hingston2022product} is pulled back from the intersection product with sign convention as in \cite{cieliebak2022topological}. Therefore the previous proposition shows that $H_*(\operatorname{CS}(\kappa_n\otimes m_2))$ recovers the sign convention of the loop product.
\end{ex}

\appendix

\section{Comparison of DG Morse homology conventions}\label{appendix: The BDHO toolset in our convention}

In Section \ref{sec:Morse hom with DG}, we define, unlike in \cite{BDHO23}, the Barraud-Cornea twisting cocycle as a family of chains of paths $m_{x,y} \in C_*(\mathcal{P}_{x \to y} X)$ for $x,y \in \Crit(f)$. Barraud, Damian, Humilière and Oancea define it as a family of chains $m_{x,y} \in C_*(\Omega X)$ over the differential graded algebra of chains on the Moore based loop space of $X$ with the Pontryagin product. The convention we use in this article was also used by Stegemeyer in \cite[\S6]{ste26}.

The construction and manipulation are easier since it does not rely on a basepoint nor on the collapse of a tree $\mathcal{Y}$ connecting all critical points to the basepoint. The drawback is that the coefficients considered no longer have a $C_*(\Omega X)$-module structure but a module over the path groupoid as defined in Definition \ref{def:path groupoid}.

This complex differs from the complex used in \cite{BDHO23} mostly by a matter of conventions on the parametrisation maps. The goal of this section of the appendix is to translate some key constructions of \cite{BDHO23} and \cite{Rie24} in this context and compare the two conventions.

Let $\Xi = (f,\xi,o,s_{x,y})$ be a set of DG Morse data (see Definition \ref{defi : DG Morse data}) and let $\fibration$ be a fibration over $X$ together with a transitive lifting function $\Phi \colon E \ftimes{\pi}{\ev_0} \mathcal{P}X \to E$. For $a \in X$, we denote $F_a = \pi^{-1}(\{a\}).$

\subsection{The Fibration Theorem}

We prove that the complex $C_*(X,C_*(F))$ computes the homology $H_*(E)$ of the total space $E$.

\begin{thm}\label{thm : fibration theorem for our complex}
    There exists a quasi-isomorphism $\Psi_E\colon C_*(X ,C_*(F)) \to C_*(E)$.
\end{thm}

\begin{proof}
    Recall that the set of DG Morse data $\Xi = (f,\xi, o,s_{x,y})$ is the data of a Morse-Smale pair $(f,\xi)$ on $X$, a choice of orientation $[\wb^u(x)]$ of the Latour cell $\wb^u(x)$ for each $x \in \Crit(f)$ and a representing chain system $(s_{x,y})_{x,y \in \Crit(f)}$ for the Morse moduli spaces in $X$, see Definition \ref{defi : representing chain system traj}.
    
    We build $\Psi_E$ in a similar fashion as the quasi-isomorphism of the Fibration Theorem \cite[Theorem 7.2.1]{BDHO23} except that we do not rely on a basepoint $\star \in X$ nor on the collapsing of a tree connecting the critical points of $f$ with the basepoint $\star$.

\begin{defi}[{\cite[Definition 7.3.1]{BDHO23}}]\label{def: repr chains Latour cells}
    A representing chain system of the Latour cells in $X$ compatible with $(s_{x,y})_{x,y \in \Crit(f)}$ is a collection $\{s_x \in C_{|x|}(\wb^u(x)) \ \lvert \  x \in \textup{Crit}(f) \}$ of chains such that:

    \begin{enumerate}
        \item each $s_x$ is a cycle relative to the boundary and represents the fundamental class $[\wb^u(x), \partial \wb^u(x)]$;\\
        \item each $s_x$ satisfies $\partial s_x = \displaystyle \sum_y s_{x,y} \times s_y$ with the product of chains defined via the inclusions $\trajb{x,y} \times \overline{W}^u(y) \subset \partial \overline{W}^u(x) \subset \overline{W}^u(x).$ 
    \end{enumerate}

Such a system always exists and is constructed inductively.
\end{defi}

    For each $x \in \Crit(f)$, consider the evaluation map $\Gamma_x\colon \wb^u(x) \to \mathcal{P}_{x\to X} X$ that parametrizes the gradient lines using the Morse function $f$. See \cite[Lemma 7.3.3]{BDHO23} for more details on this parametrisation. In particular, it satisfies 
    \begin{itemize}
         \item[a)] For any $(\lambda, a) \in \trajb{x,y} \times \wb^u(y) \subset \partial \wb^u(x)$, we have $$\Gamma_x(\lambda,a) = \Gamma_{x,y}(\lambda) \cdot \Gamma_y(a),$$ where $\Gamma_{x,y}\colon \trajb{x,y} \to \mathcal{P}_{x \to y} X$ is the parametrisation map by the values of $f$ described in Definition \ref{def : parametrisation map}.
         \item[b)] For any $x \in \Crit(f)$,  $\ev_1 \circ \Gamma_x = i_x$, where $i_x\colon  \wb^u(x) \to X$ is the attachment map associated with the cellular decomposition of $X$ by the compactified unstable manifolds (also called Latour cells \cite{AD14}). 
     \end{itemize}

    Define $m_x = \Gamma_{x,*}(s_x) \in C_{|x|}(\mathcal{P}_{x\to X} X)$ for $x \in \Crit(f)$. This family satisfies the equation

    \begin{equation}
        \partial m_x = \sum_y m_{x,y} \cdot m_y.
    \end{equation}
    
    Consider $$\deffct{\Psi_E}{C_*(X  ,C_*(F))}{C_*(E)}{x\otimes \alpha}{\Phi^{op}_*(m_x \otimes \alpha)}.$$

We prove that $\Psi_E$ is a chain map.
Let $x \in \Crit(f)$ and $\alpha \in C_*(F_x)$. 

\begin{align*}
    \partial\Psi_E(x \otimes \alpha) &= \sum_y \Phi^{op}_*(m^{CR}_{x,y} \cdot m_y \otimes \alpha) + (-1)^{|x|} \Phi^{op}_*(m_x \otimes \partial_F\alpha) \\
    &= \sum_y (-1)^{(|x|-|y|-1)|y|} \Phi^{op}_*(m_y \otimes \Phi^{op}_*(m^{CR}_{x,y} \otimes \alpha)) + (-1)^{|x|} \Phi^{op}_*(m_x \otimes \partial_F\alpha) \\
    &= \sum_y(-1)^{|x||y|} \Psi_E(y \otimes \Phi^{op}_*(m^{CR}_{x,y} \otimes \alpha)) + (-1)^{|x|} \Psi_E(x \otimes \partial_F \alpha) \\
    &= \Psi_E(\partial(x\otimes \alpha)).
\end{align*}

The rest of the proof is the same as in \cite[Theorem 7.2.1]{BDHO23}. We still spell out the arguments in our setting.

The complex $$C_*(X  ,C_*(F)) = \bigoplus_{x\in \Crit(f)} \Z\cdot x \otimes C_*(F_x)$$ is endowed with the natural filtration $$F_p\left(C_k(X  ,C_*(F))\right) = \bigoplus_{\substack{|x|+j = k \\ |x|\leq p}} \Z\cdot x \otimes C_j(F_x)$$

which induces a spectral sequence $E^r_{p,q}$ that converges to $H_*(X ,C_*(F))$ and we prove that $\Psi_E$ induces a morphism of spectral sequences between $E^r_{p,q}$
and the Leray-Serre spectral sequence $\mathcal{E}^r_{p,q}$ associated with the filtration $$F_p(C_*(E)) = C_*(\pi^{-1}(\Sk_p))$$ of $C_*(E)$, where $\Sk_p$ is the union of cells $i_x(\wb^u(x))$ with $|x| \leq p$.

Indeed, if $|x| \leq p$ and $\alpha \in C_*(F_x)$, then
$$\pi_*(\Psi_E(x \otimes \alpha)) = (\pi \circ \Phi^{op})_*(m_x \otimes \alpha) = \ev_{1,*}(m_x) = (\ev_1 \circ \Gamma_x)_*(s_x) = i_{x,*}(s_x) \in C_*(\Sk_p).$$

It now suffices to prove that $\Psi_E$ induces an isomorphism between the first pages $$\Psi_E^{(1)}\colon  \bigoplus_{x \in \Crit_p(f)} \Z \cdot x \otimes H_q(F_x) = E^1_{p,q} \overset{\simeq}{\to} \mathcal{E}^1_{p,q} = H_{p+q}(\pi^{-1}(\Sk_p),\pi^{-1}(\Sk_{p-1})).$$

For $p \in \Z_{\geq 0}$, denote $\Crit_p(f) = \{ x \in \Crit(f), \ |x|=p\}$.

\begin{lemme}
Let $x \in \Crit_p(f)$.
For all $a \in \wb^u(x)$, denote $\gamma_a = \Gamma_x(a)$ the parametrized gradient line going from $x$ to $a$.

    The maps $$\deffct{\chi_x}{\wb^u(x) \times F_x}{\pi^{-1}(\Sk_p)}{(a, v)}{\Phi^{op}(\gamma_a,v)}$$ induce an isomorphism $$\chi\colon \bigoplus_{x \in \Crit_p(f)} H_p(\wb^u(x), \partial \wb^u(x)) \otimes H_q(F) \overset{\simeq}{\to} \mathcal{E}^1_{p,q}.$$ 
\end{lemme}

To conclude the proof of Theorem \ref{thm : fibration theorem for our complex}, we remark that $\chi$ induces the same map in homology as $\Psi_E^{(1)}$. 

Indeed, at the chain-level,
$$\chi_{x,*}(s_x \otimes \alpha) = \Phi^{op}_*(\Gamma_{x,*}(s_x) \otimes \alpha) = \Psi_E(x \otimes \alpha).$$

Therefore, $\chi$ induces the same map as $\Psi_E^{(1)}$ through the bijection $$\bigoplus_{|x|=p} H_p(\wb^u(x), \partial \wb^u(x)) \simeq \Z \langle\Crit_p(f)\rangle.$$
\end{proof}

\subsection{\texorpdfstring{$\PG$}{PG}-modules and morphisms of fibrations}\label{subsec:PG-modules and morphisms of fibrations}

Let $F_0 \hookrightarrow E_0 \overset{\pi_0}{\to} X$ and $F_1 \hookrightarrow E_1  \overset{\pi_1}{\to} X$ be fibrations over $X$. Let $\varphi\colon E_0 \to E_1$ be a continuous map that preserves the fibrations, i.e., $\pi_1 \circ \varphi = \pi_0$. Such a map is called a morphism of fibrations.
Let $\Xi = (f,\xi,o,s_{x,y})$ be a set of DG Morse data.

One goal of this subsection is to extract a chain map $\widetilde{\varphi}\colon C_*(X,C_*(F_0)) \to C_*(X,C_*(F_1))$ from the morphism of fibrations $\varphi$.  

As shown in \cite[\S 5]{Rie24}, the construction of such a map in the BDHO setting requires the notion of $\Ai$-morphism of topological modules over $\Omega X$. 
The natural replacement of $\Omega X$ in our setting is the path groupoid $\PG$ of $X$ (Definition \ref{def:path groupoid}).

We then define the generalised Morse homology with coefficients in a module $A =((A_x)_{x \in \Crit(f)},\mu_A)$ over $C_*(\mathcal{C}_f)$, the chain complex associated to the flow category $\mathcal{C}_f$ of $f$. This construction encompasses coefficients that are modules over $C_*(\PG)$ using the functor $\Gamma\colon \mathcal{C}_f\to \PG$ induced by the parametrisation maps.

This generalises the construction of \S \ref{sec:Morse hom with DG} which corresponds to the $\PG$-module $C_*(F) \coloneqq(C_*(F_x), \Phi^{op}_*)$ associated to a fibration $E \overset{\pi}{\to} X$, where $F_x = \pi^{-1}(\{x\})$ for $x \in \Crit(f).$

We then replace the notion of an $\Ai$-morphism of topological modules over $\Omega X$ by the notion of an $\Ai$-morphism of $\mathcal{C}_f$-modules (Definition \ref{def: Ai morphism of Cf modules}). We show that the arguments provided for the proof of \cite[Theorem 5.8]{Rie24} apply in this setting to prove that a morphism of fibrations $\varphi\colon E_0 \to E_1$ induces an $\Ai$-morphism of $\mathcal{C}_f$-modules and thus a chain map $\widetilde{\varphi}\colon  C_*(X,C_*(F_0)) \to C_*(X ,C_*(F_1))$ that corresponds to $\varphi_*\colon C_*(E_0) \to C_*(E_1)$ via the Fibration Theorem.

\begin{defi}[\cite{cohen1995morse}]
    The \textit{flow category} of a Morse-Smale pair $(f,\xi)$, $\mathcal{C}_f$, is given by the topologically enriched category with objects $x\in \Crit(f)$ and morphism spaces 
    \begin{align*}
        \mathcal{C}_f(x,y)=\begin{cases}
            \overline{\mathcal{L}}(x,y) &\text{if } x\neq y;\\
            \{\id_x\} & \text{if }x=y;
        \end{cases}
    \end{align*}
    and composition is given by concatenation of trajectories.
\end{defi}

\begin{defi}\label{def:path groupoid}
    The \textit{path groupoid of $X$}, $\PG$, is given by the topologically enriched category with objects $x\in \Crit(f)$ and morphism spaces 
    \begin{align*}
        \PG(x,y)=\mathcal{P}_{x\to y} X
    \end{align*}
    and concatenation of paths.
\end{defi}

\begin{defi}
    A \textit{module over $\mathcal{C}_f$ (or $\PG$)} is given by a functor from $\mathcal{C}_f$ (or $\PG$, respectively) into the category of topological spaces, that is, a tuple of topological spaces $(A_x)_{x\in \Crit(f)}$ with maps
    \begin{align*}
        \mu_A\colon \mathcal{C}_f(x,y)\times A_x&\to A_y,\\
        \text{(or }\nu_A\colon \PG(x,y)\times A_x&\to A_y\text{, respectively)}
    \end{align*}
    that commute with composition.
\end{defi}

\begin{lemme}
    The evaluation maps $\Gamma_{x,y}\colon \overline{\mathcal{L}}(x,y)\to \mathcal{P}_{x\to y}X$ define a functor $\Gamma\colon \mathcal{C}_f\to \PG$. 

    In particular, precomposition with $\Gamma$ gives a functor from $\PG$-modules to $\mathcal{C}_f$-modules.
\end{lemme}
\begin{proof}
    We define $\Gamma(\id_x)=c_x$, the constant path at $x$ of length $0$. The fact that $\Gamma_{x,y}$ commutes with composition implies that $\Gamma$ is a functor.
\end{proof}

\begin{rem}
    After taking chains, we get modules over the differential graded category $C_*(\mathcal{C}_f)$. Modules over this category play a role analogous to $C_*(\Omega X)$-modules in \cite{BDHO23}. In fact $(\theta\circ p \circ \Gamma)_*\colon C_*(\mathcal{C}_f)\to C_*(\Omega X)$ commutes with composition, making any $C_*(\Omega X)$-module into a $C_*(\mathcal{C}_f)$-module (see Example \ref{ex:BDHO as PG module}).

    On the other hand, $\mathcal{C}_f$ is more flexible because there is no requirement that $A_x\simeq A_y$ for different critical points $x,y\in \Crit(f)$. In fact by work of Fourel in \cite{fourel2026morse}, $\mathcal{C}_f$ models the exit path category associated to the stable manifolds and $\mathcal{C}_f$-modules model constructible sheaves.
\end{rem}

\begin{defi}\label{defi:generalised Morse homology with DG coeffs}
    Let $A\coloneqq((A_x)_{x\in \Crit(f)},\mu_A)$ be a $C_*(\mathcal{C}_f)$-module. The \textit{generalised Morse complex with coefficients in $(A_x)_{x\in \Crit(f)}$} is given by
    \begin{align*}
        C_*(X ,A)\coloneqq \bigoplus_{x\in \Crit(f)} \Z\cdot x\otimes A_x
    \end{align*}
    and differential 
    \begin{align*}
        \partial(x \otimes \alpha) = \sum_{y \in \Crit(f)} (-1)^{|x||y|} y \otimes \mu_A(s_{x,y}\otimes \alpha) + (-1)^{|x|} x \otimes \partial(\alpha)
    \end{align*}
    for $s_{x,y}$ a representing chain system as in Definition \ref{defi : representing chain system traj}.
\end{defi}

\begin{ex}\label{ex:BDHO as PG module}
    The collection $(F_x)_{x\in \Crit(f)}$ together with $\nu_A\coloneqq\Phi^{op}$ defines a $\PG$-module. Therefore $\mu_A\coloneqq\Phi^{op}\circ(\Gamma\times \id)$ defines a $\mathcal{C}_f$-module structure. Because we defined $m_{x,y}^{CR}\coloneqq\Gamma_*(s_{x,y})$, the associated Morse complex gives $C^{CR}_*(X  ,C_*(F))$.

    The constant collection $(F_\star)_{x\in \Crit(f)}$ with $\mu_A(\lambda,\alpha)\coloneqq\Phi^{op}(\theta\circ p\circ\Gamma(\lambda),\alpha)$ defines a $\mathcal{C}_f$-module.  
\end{ex}

\begin{notation}
    For critical points $x,y\in \Crit(f)$, we denote for $k\geq 1$
    \begin{align*}
        \mathcal{C}^k_f(x,y)\coloneqq \bigcup_{z_1,\dots, z_{k-1} \in \Crit(f)} \mathcal{C}_f(z_{k-1},y)\times \dots \times \mathcal{C}_f({z_1, z_2})\times \mathcal{C}_f(x, z_1)
    \end{align*}
    and 
    \begin{align*}
        \mathcal{C}_f^0(x,y)=\begin{cases}
            *&\text{if }x=y;\\
            \varnothing &\text{else.}
        \end{cases}
    \end{align*}

    The topology on $\mathcal{C}^k_f(x,y)$ is the topology induced by $\mathcal{C}_f(x,y).$ 
\end{notation}

For $\mathcal{C}_f$-modules, we can define $\Ai$-morphisms akin to $\Ai$-morphisms of $\Omega X$-modules as in \cite[Definition 5.3]{Rie24}.

\begin{defi}\label{def: Ai morphism of Cf modules}
    Let $((A_x^1)_{x\in \Crit(f)},\mu_1)$ and $((A^2_x)_{x\in \Crit(f)},\mu_2)$ be $\mathcal{C}_f$-modules. An \textit{$\Ai$-morphism of $\mathcal{C}_f$-modules} is given by a sequence of continuous maps for each $x,y \in \Crit(f)$, for all $k\geq 0$,
    \begin{align*}
        \varphi_{k+1}^{x\to y}\colon I^k \times \mathcal{C}_f^k(x,y)\times A_x^1\to A^2_y
    \end{align*}
    such that for all $t_1,\dots, t_k\in I^k$, $\alpha\in A_x$ and $(\tau_k,\dots,\tau_1)\in \mathcal{C}_f(z_{k-1}, y)\times  \dots \times \mathcal{C}_f(x, z_1)$, it holds:
    \begin{align*}
        \varphi_{k+1}^{x\to y}(t_1,\dots,t_k,\tau_k,\dots,\tau_1,\alpha)=\begin{cases}
            \varphi_k^{\ev_1(\tau_1)\to y}(t_2,\dots,t_k,\tau_k,\dots,\tau_2,\mu_1(\tau_1,\alpha))&\text{if } t_1=1;\\
            \varphi_k^{x\to y}(t_1,\dots, \hat{t_j},\dots,t_k,\tau_k,\dots, \tau_{j-1}\cdot\tau_j,\dots,\tau_1,\alpha)&\text{if } t_j=1, j\geq 2;\\
            \mu_2(\tau_j \cdot \dots \cdot \tau_k,\varphi_j^{x\to \ev_0(\tau_j)}(t_1,\dots,t_{j-1},\tau_{j-1},\dots, \tau_{1},\alpha)) &\text{if } t_j=0.
        \end{cases}
    \end{align*}
\end{defi}

The relevance of this structure is that they induce morphisms on Morse homology:

\begin{lemme}\label{lem:Ai PC morph gives Morse morph}
    An $\Ai$-morphism of $\mathcal{C}_f$-modules, $\varphi\colon A^1\to A^2$, induces a chain map $\widetilde{\varphi}\colon C_*(X  ,C_*(A^1))\to C_*(X  ,C_*(A^2))$.
\end{lemme}
\begin{proof}
    The proof and formula are analogous to the case of $\Ai$-morphisms over $\Omega X$-modules as in \cite[Corollary 5.6]{Rie24}. 

        For $x \in \Crit(f),$ we use the notation $\mathbf{m^x} \in \Hom_{-1}(\Z\cdot x; \bigoplus_y  \Z\cdot y \otimes C_*(\mathcal{C}_f(x,y)))$ for the map
$\mathbf{m}(x) = \sum_y (-1)^{|x||y|} y \otimes s_{x,y}$. 

For $x \in \Crit(f)$, denote $$TC_*(\mathcal{C}_f)(x) = \bigoplus_{\substack{z \in \Crit(f)\\ k \geq 1}}C_*(\mathcal{C}_f^k(z,x)) \otimes C_*(A^1_z) $$ and $$A_{\mathcal{C}_f} = \bigoplus_{x \in \Crit(f)} \Z \cdot x \otimes TC_*(\mathcal{C}_f)(x). $$

We extend $\mathbf{m^x}$ to $$\mathbf{\widetilde{m}}\in \End_{-1}\left( A_{\mathcal{C}_f} \right)$$ defined by
$$\mathbf{\widetilde{m}}(x \otimes \sigma_k \otimes \dots \otimes \sigma_1 \otimes \alpha) = \sum_y (-1)^{|x||y|} y \otimes s_{x,y} \otimes \sigma_k \otimes \dots \otimes \sigma_1  \otimes \alpha.$$
Let $\left\{\varphi_{k+1}^{x \to y}\colon I^n \times \mathcal{C}_f^{k}(x,y) \times A^1_x \to A^2_y\right\}$ be the $\Ai$-morphism of $\mathcal{C}_f$-modules.
For all $k\geq 0$, define $$\phi_{k+1}\colon C_*(\mathcal{C}_f^k(x,y)) \otimes C_*(A^1_x)\to C_{k+*}(A^2_y)$$ by $$\phi_{k+1}( \sigma_k \otimes \dots \otimes \sigma_1 \otimes \alpha) = \varphi^{x \to y}_{k+1,*}(\Id_{I^k} \otimes \sigma_k \otimes \dots \otimes \sigma_1 \otimes \alpha).$$

the map $\widetilde{\varphi}$ is then defined by $$\widetilde{\varphi} = \sum_{k\geq 0}(1 \otimes \phi_{k+1}) \m^{k}. $$
\end{proof}

\begin{notation}
    For $a,b \in X$ and an integer $k \geq 1$, denote $$\PG^k(a,b) = \bigcup_{z_1, \dots, z_{k-1} \in X} \mathcal{P}_{z_{k-1} \to b} X  \times \dots \times \mathcal{P}_{z_1 \to z_2} X \times \mathcal{P}_{a \to z_1} X $$ and 
    \begin{align*}
        \PG^0(a,b)=\begin{cases}
            c_a&\text{if }a=b;\\
            \varnothing &\text{else,}
        \end{cases}
    \end{align*}
    where $c_a$ is the constant path at $a$.
The topology on $\PG^k(a,b)$ is the topology induced by $\mathcal{P}_{a \to b} X$ and the concatenation $\PG^k(a,b) \to \mathcal{P}_{a \to b}X$.

Denote also $$\PG^k= \bigcup_{a,b} \PG^k(a,b)$$ with the topology induced by the concatenation $\PG^k \to \mathcal{P}X$ and if $\tau =(\tau_k,\dots,\tau_1) \in \PG^k(a,b)$, denote $\ev_1(\tau) = \ev_1(\tau_k)=b$ and $\ev_0(\tau) = \ev_0(\tau_1)=a$.
\end{notation}

\begin{lemme}\label{lem: morphism of fibrations formula with paths}
    Let $F^1 \hookrightarrow E^1 \overset{\pi_1}{\to} X$ and $F^2 \hookrightarrow E^2 \overset{\pi_2}{\to} X$ be two fibrations over $X$ and $\Phi_1$, $\Phi_2$ be two transitive lifting functions associated with $E^1$ and $E_2$, respectively. Let $\varphi\colon E_1 \to E_2$ be a morphism of fibrations, i.e., $\varphi$ is continuous and $\pi_2 \circ\varphi = \pi_1$. Denote $\mu_1 \coloneqq \Phi_1^{op}$, $\mu_2 \coloneqq \Phi_2^{op}$ and $\varphi_1 \coloneqq \varphi$.
    
    There exists a sequence of continuous maps for $k\geq 1$,
    \begin{align*}
        \varphi_{k+1}\colon I^k\times (\PG^{k} \ftimes{\ev_1}{\pi_1} E^1) \to E^2
    \end{align*}
     such that for all $t_1,\dots, t_k\in I^k$, $e\in E^1$ and $(\tau_k,\dots,\tau_1)\in \PG^{k}$, it holds that $\pi_2\circ \varphi_{k+1}(t_1,\dots,t_k,\tau_k,\dots,\tau_1,e) = \ev_1(\tau_k)$ and

        \begin{equation}\label{eq: Ai PG module induced by morphism of fibrations}
        \varphi_{k+1}(t_1,\dots,t_k,\tau_k,\dots,\tau_1,e)=\begin{cases}
            \varphi_k(t_2,\dots,t_k,\tau_k,\dots,\tau_2,\mu_1(\tau_1,e))&\text{if } t_1=1;\\
            \varphi_k(t_1,\dots, \hat{t_j},\dots,t_k,\tau_k,\dots, \tau_{j-1}\cdot\tau_j,\dots,\tau_1,e)&\text{if } t_j=1, j\geq 2;\\
            \mu_2(\tau_j \cdot \dots \cdot \tau_k,\varphi_j(t_1,\dots,t_{j-1},\tau_{j-1},\dots, \tau_{1},e)) &\text{if } t_j=0.
        \end{cases}
    \end{equation}
\end{lemme}

\begin{proof}
    The proof is analogous to the proof of \cite[Lemma 5.9]{Rie24}
\end{proof}

\begin{cor}\label{cor: morphism of fibrations induces Ai Cf module}
    A morphism of fibrations $\varphi\colon E_1 \to E_2$ induces a chain map $\widetilde{\varphi}\colon C_*(X  ,C_*(F^1))\to C_*(X  ,C_*(F^2))$.
\end{cor}
\begin{proof}
     Using the parametrisation $\Gamma\colon \mathcal{C}_f \to \PG$, it follows from Lemma \ref{lem: morphism of fibrations formula with paths} that a morphism of fibrations induces a morphism of $\mathcal{C}_f$-modules $(\varphi_{k+1}^{x\to y})$, where $\varphi_1^{x \to x} = \varphi_{\lvert F^1_{x}}$. Lemma \ref{lem:Ai PC morph gives Morse morph} finishes the proof.
\end{proof}

We obtain the following version of \cite[Theorem 5.8]{Rie24}.

\begin{thm}\label{thm:morphism of fibrations}
    Let $\varphi\colon E_1 \to E_2$ be a morphism of fibrations. Denote $\Psi_1: C_*(X  ,C_*(F^1)) \to C_*(E_1)$ and $\Psi_2\colon C_*(X  ,C_*(F^2)) \to C_*(E_2) $ the quasi-isomorphisms given by the Fibration Theorem \ref{thm : fibration theorem for our complex}. 
    
    The following diagram commutes up to chain homotopy

    \[
    \xymatrix{
    C_*(X  ,C_*(F^1)) \ar[r]^-{\widetilde{\varphi}} \ar[d]^{\Psi_1} & C_*(X  ,C_*(F^2)) \ar[d]_{\Psi_2} \\
    C_*(E_1) \ar[r]_{\varphi_*} & C_*(E_2).
    }
    \]
\end{thm}

\begin{proof}
This proof is computationally heavy and requires a lot of notations. However, it follows the proof of \cite[Theorem 5.8]{Rie24}. We still make the computations explicit, since our change of conventions induces different notations and a few changes of signs.

    For $x \in \Crit(f),$ we use the notation $\mathbf{m^x} \in \Hom_{-1}(\Z\cdot x; \bigoplus_{y \in \Crit(f)}  \Z\cdot y \otimes C_*(\PG(x,y)))$ for the map
$\mathbf{m}(x) \coloneqq \sum_y (-1)^{|x||y|} y \otimes m_{x,y}$. 

For $x \in \Crit(f)$, denote $$TC_*(\PG)(x) = \bigoplus_{\substack{z \in \Crit(f)\\ k \geq 1}}C_*(\PG^k(z,x)) \otimes C_*(F^1_z) $$ and $$A_{\PG} = \bigoplus_{x \in \Crit(f)} \Z \cdot x \otimes TC_*(\PG)(x). $$

We extend $\mathbf{m^x}$ by $$\mathbf{\widetilde{m}}\in \End_{-1}\left( A_{\PG} \right)$$ defined by
$$\mathbf{\widetilde{m}}(x \otimes \sigma_k \otimes \dots \otimes \sigma_1 \otimes \alpha) = \sum_y (-1)^{|x||y|} y \otimes m_{x,y} \otimes \sigma_k \otimes \dots \otimes \sigma_1  \otimes \alpha.$$

We use the representing chain system $\left\{s_x \in C_{|x|}\left(\wb^u(x)\right)\right\}_{x \in \Crit(f)}$  of the Latour cell $\wb^u(x)$  (see Definition \ref{def: repr chains Latour cells}). Consider for $x \in \Crit(f)$ the evaluation maps $\Gamma_x\colon \wb^u(x) \to \mathcal{P}_{x \to X}X$ by the values of $f$ and define $m_x = \Gamma_{x,*}(s_x) \in C_{|x|}\left(\mathcal{P}_{x \to X} X\right).$ It satisfies the equation

\begin{equation}\label{eq : m_x}
    \partial m_x = \sum_{y \in \Crit(f)} m_{x,y}\cdot  m_y.
\end{equation}

We denote $\mathbf{m^{L,x}} \in \Hom_0\left(\Z\cdot x, C_*\left(\mathcal{P}_{x \to X}X\right)\right)$ the morphism defined by $$\mathbf{m^{L,x}}(x) = m_x$$ for all $x \in \textup{Crit}(f)$ and we extend it to 
$$\mathbf{\widetilde{m}^L} \in \Hom_0\left(A_{\PG} , \bigoplus_{x \in \Crit(f)}C_*(\mathcal{P}_{x \to X}X) \otimes TC_*(\PG)(x) \right)$$
by
$$\mathbf{\widetilde{m}^L}(x \otimes \sigma_k \otimes \dots \otimes \sigma_1 \otimes \alpha) = m_x \otimes \sigma_k \otimes \dots \otimes \sigma_1 \otimes \alpha.$$

Then, for $i \in \{1,2\}$, the definition of the quasi-isomorphism $\Psi_i\colon C_*(X  ,C_*(F^i)) \to C_*(E_i)$ can be written
\begin{equation}\label{eq : preuve theorem C def iso fibration thm}
    \Psi_i = \Phi^{op}_{i,*} \mathbf{\widetilde{m}^L}.
\end{equation}

Let $\left\{\varphi_{k+1}\colon I^k \times (\PG^{k} \ftimes{\ev_0}{\pi} E^1) \to E^2\right\}$ be the sequence of maps given by Lemma \ref{lem: morphism of fibrations formula with paths}.
We denote $\phi_1 = \varphi_{*}\colon C_*(E_1) \to C_*(E_2)$ and for all $k\geq 1$, define $$\phi_{k+1}\colon C_*(\PG^k \ftimes{\ev_1}{\pi} E^1)\to C_*(E^2)$$ by $$\phi_{k+1}( \sigma_k \otimes \dots \otimes \sigma_1 \otimes \alpha) = \varphi_{k+1,*}(\Id_{I^n} \otimes \sigma_k \otimes \dots \otimes \sigma_1 \otimes \alpha).$$

In particular $$\widetilde{\varphi} = \sum_{k\geq 0}(1 \otimes \phi_{k+1}) \m^{k}. $$

\begin{notation}
    \begin{itemize}
        \item Denote $\nu_2\colon C_*(\PG^2(x,a)) \to C_*(\mathcal{P}_{x \to a} X)$, $\nu_2(\sigma_2 \otimes \sigma_1) = (-1)^{|\sigma_1||\sigma_2|} \sigma_1 \cdot \sigma_2 $.
        \item Denote $\nu_1$ the differential of $C_*(\mathcal{P}X).$
    \end{itemize}
    
\end{notation}

Equation \eqref{eq: Ai PG module induced by morphism of fibrations} at the chain-level gives

\begin{equation}\label{eq : preuve Theorem C}
\begin{split}
    \partial \phi_{N+1} &= \phi_{N}( 1^{\otimes N-1} \otimes \Phi_{1,*}^{op} ) + \sum_{r=1}^{N-1} (-1)^r \phi_{N}(1^{\otimes N-r-1} \otimes \nu_2 \otimes 1^{\otimes r})+ (-1)^{N}\Phi_{2,*}^{op}(1 \otimes \phi_{N}) \\
    &+   \sum_{r=1}^{N} (-1)^N \phi_{N+1}(1^{\otimes N-r} \otimes \nu_1 \otimes 1^{\otimes r}) + (-1)^{N} \phi_{N+1} (1^{\otimes N} \otimes \partial).
\end{split}
\end{equation}

Define
$$v = \sum_{n \geq 1} (-1)^{n+1} \phi_{n+1} \mL \m^{n-1}\colon C_*(X  ,C_*(F^1)) \to C_{1+*}(E_2).$$ 

We prove that $$\varphi_*\Psi_1 - \Psi_2 \widetilde{\varphi} = \partial v + v\partial.$$

Let us first remark that the equation \eqref{eq : m_x} can be written \begin{equation}\label{eq : sys repr cell Ai}
    \nu_1 \mL = \nu_2\mathbf{m^L} \m.
\end{equation}

We now use Equation \eqref{eq : preuve Theorem C} to compute

\begin{align*}
    \partial v &=  \sum_{k \geq 1} (-1)^{k+1} \partial \phi_{k+1} \mL \m^{k-1} \\
    &\overset{\eqref{eq : preuve Theorem C}}{=}\sum_{k\geq 1} (-1)^{k+1} \phi_{k}( 1^{\otimes k-1} \otimes \Phi_{1,*}^{op} ) \mL \m^{k-1} + \sum_{k\geq 1} (-1)^{k+1} \sum_{r=1}^{k-1} (-1)^r \phi_{k}(1^{\otimes k-r-1} \otimes \nu_2 \otimes 1^{\otimes r})  \mL \m^{k-1} \\
    &+ \sum_{k\geq 1} (-1)^{k+1} (-1)^{k}\Phi_{2,*}^{op}(1 \otimes \phi_{k}) \mL \m^{k-1} +   \sum_{k\geq 1} (-1)^{k+1}\sum_{r=1}^{k} (-1)^k \phi_{k+1}(1^{\otimes k-r} \otimes \nu_1 \otimes 1^{\otimes r}) \mL \m^{k-1} \\
    &+ \sum_{k\geq 1} (-1)^{k+1}(-1)^{k} \phi_{k+1} (1^{\otimes k} \otimes \partial) \mL \m^{k-1}.
\end{align*}

To continue the computation, we perform the following operations:

\begin{enumerate}
    \item Split \begin{align*}
        \sum_{k\geq 1} (-1)^{k+1} \phi_{k}( 1^{\otimes k-1} \otimes \Phi_{1,*}^{op} ) \mL \m^{k-1} &= \underbrace{\phi_1}_{= \varphi_*} \underbrace{\Phi_{1,*}^{op} \mL}_{\overset{\eqref{eq : preuve theorem C def iso fibration thm}}{=}{\Psi_1}} + \sum_{k\geq 2} (-1)^{k+1} \phi_k ( 1^{\otimes k-1} \otimes \Phi_{1,*}^{op}) \mL \m^{k-1}\\
        &= \varphi_* \Psi_1 - \sum_{k \geq 1} (-1)^{k+1} \phi_{k+1}( 1^{\otimes k} \otimes \Phi_{1,*}^{op})\mL \m^{k}
    \end{align*}
    \item Notice that $\Phi_{2,*}(1 \otimes \phi_k) \mL \m^{k-1} = \Phi_{2,*}^{op} \mL (1 \otimes \phi_k) \m^{k-1} $ for all $k \geq 1$ since $\phi_k$ only takes the last $k$ arguments. By definition, $$\sum_{k \geq 1}\Phi_{2,*}^{op} \mL (1 \otimes \phi_k) \m^{k-1} = \Phi_{2,*} \mL\left( \sum_{k \geq 1} (1 \otimes \phi_k ) \m^{k-1}\right) = \Psi_2 \widetilde{\varphi}.$$
    \item Notice that $$\phi_{k+1} (1^{\otimes k} \otimes \partial) \mL \m^{k-1} = - (-1)^k \phi_{k+1}\mL \m^{k-1} (1 \otimes \partial)$$ for all $k \geq 1.$
\end{enumerate}

We thus obtain 

\begin{align*}
    \partial v &= \varphi_*\Psi_1 - \sum_{k \geq 1} (-1)^{k+1} \phi_{k+1}(\Phi_{1,*} \otimes 1^{\otimes k})\mL \m^{k} - \Psi_2 \widetilde{\varphi} - \sum_{k \geq 1} (-1)^{k+1} \phi_{k+1} \mL \m^{k-1}(1 \otimes \partial) \\
    &- \sum_{k\geq 1} (-1)^{k} \sum_{r=1}^{k-1} (-1)^r \phi_{k}(1^{\otimes k-r-1} \otimes \nu_2 \otimes 1^{\otimes r})  \mL \m^{k-1} \\
    &-  \sum_{k\geq 1} \sum_{r\geq 1} \phi_{k+1}(1^{\otimes k-r} \otimes \nu_1 \otimes 1^{\otimes r}) \mL \m^{k-1}.
\end{align*}

Using equation \eqref{eq : sys repr cell Ai}, we obtain the following equation: \begin{align*}
    (1^{\otimes k-r} \otimes \nu_1 \otimes 1^{\otimes r}) \mL \m^{k-1} &= (-1)^{r+k+1} \mL \m^{k-r-1}(1 \otimes \nu_1 \otimes 1^{\otimes r})\m^r\\
    &= (-1)^{r+k} \mL \m^{k-r-1}(1 \otimes \nu_2 \otimes 1^{\otimes r})\m^{r+1}\\
    &=  (-1)^{r+k}(1^{\otimes r} \otimes \nu_2 \otimes 1^{\otimes k-r}) \mL \mathbf{\widetilde{m}}^{k+1}.
\end{align*}

 Indeed, a direct computation proves that for all $z,x \in \Crit(f)$, $\sigma_{r-1}\otimes \dots\otimes \sigma_1 \otimes \alpha \in C_*(\PG^{r-1} \ftimes{\ev_1}{\pi} E^1)$,

$$(1 \otimes \nu_1 \otimes 1^{\otimes r})\m(x \otimes \sigma_{r-1} \otimes \dots \otimes \sigma_1 \otimes \alpha) = -(1 \otimes \nu_2 \otimes 1^{\otimes r})\m^2(x \otimes \sigma_{r-1} \otimes \dots \otimes \sigma_1 \otimes \alpha).$$

We conclude the computation of $\partial v$ with the formula

$$\partial v =  \varphi_*\Psi_1 - \sum_{k \geq 1} (-1)^{k+1} \phi_{k+1}(\Phi_{1,*}^{op} \otimes 1^{\otimes k})\mL \m^{k} - \Psi_2 \widetilde{\varphi} - \sum_{k \geq 1} (-1)^{k+1} \phi_{k+1} \mL \m^{k-1}(1 \otimes \partial).$$

Moreover,
\begin{align*}
    v \partial &= v( \Phi^{op}_{1,*} \m + (1 \otimes \partial)) \\
    & = \sum_{n \geq 1} (-1)^{n+1} \phi_{n+1} \mL \m^{n-1} \Phi_{1,*}^{op}\m + \sum_{k \geq 1} (-1)^{k+1} \phi_{k+1} \mL \m^{k-1} (1 \otimes \partial)\\
    &= \sum_{n \geq 1} (-1)^{n+1} \phi_{n+1}(1^{\otimes k} \otimes \Phi^{op}_{1,*} ) \mL \m^{n} + \sum_{n \geq 1} (-1)^{n+1} \phi_{n+1} \mL \m^{n-1} (1 \otimes \partial).
\end{align*}

We therefore get the relation 

$$\partial v + v \partial = \varphi_*\Psi_1 - \Psi_2\widetilde{\varphi}.$$
\end{proof}

\subsection{Künneth map}\label{subsec:Künneth map}

The goal of this subsection is to recall the construction of the Künneth twisting cocycle and the Künneth map. We follow \cite[\S 6]{Rie24} and change only the evaluation maps by not applying $\theta \circ p$.

We consider two  closed, connected, oriented, smooth manifolds $X,Y$ endowed with DG Morse data $\Xi_X = (f,\xi_X,o_X,s_{x,y}^X)$ and $\Xi_Y=(g,\xi_Y,o_Y,s_{x,y}^Y)$, respectively. We denote $(F_x)_{x \in \Crit(f)}$ the fibres of a fibration $ E_X \overset{\pi_X}{\to} X$ and we denote $(G_x)_{x \in \Crit(g)}$ the fibres of a fibration $ E_Y \to Y$. 

If $\Phi_X$ and $\Phi_Y$ are transitive lifting function associated to $ \pi_X\colon E_X \to X$ and $\pi_Y\colon E_Y \to Y$ respectively, then $\Phi_{X \times Y} = (\Phi_X, \Phi_Y)$ is a transitive lifting function associated to the product fibration $E_X \times E_Y \to X \times Y.$

We remark that $f+g\colon X \times Y \to \R$, $(f+g)(x,y) = f(x) + g(y)$ is a Morse function on $X \times Y$ and $\xi_X + \xi_Y$ is an adapted pseudo-gradient.

It is possible to identify the moduli spaces of trajectories of the cartesian product $X\times Y$ with the cartesian product of the moduli spaces of trajectories on $X$ and on $Y$ only in very specific cases. More precisely, 
$$\trajbi{f+g}{(x,y),(x',y')} \ncong \trajbi{f}{x,x'} \times \trajbi{g}{y,y'}$$ except if $x=x'$ or $y=y'$. There exists a projection $$\pi\colon \trajbi{f+g}{(x,y),(x',y')} \to \trajbi{f}{x,x'} \times \trajbi{g}{y,y'}.$$

\begin{lemme}\cite[Lemma 6.2]{Rie24}\label{lemme : representing chain system of a Cartesian product}
    Let $\{s^X_{x,x'}\}$ and $\{s^Y_{y,y'}\}$ be representing chain systems for the moduli spaces of trajectories $\trajbi{f}{x,x'}$ and $\trajbi{g}{y,y'}$ respectively.
    There exists a representing chain system $\{ s_{(x,y),(x',y')} \}$ for the moduli spaces $\trajbi{f+g}{(x,y),(x',y')}$ such that 

    \begin{enumerate}
        \item $s_{(x,y),(x',y')} = (-1)^{|x|(|y|-|y'|)} (\{x\}, s^Y_{y,y'})$ if $x=x'$.\\
        \item $s_{(x,y),(x',y')} = (s_{x,x'}^X, \{y\})$ if $y=y'$\\
        \item $\pi_{*}(s_{(x,y),(x',y')}) = 0$ if $x \neq x'$ and $y\neq y'$.
    \end{enumerate}
\end{lemme}

\begin{lemme}\cite[Lemma 6.3]{Rie24}\label{lemma : q^K}
The family 
   $$\xymatrix{\Gamma^K_{(x,y),(x',y')}\colon \trajbi{f+g}{(x,y),(x',y')} \ar[r]^-{\pi} & \trajbi{f}{x,x'} \times \trajbi{g}{y,y'} \ar[rr]^{\quad (\Gamma^X_{x,x'}, \Gamma^Y_{y,y'})} & & \mathcal{P}_{x \to x'}X \times \mathcal{P}_{y \to y'}Y }$$ 
    satisfies that if $(\lambda, \lambda') \in \trajb{(x,y),(z,w)} \times \trajb{(z,w),(x',y')} \subset \partial\trajb{(x,y),(x',y')}$, then $$\Gamma^K_{(x,y),(x',y')}(\lambda, \lambda') = \Gamma^K_{(x,y),(z,w)}(\lambda) \cdot \Gamma^K_{(z,w),(x',y')}(\lambda').$$
\end{lemme}

\begin{defi}
    Define the Künneth twisting cocycle $m^K_{(x,y),(x',y')} = \Gamma^K_{(x,y),(x',y'),*}(s_{(x,y),(x',y')})$. This family satisfies $$m^K_{(x,y),(x',y')} = \left\{ \begin{array}{ll}
   (-1)^{|x|(|y|-|y'|)} (c_x, m_{y,y'}) & \textup{if} \ x=x',  \\
     (m_{x,x'}, c_y)& \textup{if} \ y=y', \\
     0 & \textup{otherwise}
\end{array} \right.$$ and $$\partial m^K_{(x,y),(x',y')} = \sum_{(z,w)} (-1)^{|x|+|y|-|z|-|w|} m^K_{(x,y),(z,w)} m^K_{(z,w),(x',y')}.$$ We denote $c_x$ and $c_y$ the constant paths of length 0 at the points $x$ and $y$ respectively seen as $0$-chains.
\end{defi}

\begin{defi}
Denote $C_*(X \times Y, m^K ,C_*(F\times G))$ the complex $$\bigoplus_{(x,y) \in \Crit(f+g)}\Z\cdot(x,y) \otimes C_*(F_x \times G_y)$$ with differential given on the generators by 

\begin{align*}
    \partial((x,y) \otimes (\alpha,\beta)) &= \sum_{(x',y') \in \Crit(f+g)} (-1)^{(|x|+|y|)(|x'|+|y'|)}(x',y') \otimes \Phi_{X,Y,*}^{op}(m^K_{(x,y),(x',y')} \otimes (\alpha, \beta))\\
    &+ (-1)^{|x|+|y|} (x,y) \otimes \partial(\alpha, \beta).
\end{align*}

\end{defi}

\begin{prop}\label{prop : Kunneth map}
    The map $$\deffct{K}{C_*(X ,C_*(F)) \otimes C_*(Y ,C_*(G))}{C_*(X \times Y, m^K ,C_*(F \times G))}{(x \otimes \alpha) \otimes (y\otimes \beta)}{(-1)^{|\alpha||y|} (x,y) \otimes (\alpha, \beta)}$$ is a chain map and corresponds to the Eilenberg-Zilber map via the Fibration Theorem. More precisely, the following diagram commutes

    \[
    \xymatrix{
    H_*(X, C_*(F)) \otimes H_*(Y ,C_*(G)) \ar[r]^-{K} \ar[d]^{\Psi_{E_X} \otimes \Psi_{E_Y}} & H_*(X \times Y, C_*(F \times G)) \ar[d]^{\Psi_{E_{X} \times E_{Y}}} \\
    H_*(E_X) \otimes H_*(E_Y) \ar[r]^-{\EZ} & H_*(E_X \times E_Y),
    }
    \]
    where $\Psi_{E_X}$, $\Psi_{E_Y}$ and $\Psi_{E_X \times E_Y}$ are the isomorphisms given by the Fibration Theorem.
\end{prop}

\subsection{Shriek map}\label{subsec:Shriek map}

In this subsection, we follow \cite[\S 10.4]{BDHO23} and change the evaluation maps by not applying $\theta \circ p$.

We consider  two closed, connected, oriented, smooth manifolds  $X^n,Y^m$ endowed with DG Morse data $\Xi_X = (f,\xi_X,o_X,s_{x,y}^X)$ and $\Xi_Y=(g,\xi_Y,o_Y,s_{x,y}^Y)$, respectively. We denote $(G_x)_{x \in \Crit(g)}$ the fibres of a fibration $G \hookrightarrow E \to Y$. We also denote $\Phi$ a transitive lifting function associated to this fibration.

We consider $\varphi\colon X \to Y$ a continuous map.

\begin{defi}\label{defi: pull back module}
    Denote $\varphi^* C_*(G)$ the $C_*(\PG)$-module pulled back by $\varphi$. More precisely, if $x \in \Crit(f)$, $\varphi^*C_*(G)_x = C_*(G_{\varphi(x)})$ endowed with the $C_*(\PG)$-module structure, 

$$\forall \alpha \in C_*(G_{\varphi(x)}), \ \forall \sigma \in C_*(\mathcal{P}_{x \to y}X), \  \sigma\cdot\alpha \coloneqq \Phi^{op}_*(\varphi_*\circ \sigma \otimes \alpha) \in C_*(G_{\varphi(y)}) = \varphi^*C_*(G)_y.$$
\end{defi}

\begin{thm}\cite[Theorem 8.1.1]{BDHO23}\label{thm : Functoriality}
There exist two maps called the direct map and the shriek map, 

$$\varphi_*\colon H_*(X, \varphi^*C_*(G)) \to H_*(Y, C_*(G)) \ \textup{and} \ \varphi_!\colon H_*(Y ,C_*(G) ) \to H_{n-m+*}(X, \varphi^* C_*(G))$$ 

that satisfy the following properties:

\begin{enumerate}
    \item \textbf{Identity:} $\textup{Id}_* = \textup{Id}_!= \textup{Id}$ in homology.
    \item \textbf{Composition:} Let $\psi\colon Z^k \to X $ be a continuous map. Then $$(\varphi \circ \psi)_* = \varphi_* \circ \psi_*\colon H_*(Z, \varphi^*\psi^*C_*(G)) \to H_*(X, \varphi^*C_*(G)) \to H_*(Y, C_*(G))  $$ and $$(\varphi \circ \psi)_! = \psi_! \varphi_!\colon H_*(Y, C_*(G)) \to H_{n-m+*}(X, \varphi^*C_*(G)) \to H_{k-m+*}(Z, \varphi^*\psi^*C_*(G)).$$
    \item \textbf{Homotopy:} Two homotopic maps induce the same direct and shriek maps in homology.
    \item \textbf{Spectral sequence:} The morphisms $\varphi_*$ and $\varphi_!$ in homology are limits of morphisms between spectral sequences associated with the corresponding enriched complexes, given at the second page by 
    $$\varphi_{p,*}\colon H_p(X, \varphi^* H_q(G)) \to H_p(Y, H_q(G)) \ \textup{and} \ \varphi_{p,!}\colon H_p(Y,H_q(G)) \to H_{p+n-m}(X, \varphi^* H_q(G)). $$
\end{enumerate}
\end{thm}

We make the construction of shriek maps explicit when $\varphi$ is smooth. Assume that \begin{equation*}
   \varphi \lvert_{W^s(x)} \pitchfork W^u(y)
\end{equation*}
for $x \in \Crit(f)$ and $y \in \Crit(g)$.

Consider the manifold of dimension $|y|-|x|+n-m$ $$\overline{\mathcal{M}}^{\varphi_!}(y,x) = \wb^s(x) \cap \varphi^{-1}(\wb^u(y)),$$

oriented by the rule $$\left( \Or \ \overline{W}^s_g(y), \ \Or \ \overline{\mathcal{M}}^{\varphi_!}(y,x) \right) = \Or \ \overline{W}^s_f(x).$$

\begin{rem}
    There is an identification $$\overline{\mathcal{M}}^{\varphi_!}(y,x) \cong \wbi{s}{f}{x} \ftimes{\varphi}{} \wbi{u}{g}{y}  = \left\{(a,b) \in \wbi{s}{f}{x} \times \wbi{u}{g}{y}, \ b = \varphi(a)\right\}$$ given by $$a \in \wb^s(x) \cap \varphi^{-1}(\wb^u(y)) \Leftrightarrow (a, \varphi(a)) \in \wbi{s}{f}{x} \ftimes{\varphi}{} \wbi{u}{g}{y}.$$ We see an element $\lambda$ of this moduli space as a half-infinite trajectory $\lambda_X$ going from a point $a \in \wb^s(x)$ to $x$ followed by a half-infinite trajectory $\lambda_Y$ going from $y$ to $\varphi(a) \in \wb^u(y)$.
\end{rem}

We consider  a representing chain system  $\left\{\sigma^{\varphi_!}_{y,x}\in C_{|y|-|x|+n-m}(\overline{\mathcal{M}}^{\varphi_!}(y,x))\right\}$ of the moduli spaces $\overline{\mathcal{M}}^{\varphi_!}(y,x)$.

This family satisfies the equation \begin{equation*}
    \partial \sigma^{\varphi_!}_{y,x} = \sum_{w \in \Crit(g)} s^Y_{y,w}\times \sigma^{\varphi_!}_{w,x} - \sum_{z \in \Crit(f)} (-1)^{|y|-|z|+n-m} \sigma^{\varphi_!}_{y,z}\times s^X_{z,x}.
\end{equation*}

The family of evaluation maps $$\Gamma^{\varphi_!}_{y,x}\colon \overline{\mathcal{M}}^{\varphi_!}(y,x) \to \mathcal{P}_{y \to \varphi(x)} Y$$ is defined as follows

An element $\lambda = (\lambda_X,\lambda_Y) \in \overline{\mathcal{M}}^{\varphi_!}(y,x)$ is a trajectory consisting of two half-infinite trajectories $\lambda_X$ going from a certain point $a \in \wbi{s}{f}{x}$ to $x$ and $\lambda_Y$ going from $y$ to $\varphi(a) \in \wbi{u}{g}{y}$. We define 

$$\Gamma^{\varphi_!}_{y,x}(\lambda_X, \lambda_Y) = \Gamma^Y_y(\lambda_Y) \cdot \varphi (\Gamma^X_{x}(\lambda_X)).$$

We consider the family $\left\{\nu^{\varphi_!}_{y,x} = \Gamma^{\varphi_!}_{y,x,*}( \sigma^{\varphi_!}_{y,x}) \in C_{|y|-|x|+n-m}(\mathcal{P}_{y \to \varphi(x)} Y)\right\}.$

This family satisfies \begin{equation*}
    \partial \nu^{\varphi_!}_{y,x} = \sum_{w \in \Crit(g)} m^Y_{y,w}\cdot \nu^{\varphi_!}_{w,x} - \sum_{z \in \Crit(f)} (-1)^{|y|-|z|+n-m} \nu^{\varphi_!}_{y,z}\cdot \varphi_*\circ m^X_{z,x}.
\end{equation*}

The shriek map $\varphi_!\colon C_*(Y ,C_*(G)) \to C_*(X,\varphi^*C_*(G))$ is then defined by

$$\varphi_!(y \otimes \alpha) = \sum_{x} x \otimes \nu^{\varphi_!}_{y,x} \cdot \alpha.$$

We can extend this definition to continuous maps using the identification morphism described in \cite[Proposition 5.1]{PPDG} and \cite[Proposition 8.2.1]{BDHO23} as well as the fact that any continuous map can be approximated by smooth maps.

\subsection{Equivalence of the two conventions}\label{subsec: Equivalence of the two conventions}

We now build a quasi-isomorphism between our complex (denoted here $C^{CR}_*(X ,C_*(F))$) and the complex $C^{BDHO}_*(X ,C_*(F))$ the conventions found in \cite{BDHO23}. We fix a fibration $F \hookrightarrow \fibration$ and a transitive lifting function $\Phi \colon  E \ftimes{\pi}{\ev_0} \mathcal{P}X \to E$  associated with this fibration. The set $\Xi = (f,\xi,o,s_{x,y})$ denotes a set of DG Morse data.

\subsubsection{The original convention}
We recall the definition of Morse homology with DG coefficients as in \cite{BDHO23}.

In order to define $C^{BDHO}_*(X ,C_*(F))$, we also need to fix the following data:
\begin{itemize}
    \item we fix $\star\in X$ a preferred basepoint and denote $\Omega X = \left\{ \gamma \in \mathcal{P}X\mid  \ev_0(\gamma)=\ev_1(\gamma)=\star \right\}$ the space of Moore based loops;
    \item we denote $\mathcal{Y}$ for a tree in $X$ connecting all critical points to the basepoint $\star$;
    \item we define $p \colon  X \to \Xq$ as the canonical projection that collapses the tree $\mathcal{Y}$;
    \item we choose a homotopy inverse $\theta \colon \Xq \to X$ of $p$ such that $\theta([\mathcal{Y}]) = \star$.
\end{itemize}

By definition, for $x,y \in \Crit(f)$, we have  $$m^{BDHO}_{x,y}=(\theta\circ p)_*(m_{x,y}) \in C_*(\Omega X).$$

To differentiate the two conventions, we denote $m_{x,y}^{CR}\coloneqq m_{x,y}$.

By \cite{BDHO23}, we have $$ C_*^{BDHO}(X ,C_*(F)) = C_*(F_{\star}) \otimes \Z \cdot\Crit(f)$$ with the differential given on the generators by 
$$\partial(\alpha \otimes x) = \partial \alpha \otimes x + (-1)^{|\alpha|} \sum_{y \in \Crit(f)} \alpha \cdot m_{x,y}^{BDHO} \otimes y.$$

Our convention differs from the original convention in two ways: we tensor the degree shifting $\Z\cdot x$ on the left instead of the right and we do not identify all fibres with $F_\star$. The first difference is rather superficial, but the second difference is more involved.

\subsubsection{The quasi-isomorphism}

We now build the wanted quasi-isomorphism.

\begin{prop}\label{prop: qiso BDHO CR}
    There exists a quasi-isomorphism $$C_*^{CR}(X ,C_*(F)) \to C_*^{BDHO}(X ,C_*(F)).$$
\end{prop}
\begin{proof}

Let $E' \overset{\pi'}{\to} X$ be the pull-back fibration of $\fibration$ by $\theta \circ p \colon X \to X.$ Recall that $E' = \{(a,e) \in X \times E, \ \theta \circ p(a) = \pi(e)\}.$

Let $\Phi'\colon E' \ftimes{\pi}{\ev_0} \mathcal{P} X \to E'$ defined by $$\Phi'((a,e),\gamma) = (\ev_1(\gamma), \Phi(e,\theta \circ p \circ \gamma))$$ be a transitive lifting function associated with this fibration.

This fibration is a particularly well-suited intermediate object, since if $x\in \Crit(f)$, there is an identification $$F'_x \coloneqq \pi'^{-1}(x) =\{(x,e), \ \pi(e)=\theta\circ p(x) = \star\} \cong F_\star.$$
    
    The map $\theta \circ p$ is homotopic to the identity of $X$ and consider $H\colon [0,1] \times X \to X$ a homotopy such that $H(0,\cdot) = \Id_X$ and $H(1,\cdot) = \theta \circ p$. For $a \in X$, denote $H(\cdot,a) \in \mathcal{P}_{a \to \theta \circ p(a)} X$.

    In particular, $\deffct{\kappa}{E}{E'}{e}{(\pi(e), \Phi(e, H(\cdot,\pi(e))))}$ is a homotopy equivalence and  a morphism of fibrations. Therefore Theorem \ref{thm:morphism of fibrations} gives a quasi-isomorphism $$\widetilde{\kappa} \colon C^{CR}_*(X  ,C_*(F)) \to C^{CR}_*(X  ,C_*(F')).$$

    This is a simple computation to prove that the complex $$C^{mid}_*(X , C_*(F')) = \bigoplus_{x \in \Crit(f)} C_{*}(F'_x) \otimes \Z \cdot x$$ with differential given on the generators by $$\partial(\alpha \otimes x) = \partial' \alpha \otimes x + (-1)^{|\alpha|} \sum_y \Phi'_*(\alpha \otimes m^{CR}_{x,y}) \otimes y$$ is isomorphic to $C^{CR}_*(X , C_*(F'))$ via the isomorphism $$\deffct{\textup{sw}}{C_*^{CR}(X , C_*(F'))}{C^{mid}_*(X , C_*(F'))}{x \otimes \alpha}{(-1)^{|\alpha||x|}\alpha \otimes x}.$$

    Now, we remark that $C^{mid}_*(X ,C_*(F')) \cong C^{BDHO}_*(X ,C_*(F))$ under the identification $F'_x \cong F_\star$ since $$\Phi'_*(\alpha \otimes m^{CR}_{x,y}) = \Phi_*(\alpha \otimes (\theta\circ p)_*(m^{CR}_{x,y})) = \Phi_*(\alpha \otimes m^{BDHO}_{x,y}).$$

    Therefore, $$\textup{sw} \circ \widetilde{\kappa}\colon C^{CR}_*(X  ,C_*(F)) \to C^{BDHO}_*(X ,C_*(F))$$ is a quasi-isomorphism.
\end{proof}

\subsection{Chas-Sullivan product}\label{appendix subsec: CS product}

The product $\CS\colon H^{BDHO}_{n+*}(X ,C_*(F))^{\otimes 2} \to H^{BDHO}_{n+*}(X ,C_*(F))$ is defined in \cite{Rie24} as the composition of three maps:

    \begin{itemize}
        \item[$\bullet$] The Künneth map $K\colon H^{BDHO}_{n+*}(X ,C_*(F))^{\otimes 2} \to H^{BDHO}_{2n+*}(X \times X ,C_*(F \times F))$ defined at the chain-level by $$\deffct{K}{C_*(X  ,C_*(F))^{\otimes 2}}{C_{*}(X \times X, m^K ,C_*(F \times F))}{(\alpha \otimes x) \otimes (\beta \otimes y)}{(-1)^{|\beta||x|} (\alpha \times \beta) \otimes (x,y)},$$ where $m^K$ is the Künneth twisting cocycle. See \cite[\S 6]{Rie24} for more details on this map and this twisting cocycle.

        \item[$\bullet$] The shriek map $\Delta_!\colon H^{BDHO}_{2n+*}(X \times X ,C_*(F \times F)) \to H^{BDHO}_{n+*}(X , \Delta^*C_*(F \times F))$ defined at the chain-level  by 
        $$\deffct{\Delta_!}{C^{BDHO}_*(X \times X, m^K ,C_*(F \times F))}{C^{BDHO}_*(X , \Delta^*C_*(F \times F))}{\alpha \otimes(x,y)}{ \sum_z\alpha \cdot \nu_{(x,y),z} \otimes z},$$
        where $\nu_{(x,y),z} \in C_{|x|+|y|-|z|}(\Omega X)$ is a family of chains obtained by the evaluations of a family of chain representatives of the fundamental classes of the moduli spaces $$\overline{\mathcal{M}}^{\Delta_!}((x,y),z) = \Delta^{-1}(\wb^u(x,y)) \cap \wb^s(z)$$ through the evaluation maps $q^{\Delta}_{(x,y),z}\colon \mathcal{M}^{\Delta_!}((x,y),z) \to \Omega X \times \Omega X$. See \cite[Section 10.4]{BDHO23} for more details on this map.

        \item[$\bullet$] The map $\widetilde{\mu}\colon H^{BDHO}_{n+*}(X, \Delta^*C_*(F\times F)) \to H^{BDHO}_{n+*}(X ,C_*(F))$ induced by the morphism of fibrations $\mu\colon  E \ftimes{\pi}{\pi} E \to E$. See \cite[Theorem 5.8]{Rie24} for more details on this map.
    \end{itemize}

Since each of these maps is defined in our conventions (see \S \ref{subsec:PG-modules and morphisms of fibrations}, \S \ref{subsec:Künneth map} and \S \ref{subsec:Shriek map}), we define the Chas-Sullivan product $$ \CS= \widetilde{\mu}\circ \Delta_! \circ K\colon H^{CR}_{n+*}(X,C_*(F))^{\otimes 2} \to H^{CR}_{n+*}(X,C_*(F)).$$

\begin{prop}
    This definition of $\CS$ coincides with the Chas-Sullivan product defined in \cite{Rie24} via the quasi-isomorphism of Proposition \ref{prop: qiso BDHO CR}.
\end{prop}

\begin{proof}

  We denote $\mu'\colon E' \ftimes{\pi'}{\pi'} E' \to E'$ the morphism of fibrations $\mu'((a,e),(a,e')) = (a,\mu(e,e'))$.

Consider the morphism of fibrations $\Delta^*\kappa \times \kappa\colon E' \ftimes{\pi'}{\pi'} E' \to E' \ftimes{\pi'}{\pi'} E'$, $(e,e') \mapsto (\kappa(e), \kappa(e'))$.




The following diagram commutes.
    \[\xymatrix@C=1em{
    H^{CR}_{n+*}(X ,C_*(F))^{\otimes 2} \ar[r]^-{K} \ar[d]^{\widetilde{\kappa}\otimes \widetilde{\kappa}} \ar@{}[rd]|{\blue{(1)}} & H^{CR}_{2n+*}(X \times X ,C_*(F\times F)) \ar@{}[rd]|{\blue{(2)}} \ar[r]^-{\Delta_!} \ar[d]^{\widetilde{\kappa \times \kappa}} & H^{CR}_{n+*}(X , \Delta^*C_*(F \times F))\ar@{}[rd]|{\blue{(3)}} \ar[r]^-{\widetilde{\mu}} \ar[d]^{\Delta^*\widetilde{\kappa \times \kappa}} & H^{CR}_{n+*}(X ,C_*(F)) \ar[d]^{\widetilde{\kappa}} \\
    H^{CR}_{n+*}(X ,C_*(F'))^{\otimes 2} \ar@{}[rd]|{\blue{(4)}}\ar[r]^-{K} \ar[d]^{\textup{sw} \otimes \textup{sw} } & H_{2n+*}^{CR}(X \times X ,C_*(F'\times F')) \ar[r]^-{\Delta_!} \ar[d]^{\textup{sw}} \ar@{}[rd]|{\blue{(5)}} & H^{CR}_{n+*}(X , \Delta^*C_*(F' \times F')) \ar[d]^{\textup{sw}} \ar@{}[rd]|{\blue{(6)}} \ar[r]^-{\widetilde{\mu}'} & H^{CR}_{n+*}(X ,C_*(F')) \ar[d]^{\textup{sw}} \\
    H_{n+*}^{BDHO}(X ,C_*(F))^{\otimes 2} \ar[r]^-{K} & H_{2n+*}^{BDHO}(X \times X ,C_*(F\times F)) \ar[r]^-{\Delta_!} & H_{n+*}^{BDHO}(X ,\Delta^*C_*(F \times F)) \ar[r]^-{\widetilde{\mu}} & H_{n+*}^{BDHO}(X ,C_*(F)).
    }\]

  \begin{itemize}
      \item[\blue{(1)}] This is a direct application of Theorem \ref{thm:morphism of fibrations} and Proposition \ref{prop : Kunneth map} that this square commutes.
      \item[\blue{(2)}] The same computation as \cite[Proposition 5.19]{Rie24} proves that this square commutes.
      \item[\blue{(3)}]The maps $\mu' \circ (\kappa \times \kappa)\colon E \ftimes{\pi}{\pi}E \to E'$ and $\kappa \circ \mu\colon E \ftimes{\pi}{\pi}E \to E'$ are homotopic. Indeed, let $e,v \in E$ such that $\pi(e)=\pi(v)\eqqcolon a$. Then, $$\mu' \circ (\kappa \times \kappa)(e,v) = \left(a, \mu(\Phi(e,H(\cdot,a)), \Phi(v, H(\cdot,a))) \right)$$ and $$\kappa \circ \mu(e,v) = (a,\Phi(\mu(e,v),H(\cdot,a))).$$

      This is a corollary of Lemma \ref{lem: morphism of fibrations formula with paths} that there exists $\mu_2\colon I \times (\mathcal{P}X \ftimes{\ev_0}{\pi} (E\ftimes{\pi}{\pi}E)) \to E $ such that $$(t,(e,e')) \mapsto (\pi(e),\mu_2(t, H(\cdot,a), (e,e')))$$ is a homotopy between $\mu' \circ (\kappa \times \kappa)$ and $\kappa \circ \mu$. The commutativity of this diagram is then a direct application of Theorem \ref{thm:morphism of fibrations}.
      \item[\blue{(4)}] A direct computation on chains shows that this diagram commutes already at the chain-level. 
      \item[\blue{(5)}] Again, a direct computation using the shifted grading $|x|-n$ proves that this square commutes at the chain-level.
      \item[\blue{(6)}] A more involved sign computation that we do not provide here shows that this diagram commutes at the chain-level.     
  \end{itemize}   
\end{proof}

\section{Oriented trees}\label{sec:or trees}
In this section, we construct explicit orientations for $\Tree_k(\ux;y)$, the moduli spaces of perturbed Morse gradient trees. 

The goal is in particular to prove the following lemma:

\begin{lemme}\label{lem:existence of orientations appendix}
    There exist orientations $[\Tree_k(\ux;y)]\in \Or(\Tree_k(\ux;y))$ such that, for all points $p\in \mathcal{T}_l^{\mathbb{X}}(\uw;y)\times \mathcal{T}^{\mathbb{X}}_m(\uz;w_{i})\subseteq \partial \Tree_k(\ux;y)$, it holds
    \begin{align*}
        [\Tree_k(\ux;y)]_p=\varepsilon_0(\uw;\uz;i)[\vec{n}]\times [\Tree_l(\uw;y)]_p\times [\Tree_m(\uz;w_{i})]_p\in \Or_p(\Tree_{k}(\ux;y))
    \end{align*}
    where $\vec{n}\in T_p\Tree_{k}(\ux;y)$ is a vector pointing out of the manifold and $\varepsilon_0(\uw;\uz;i)$ is the sign appearing in the equality
    \begin{align*}
        \partial(\tau_{\ux,y}\otimes m_{\ux,y})=\sum_{\substack{\uw\circ_{i}\uz=\ux}} \varepsilon_0(\uw;\uz;i) (\tau_{\uw,y}\otimes m_{\uw,y})\circ _{i}(\tau_{\uz,w_{i}}\otimes m_{\uz,w_{i}}).
    \end{align*}
    Moreover for $k=1$, we have $[\Tree_1(x;y)]=(-1)^{|x|+|x||y|}[\overline{\mathcal{L}}(x,y)]$ for $[\overline{\mathcal{L}}(x,y)]$ as in \eqref{eq: or L(x,y)}.
\end{lemme}

\subsection{Picking orientations}\label{subsec:def orient}
In this subsection, we give the definitions of the orientations $[\Tree_k(\ux;y)]$ as in Lemma \ref{lem:existence of orientations appendix}. The orientation for $k=1$ is already explicitly defined in Lemma \ref{lem:existence of orientations appendix}. We thus assume $k>1$.

The key to choosing the orientation of $\Tree_k(\ux;y)$ for $k>1$ is to construct an isomorphisms of graded abelian groups
\begin{align}\label{eq:iso for k>1}
    \Or(\Tree_k(\ux;y))\cong \underline{\Sigma}^{\|\cdot \|}(\ux;y)\otimes \Or(\overline{\mathcal{T}}_k)
\end{align}
that commute with boundary inclusions in an appropriate sense.

We recall that we fixed orientations $[X]\in \Or(X)$ and $[\wb^u(x)]\in \Or(\wb^u(x))$ for all critical points $x\in \Crit(f)$. We denote by $[\wb^s(x)]\in \Or([\wb^s(x)])$ the orientation that satisfies:
\begin{align}\label{eq:def of [ws]}
    [\wb^u(x)]_x\times [\wb^s(x)]_x=[X]_x\in \Or_x(X).
\end{align}
Moreover, let $[\mathcal{T}_k]\in \Or(\overline{\mathcal{T}_k})$ be the orientations corresponding to $m_k$ as in Lemma \ref{lem:or Ai to cell Ai}. 

For shorthand, we denote by $\wb^u(\ux)\coloneqq\wb^u(x_1)\times \dots \times \wb^u(x_k)$ and $[\wb^u(\ux)]\in \Or(\wb^u(\ux))$ the orientation $[\wb^u(x_1)]\times \dots \times [\wb^u(x_k)]$. Analogously for $\wb^s$.

In order to produce \eqref{eq:iso for k>1}, we use the morphism 
\begin{align*}
    \Phi_{\ux}\colon  \wb^u(\ux)\times \mathcal{T}_k\to X^k
\end{align*}
as in Definition \ref{defi : Phi_X} and the fact that $\Phi_{\ux}$ intersects $\Delta^k(\wb(y))$ transversely and $\Phi^{-1}_{\ux}(\Delta^k(W^s(y)))\cong \mathcal{T}^{\X}_k(\ux;y)$. 

In particular, the following lemma relates $\Tree_k(\ux;y)$ to the manifolds $\wb^u(\ux)$, $\wb^s(y)$, $X$ and $\overline{\mathcal{T}_k}$. These are all manifolds we have already oriented.

\begin{lemme}\label{lem:transport orient}
    The map $\Phi_{\ux}$ induces a degree-preserving map
    \begin{align*}
        f_{\ux,y}\colon\Hom(\Or(\wb^s(y)),\Or(X^k))\to \Hom(\Or(\Tree_k(\ux;y)),\Or(\wb^u(\ux)\times \overline{\mathcal{T}_k})).
    \end{align*}
\end{lemme}
\begin{proof}
    We first assume that $\Tree_k(\ux;y)$ is connected. An orientation $[\nu]$ of the normal bundle of the inclusion $\mathcal{T}^\X_k(\ux;y)\subseteq W^u(\ux)\times \mathcal{T}_k$ defines an element in $\Hom(\Or(\Tree_k(\ux;y)),\Or(\wb^u(\ux)\times\overline{\mathcal{T}_k}))$ by
    \begin{align*}
        \Or(\Tree_k(\ux;y))\ni [\Tree_k(\ux;y)]\mapsto [\nu]\times [\Tree_k(\ux;y)]\in  \Or(\wb^u(\ux)\times \overline{\mathcal{T}_k}).
    \end{align*}
    By transversality, $\Phi_{\ux}$ gives an isomorphism of this normal bundle to the normal bundle of the inclusion $\Delta^k(W^s(y))\subseteq X^k$. 

    By an analogous argument to the one given above, an orientation of the normal bundle $\Delta^k(W^s(y))\subseteq X^k$ defines an element in $\Hom(\Or(\wb^s(y)),\Or(X^k))$. Combining this, we get the desired isomorphism
    \begin{align*}
        \Hom(\Or(\wb^s(y)),\Or(X^k))\cong \Hom(\Or(\Tree_k(\ux;y)),\Or(\wb^u(\ux)\times \overline{\mathcal{T}_k})).
    \end{align*}

    If $\Tree_k(\ux;y)$ is not connected, we take the sum over all connected components.
\end{proof}

\begin{rem}\label{rem:cancelling orient}
    We work with $\Hom(\Or(A),\Or(B))$ rather than the notation $\Or^{-1}(A)\coloneqq\Hom(\Or(A),\mathbbm{1})$ as in \cite{Abo09}. 
    
    We do this to notationally emphasise how orientation spaces cancel: in the following, we use maps of the form
    \begin{align*}
        \Hom(\Or(A),\Or(A))\to \mathbbm{1}.
    \end{align*}
    If $A$ is connected, this is an isomorphism. 

    In the following, we thus may cancel terms of the form $\Hom(C,D)$ with terms $\Hom(D,C)$ for $C$ and $D$ tensor products of orientation spaces. 

    We may also add terms of the form $\Hom(C,D)\otimes \Hom(D,C)$ as long as all the manifolds involved in the terms $C$ and $D$ are connected. The only manifolds involved in our calculation that may not be connected are $\Tree_k(\ux;y)$.

    In accordance with earlier conventions (see Remark \ref{rem:ev and coev}), we use only the evaluation map to cancel and add term, rather than the coevaluation maps. We explain this subtlety in the base case where $C=\Or(A)$ for $A$ connected and $D=\mathbbm{1}$.
    
    Evaluation is given by the following chain of isomorphisms:
    \begin{align*}
        \Hom(\Or(A),\mathbbm{1})\otimes \Or(A)&\cong \mathbbm{1},\\
        f\otimes [A]&\mapsto f([A]).
    \end{align*}

    There is another tempting map that we do \textit{not} use: coevaluation is given by the following isomorphism
    \begin{align*}
        \mathbbm{1}&\cong \Hom(\Or(A),\Or(A))\cong \Or(A)\otimes \Hom(\Or(A),\mathbbm{1}),\\
        1&\mapsto\id_{\Or(A)}.
    \end{align*}
    Both isomorphisms are natural and symmetric monoidal, that is, they commute with diffeomorphisms, cartesian and tensor products and permutation of factors. However, the two isomorphisms only agree up to the sign $(-1)^{\dim(A)}$ (see Remark \ref{rem:ev and coev}).
\end{rem}

We now construct a degree-preserving map:
\begin{align*}
    \alpha_{\ux,y}\colon \Or(\Tree_k(\ux;y))\to \underline{\Sigma}^{\|\cdot \|}(\ux;y)\otimes \Or(\overline{\mathcal{T}_k}).
\end{align*}
It is given by the following composition
\begin{align}\label{eq:alpha x y}
    \begin{split}
        \Or(\Tree_k(\ux;y))
        \cong&\ \Or(\Tree_k(\ux;y))\otimes \Hom(\Or(\wb^s(y)),\Or(\wb^u(\ux))\otimes\Or(\wb^s(\ux)))\\
        &\ \otimes\Hom(\Or(\wb^s(\ux))\otimes \Or(\wb^u(\ux)),\Or(\wb^s(y)))\\
        \cong&\   \Or(\Tree_k(\ux;y))\otimes \Hom(\Or(\wb^s(y)),\Or(X^k))\\
        &\ \otimes\Hom(\Or(\wb^s(\ux))\otimes \Or(\wb^u(\ux)),\Or(\wb^s(y)))\\
        \overset{f_{\ux,y}}{\to}\!\! &\ \  \Or(\Tree_k(\ux;y))\otimes\Hom(\Or(\Tree_k(\ux;y)),\Or(\wb^u(\ux)\times \overline{\mathcal{T}_k}))\\
        &\ \ \otimes\Hom(\Or(\wb^s(\ux))\otimes \Or(\wb^u(\ux)),\Or(\wb^s(y)))\\
        \to&\  \Hom(\Or(\wb^s(\ux)),\Or(\wb^s(y)))\otimes \Or(\overline{\mathcal{T}_k})\\
        \cong&\  \underline{\Sigma}^{\|\cdot \|}(\ux;y)\otimes \Or(\overline{\mathcal{T}_k}),
    \end{split}
\end{align}
where the last isomorphism $\Hom(\Or(\wb^s(\ux)),\Or(\wb^s(y)))\cong \underline{\Sigma}^{\|\cdot \|}(\ux;y)$ comes from the orientations $[\wb^s(x)]$.

\begin{defi}
    For $k>1$, $[\Tree_k(\ux;y)]\in \Or(\Tree_k(\ux;y))$ is defined as the orientation such that for all points $p\in \Tree_k(\ux;y)$, it holds 
    \begin{align*}
        \alpha_{\ux,y}( [\Tree_k(\ux;y)]_p)= \tau_{\ux,y}\otimes [\overline{\mathcal{T}_k}].
    \end{align*}
\end{defi}

\subsection{Sign computations}\label{subsec:sign comp}
In this subsection, we show that the orientations defined in Subsection \ref{subsec:def orient} indeed satisfy the relation as in Lemma \ref{lem:existence of orientations appendix}.

For this, it is useful to have a more explicit description of the orientation $[\Tree_k(\ux;y)]$ for the case $k>1$.

\begin{lemme}\label{lem:explicit orient k>1}
    For $\nu_{\ux,y}$ the normal bundle of $\Delta_*^k([\wb^s(y)])\subseteq [X]^k$ and $[\nu_{\ux,y}]$ the orientation that satisfies
    \begin{align*}
        \nu_{\ux,y}\times \Delta_*^k([\wb^s(y)])=[X]^k,
    \end{align*}
    it holds for all $p\in \Tree_k(\ux;y)$
    \begin{align*}
        f_{\ux,y}([\nu_{\ux,y}])_p\times [\Tree_k(\ux;y)]_p=(-1)^{\eta(\ux;y)}([\wb^u(\ux)]\times [\overline{\mathcal{T}_k}])_p\in \Or_p(\wb^u(\ux)\times \overline{\mathcal{T}_k})
    \end{align*}
    with
    \begin{align*}
        \eta(\ux;y)=\sum_j (k-j)n \|x_j\|+\sum_{j>j'} \|x_j\|\|x_{j'}\|+\|y\|.
    \end{align*}
\end{lemme}
\begin{proof}
    For an orientation $[A]\in \Or(A)$, we denote by $[A]^\lor\in \Hom(\Or(A),\mathbbm{1})$ the map that sends $[A]$ to $1$. This notation reverses order of the tensor products, that is, under the natural isomorphism, we have:
    \begin{align*}
        \Hom(\Or(A)\otimes \Or(B),\mathbbm{1})&\to \Hom(\Or(B),\mathbbm{1})\otimes \Hom(\Or(A),\mathbbm{1}),\\
        ([A]\otimes [B])^\lor&\mapsto [B]^\lor\otimes [A]^\lor.
    \end{align*}
    This holds because both sides evaluate $[A]\otimes [B]$ to $1$.

    We fix a point $p\in \Tree_k(\ux;y)$. Under the first map from the composition \eqref{eq:alpha x y} defining $\alpha_{\ux,y}$, we thus have
    \begin{align*}
        [\Tree_k(\ux;y)]_p&\mapsto [\Tree_k(\ux;y)]_p\otimes [\wb^s(\ux)]^\lor\otimes [\wb^u(\ux)]^\lor\otimes [\wb^s(y)]^\lor\otimes [\wb^s(y)]\otimes [\wb^u(\ux)]\otimes [\wb^s(\ux)]\\
        &\mapsto (-1)^{\eta_1}[\Tree_k(\ux;y)] \otimes [\wb^u(\ux)]\otimes [\wb^s(\ux)]\otimes [\wb^s(y)]^\lor \otimes [\wb^s(y)]\otimes[\wb^u(\ux)]^\lor\otimes [\wb^s(\ux)]^\lor
    \end{align*}
    where 
    \begin{align*}
        \eta_1&\equiv (|\wb^u(\ux)|+|\wb^s(\ux)|)(|\wb^u(\ux)|+|\wb^s(\ux)|)+|\wb^s(\ux)||\wb^u(\ux)|\\
        &\equiv nk+\|\ux\||\ux|\equiv nk+\|\ux\|nk+\|\ux\|
    \end{align*}
    in $\Z/2\Z$.

    The second map in the composition \eqref{eq:alpha x y} identifies $\Or(\wb^u(\ux))\otimes \Or(\wb^s(\ux))$ with $\Or(X^k)$ which sends
    \begin{align*}
        [\wb^u(\ux)]\otimes [\wb^s(\ux)]&\mapsto (-1)^{\eta_2}[\wb^u(x_1)]\times [\wb^s(x_1)]\otimes \dots \otimes [\wb^u(x_k)]\times [\wb^s(x_k)]=(-1)^{\eta_2}[X]^k
    \end{align*}
    where 
    \begin{align*}
        \eta_2\equiv \sum_{j<j'} |\wb^s(x_{j})||\wb^u(x_{j'})|\equiv \sum_{j<j'}\|x_j\| |x_{j'}|\equiv \sum_{j<j'} \|x_j\|\|x_{j'}\|+\sum_j (k-j)n\|x_j\|
    \end{align*}
    by $[\wb^u(x_i)]\times [\wb^s(x_i)]=[X]$ as in \eqref{eq:def of [ws]}.

    In our notation, the element $[X]^k\otimes [\wb^s(y)]^\lor\in \Hom(\Or(\wb^s(y),\Or(X^k)))$ corresponds to $[\nu_{\ux,y}]$. Under the composition of the first three maps in \eqref{eq:alpha x y}, we thus have:
    \begin{align*}
        [\Tree_k(\ux;y)]_p\mapsto (-1)^{\eta_1+\eta_2}[\Tree_k(\ux;y)]_p \otimes f_{\ux,y}([\nu_{\ux,y}]) \otimes [\wb^s(y)]\otimes[\wb^u(\ux)]^\lor\otimes [\wb^s(\ux)]^\lor.
    \end{align*}
    The next isomorphism is given by shuffling and then evaluating $f_{\ux,y}([\nu_{\ux,y}])$ at $[\Tree_k(\ux;y)]_p$:
    \begin{align*}
        [\Tree_k(\ux;y)]&\mapsto (-1)^{\eta_1+\eta_2+\eta_3} f_{\ux,y}([\nu_{\ux,y}])\otimes [\Tree_k(\ux;y)]_p\otimes [\wb^s(y)]\otimes[\wb^u(\ux)]^\lor\otimes [\wb^s(\ux)]^\lor\\
        &\mapsto (-1)^{\eta_1+\eta_2+\eta_3+\eta} [\wb^u(\ux)]\times [\overline{\mathcal{T}_k}]\otimes [\wb^s(y)]\otimes[\wb^u(\ux)]^\lor\otimes [\wb^s(\ux)]^\lor
    \end{align*}
    for $\eta\coloneqq\eta(\ux;y)$ as in the statement and 
    \begin{align*}
        \eta_3\equiv |\Tree_k(\ux;y)||\nu_{\ux,y}|\equiv (k+\|\ux\|+\|y\|)(nk+\|y\|)\equiv nk+\|y\|(nk+k+1)+ \|\ux\|nk+\|\ux\|\|y\|.
    \end{align*}
    Shuffling and evaluating at $[\wb^u(\ux)]$, we find:
    \begin{align*}
        [\Tree_k(\ux;y)]&\mapsto (-1)^{\eta_1+\eta_2+\eta_3+\eta_4+\eta} [\wb^u(\ux)]^\lor\otimes[\wb^u(\ux)]\otimes [\overline{\mathcal{T}_k}]\otimes [\wb^s(y)]\otimes [\wb^s(\ux)]^\lor\\
        &\mapsto (-1)^{\eta_1+\eta_2+\eta_3+\eta_4+\eta}  [\overline{\mathcal{T}_k}]\otimes [\wb^s(y)]\otimes [\wb^s(\ux)]^\lor
    \end{align*}
    for
    \begin{align*}
        \eta_4\equiv |\wb^u(\ux)|(|\wb^u(\ux)|+|\overline{\mathcal{T}_k}|+|\wb^s(y)|)\equiv |\ux|+k|\ux|+|\ux|\|y\|\equiv \|\ux\|+k\|\ux\|+\|\ux\|\|y\|+nk\|y\|.
    \end{align*}
    The last isomorphism 
    \begin{align*}
        \Hom(\Or(\wb^s(\ux)),\Or(\wb^s(y)))\cong \underline{\Sigma}^{\|\cdot \|}(\ux;y)
    \end{align*}
    sends $ [\wb^s(y)]\otimes [\wb^s(\ux)]^\lor$ to $\tau_{\ux,y}$. We thus find after this identification and shuffling that
    \begin{align*}
        [\Tree_k(\ux;y)]_p&\mapsto(-1)^{\eta_1+\eta_2+\eta_3+\eta_4+\eta}  [\overline{\mathcal{T}_k}]\otimes\tau_{\ux,y}\\
        &\mapsto(-1)^{\eta_1+\eta_2+\eta_3+\eta_4+\eta_5+\eta}  \tau_{\ux,y}\otimes[\overline{\mathcal{T}_k}]
    \end{align*}
    for 
    \begin{align*}
        \eta_5\equiv |\tau_{\ux,y}||\mathcal{T}_k|=(\|\ux\|+\|y\|)k
    \end{align*}
    The upshot thus is that
    \begin{align*}
        \alpha_{\ux,y}([\Tree_k(\ux;y)]_p)=(-1)^{\eta_1+\eta_2+\eta_3+\eta_4+\eta_5+\eta}\tau_{\ux,y}\otimes [\overline{\mathcal{T}_k}].
    \end{align*}
    By definition of $[\Tree_k(\ux;y)]$, we find
    \begin{align*}
        \eta&\equiv \eta_1+\eta_2+\eta_3+\eta_4+\eta_5\\
        &\equiv \sum_j (k-j)n \|x_j\|+\sum_{j>j'} \|x_j\|\|x_{j'}\|+\|y\|.
    \end{align*}
\end{proof}

We fix $\ux\in \Crit(f)^k$, $y\in \Crit(f)$, $\uw\in \Crit(f)^l$, $\uz\in \Crit(f)^m$ and $1\leq i\leq l$ such that $\Tree_l(\uw;y)\times \Tree_m(\uz;w_i)\subseteq \partial\Tree_k(\ux;y)$. 

If any of the three manifolds is empty, there is nothing to prove. Otherwise, we fix a point $p\in  \partial\Tree_k(\ux;y)$ in the interior of $\Tree_l(\uw;y)\times \Tree_m(\uz;w_i)$ and we fix $\vec{n}\in T_p \Tree_k(\ux;y)$ pointing out of the manifold. 

We want to show that 
\begin{align*}
        [\Tree_k(\ux;y)]_p=\varepsilon_0[\vec{n}]\times [\Tree_l(\uw;y)]_p\times [\Tree_m(\uz;w_{i})]_p\in \Or_p(\Tree_{k}(\ux;y))
\end{align*}
for $\varepsilon_0\coloneqq\varepsilon_0(\uw;\uz;i)$ the sign appearing in the equality
\begin{align*}
    \partial(\tau_{\ux,y}\otimes m_{\ux,y})=\sum_{\substack{\uw\circ_{i}\uz=\ux}} \varepsilon_0(\uw;\uz;i) (\tau_{\uw,y}\otimes m_{\uw,y})\circ _{i}(\tau_{\uz,w_{i}}\otimes m_{\uz,w_{i}}).
\end{align*}

We first compute $\varepsilon_0$:
\begin{lemme}\label{lem:eps0}
    It holds $\varepsilon_0=(-1)^{\eta_0}$ where
    \begin{align*}
        \eta_0\coloneqq\|\ux\|+\|y\|+i+m(l-i)+(\|\uz\|+\|w_i\|)\left(\sum_{j=1}^{i-1} \|w_j\|+l\right).
    \end{align*}
\end{lemme}
\begin{proof}
    We compute:
    \begin{align*}
        \partial(\tau_{\ux,y}\otimes m_{\ux,y})=(-1)^{\|\ux\|+\|y\|} \tau_{\ux,y}\otimes \partial(m_{\ux,y}).
    \end{align*}
    The term corresponding to $\uw\circ_i\uz=\ux$ that we are interested in is given by
    \begin{align*}
        (-1)^{\|\ux\|+\|y\|+i+m(l-i)} \tau_{\ux,y}\otimes (m_{\uw,y}\circ_im_{\uz,w_i}).
    \end{align*}
    By the computation 
    \begin{align*}
        \tau_{\uw,y}\circ_i\tau_{\uz,w_i}(\sigma^{\|x_1\|}\otimes \dots \otimes \sigma^{\|x_k\|})=&\ (-1)^{(\|\uz\|+\|w_i\|)(\sum_{j<i} \|w_j\|)} \\
        &\ \tau_{\uw,y}(\sigma^{\|x_1\|}\otimes \dots\otimes \tau_{\uz,w_i}(\sigma^{\|z_1\|}\otimes \dots\otimes \sigma^{\|z_m\|})\otimes \dots\otimes \sigma^{\|x_k\|})\\
        =&\ (-1)^{(\|\uz\|+\|w_i\|)(\sum_{j<i} \|w_j\|)} \tau_{\uw,y}(\sigma^{\|x_1\|}\otimes \dots\otimes \sigma^{\|w_i\|}\otimes \dots\otimes \sigma^{\|x_k\|})\\
        =&\ (-1)^{(\|\uz\|+\|w_i\|)(\sum_{j<i} \|w_j\|)} \sigma^{\|y\|},
    \end{align*}
    we deduce that $\tau_{\uw,y}\circ_i\tau_{\uz,w_i}=(-1)^{(\|\uz\|+\|w_i\|)(\sum_{j<i} \|w_j\|)}\tau_{\ux,y}$.

    We thus have
    \begin{align*}
        (-1)^{\|\ux\|+\|y\|+i+m(l-i)} \tau_{\ux,y}\otimes (m_{\uw,y}\circ_im_{\uz,w_i})=&\ (-1)^{\|\ux\|+\|y\|+i+m(l-i)+(\|\uz\|+\|w_i\|)(\sum_{j<i} \|w_j\|)}\\
        &\ (\tau_{\uw,y}\circ_i\tau_{\uz,w_i})\otimes (m_{\uw,y}\circ_im_{\uz,w_i})\\
        =&\ (-1)^{\|\ux\|+\|y\|+i+m(l-i)+(\|\uz\|+\|w_i\|)(\sum_{j<i} \|w_j\|+l)}\\
        &\ (\tau_{\uw,y}\otimes m_{\uw,y})\circ_i (\tau_{\uz,w_i}\otimes m_{\uz,w_i}).
    \end{align*}
\end{proof}

\subsubsection{The case \texorpdfstring{$l=m=1$}{l=m=1}}
We first assume that $l=m=1$ and thus $k=1$. In particular, we have $\Tree_1(x;y)=\overline{\mathcal{L}}(x,y)$, $\Tree_1(x;z)=\overline{\mathcal{L}}(x,z)$ and $\Tree_1(z;y)=\overline{\mathcal{L}}(z,y)$. For this case, we recall that 
\begin{align*}
    [\overline{\mathcal{L}}(x,y)]&=(-1)^{|z|+|x|}[\vec{n}]\times [\overline{\mathcal{L}}(x,z)]\times [\overline{\mathcal{L}}(z,y)]\\
    &=(-1)^{|y|+|x||y|+|x||z|+|z||y|+1}[\vec{n}]\times [\overline{\mathcal{L}}(z,y)]\times [\overline{\mathcal{L}}(x,z)] 
\end{align*}
see \eqref{eq: or L(x,y)}. By the second part of Lemma \ref{lem:existence of orientations appendix}, we have that
\begin{align*}
    [\Tree_1(x;y)]=(-1)^{|x|+|x||y|}[\overline{\mathcal{L}}(x,y)].
\end{align*}
We thus find 
\begin{align*}
    (-1)^{|z|+|z||y|+|x|+|x||z|}[\vec{n}]\times [\Tree_1(z;y)]\times [\Tree_1(x;z)]=(-1)^{|y|+|x||y|+|x||z|+|z||y|+|x|+|x||y|+1} [\Tree_1(x;y)].
\end{align*}

Simplifying the above equation, we find
\begin{align*}
    [\vec{n}]\times [\Tree_1(z;y)]\times [\Tree_1(x;z)]=(-1)^{|y|+|z|+1} [\Tree_1(x;y)]=\varepsilon_0[\Tree_1(x;y)].
\end{align*}

\subsubsection{The case \texorpdfstring{$k=l>1=m$}{k=l>1=m}} 
For this case, we use the notation that $\uw_{<i}\coloneqq(w_1,\dots, w_{i-1})$, $\uw_{\leq i}\coloneqq(w_1,\dots, w_{i})$ and so on. We can extend $\Phi_{\ux}$ to $W^u(\uw_{<i})\times \mathcal{L}(x_i,w_i)\times W^u(\uw_{\geq i})\times \mathcal{T}_k\subseteq \wb^u(\ux)\times  \mathcal{T}_k$ via 
\begin{align*}
    W^u(\uw_{<i})\times \mathcal{L}(x_i,w_i)\times W^u(\uw_{\geq i})\times \mathcal{T}_k\to W^u(\uw)\times \mathcal{T}_k\overset{\Phi_{\uw}}{\to} X^k.
\end{align*}

This shows that the maps that $\Phi_{\ux}$ and $\Phi_{\uw}$ induce as in Lemma \ref{lem:transport orient} are compatible, that is, the following diagram commutes:
\begin{center}
    \begin{tikzcd}
        \Hom(\Or(\wb^s(y)),\Or(X^k))\ar[dd, equal]\ar[r, "f_{\ux,y}"] &\Hom(\Or(\Tree_k(\ux;y)),\Or(\wb^u(\ux)\times \overline{\mathcal{T}_k}))\ar[d,"\cong"]\\
        &\makecell{\Hom(\Or(\vec{n})\otimes \Or(\mathcal{L}(x_i,w_i))\otimes \Tree_k(\uw;y),\\ \Or(\vec{n})\otimes \Or(\mathcal{L}(x_i,w_i))\otimes\Or(\wb^u(\uw)\times \overline{\mathcal{T}_k}))}\ar[d, "\cong ", "\mathrm{ev}"']\\
        \Hom(\Or(\wb^s(y)),\Or(X^k))\ar[r, "f_{\uw,y}"] &\Hom(\Or(\Tree_k(\uw;y)),\Or(\wb^u(\uw)\times \overline{\mathcal{T}_k}))
    \end{tikzcd}
\end{center}

We recall that the second vertical map on the right is given by evaluation at $[\vec{n}]\times [\mathcal{L}(x_i,w_i)]$. By Remark \ref{rem:cancelling orient}, this map is $(-1)^{|\vec{n}\times \mathcal{L}(x_i,w_i)|}\mathrm{coev}^{-1}=(-1)^{\|x_i\|+\|w_i\|}\mathrm{coev}^{-1}$ and $\mathrm{coev}(1)=\id_{\Or(\vec{n}\times \mathcal{L}(x_i,w_i))}$. Thus under the above identification, we have that a map $\alpha\in \Hom(\Or(\Tree_k(\uw;y)),\Or(\wb^u(\uw)\times \overline{\mathcal{T}_k}))$ is identified with the map $(-1)^{\|x_i\|+\|w_i\|}\mathrm{coev}(1)\otimes \alpha=(-1)^{\|x_i\|+\|w_i\|}\id\otimes \alpha$ in $\Hom(\Or(\Tree_k(\ux;y)),\Or(\wb^u(\ux)\times \overline{\mathcal{T}_k}))$.

Let $[\nu_{\ux,y}],[\nu_{\uw,y}]\in \Hom(\Or(\wb^s(y)),\Or(X^k))$ be as in Lemma \ref{lem:explicit orient k>1}. We note that $[\nu_{\ux,y}]=[\nu_{\uw,y}]$. By the above diagram, we thus have that $f_{\ux,y}([\nu_{\ux,y}])$ is sent to $(-1)^{\|x_i\|+\|w_i\|}\id\otimes f_{\uw,y}([\nu_{\uw,y}])$ in
\begin{align*}
    \Hom(\Or(\vec{n}\times \Tree_1(z;w_i)),\Or(\vec{n}\times \Tree_1(z;w_i)))\otimes  \Hom(\Or(\Tree_k(\uw;y)),\Or(\wb^u(\uw)\times \overline{\mathcal{T}_k})).
\end{align*}

We can now compute:
\begin{align*}
    f_{\ux,y}([\nu_{\ux,y}])_p&\times [\vec{n}]\times [\Tree_k(\uw;y)]_p\times[\Tree_1(x_i;w_i)]_p\\
    &=(-1)^{\eta_1}f_{\ux,y}([\nu_{\ux,y}])_p\times [\vec{n}]\times[\Tree_1(x_i;w_i)]_p\times [\Tree_k(\uw;y)]_p\\
    &\mapsto(-1)^{\eta_1+\eta_2}\id\otimes f_{\uw,y}([\nu_{\uw,y}])([\vec{n}]\times[\Tree_1(x_i;w_i)]_p\otimes [\Tree_k(\uw;y)]_p)\\
    &\mapsto(-1)^{\eta_1+\eta_2+\eta_3}[\vec{n}]\times [\Tree_1(x_i;w_i)]_p\times f_{\uw,y}(\nu_{\uw,y})([\Tree_k(\uw;y)]_p)\\
    &=(-1)^{\eta_1+\eta_2+\eta_3+\eta(\uw;y)}[\vec{n}]\times [\Tree_1(x_i;w_i)]_p\times([\wb^u(\uw)]\times [\overline{\mathcal{T}_k}])_p\\
    &=(-1)^{\eta_1+\eta_2+\eta_3+\eta(\uw;y)+\eta_4}([\wb^u(\uw_{< i})]\times [\vec{n}]\times [\Tree_1(x_i;w_i)]\times[\wb^u(\uw_{\geq i})]\times [\overline{\mathcal{T}_k}])_p
\end{align*}
where $\eta(\uw;y)$ as in Lemma \ref{lem:explicit orient k>1} and
\begin{align*}
    \eta_1&\equiv |\Tree_k(\uw;y)||\Tree_1(x_i;w_i)|\equiv (k+\|\uw\|+\|y\|)(1+\|x_i\|+\|w_i\|),\\
    \eta_2&\equiv \|x_i\|+\|w_i\|,\\
    \eta_3&\equiv (\|x_i\|+\|w_i\|)|\nu_{\uw,y}|\equiv (\|x_i\|+\|w_i\|)(nk+\|y\|),\\
    \eta_4&\equiv (\|x_i\|+\|w_i\|)|\wb^u(\uw_{<i})|=(\|x_i\|+\|w_i\|)|\uw_{<i}|=(\|x_i\|+\|w_i\|)(\|\uw_{<i}\|+n(i-1)).
\end{align*}
We further note that by \eqref{eq: or L(x,y)}
\begin{align*}
     [\vec{n}]\times [\Tree_1(x_i;w_i)]\times[\wb^u(w_i)]&=(-1)^{|x_i|+|x_i||w_i|}[\vec{n}]\times [\overline{\mathcal{L}}(x_i,w_i)]\times[\wb^u(w_i)]\\
     &=(-1)^{|x_i|+|x_i||w_i|+|x_i|+|w_i|}[\wb^u(x_i)].
\end{align*}
We set 
\begin{align*}
    \eta_5\equiv |x_i|+|x_i||w_i|+|x_i|+|w_i|\equiv |x_i||w_i|+|w_i|\equiv (\|x_i\|+\|w_i\|)(n+\|w_i\|)
\end{align*}
and thus have
\begin{align*}
    [\nu_{\ux,y}]_p&\times [\vec{n}]\times [\Tree_k(\uw;y)]_p\times[\Tree_1(x_i;w_i)]_p=(-1)^{\eta_1+\eta_2+\eta_3+\eta_4+\eta_5+\eta(\uw;y)}([\wb^u(\ux)]\times [\overline{\mathcal{T}_k}])_p.
\end{align*}
Finally, we recall that 
\begin{align*}
    f_{\ux,y}([\nu_{\ux,y}])_p\times[\Tree_k(\ux;y)]_p=(-1)^{\eta(\ux;y)}([\wb^u(\ux)]\times [\overline{\mathcal{T}_k}])_p.
\end{align*}
We thus find
\begin{align*}
    [\Tree_k(\ux;y)]_p=(-1)^{\eta_1+\eta_2+\eta_3+\eta_4+\eta_5+\eta(\uw;y)+\eta(\ux;y)}[\vec{n}]\times [\Tree_k(\uw;y)]\times [\Tree_1(x_i;w_i)]_p.
\end{align*}
By Lemma \ref{lem:explicit orient k>1}, we find $\eta(\uw;y)+\eta(\ux;y)\equiv (\|x_i\|+\|w_i\|)(n(k-i)+\sum_{j\neq i} \|w_j\|)$. We conclude by computing:
\begin{align*}
    \eta_1+\eta_2&+\eta_3+\eta_4+\eta_5+\eta(\uw;y)+\eta(\ux;y)\\
    \equiv&\  k+ \|\uw\|+\|y\|+ (\|x_i\|+\|w_i\|)(k+\|\uw\|+\|y\|+1+nk+\|y\|+\|\uw_{<i}\|+n(i-1)+n+\|w_i\|)\\
    &\ +(\|x_i\|+\|w_i\|)(n(k-i)+\sum_{j\neq i} \|w_j\|)\\
    \equiv&\  k+ \|\uw\|+\|y\|+ (\|x_i\|+\|w_i\|)(k+1+\|\uw_{< i}\|)\\
    \equiv&\  k+ \|\ux\|+\|y\|+ (\|x_i\|+\|w_i\|)(k+\|\uw_{< i}\|)\\
    \equiv&\  \varepsilon_0.
\end{align*}

\subsubsection{The case \texorpdfstring{$k=m>1=l$}{k=m>1=l}}
For this case, we need to characterise $[\Tree_1(w;y)]$ in terms of the orientations $[\wb^s(w)]$ and $[\wb^s(y)]$:

\begin{lemme}\label{lem:T1 in terms of Ws}
    As local orientations at $w$, it holds:
    \begin{align*}
        [\wb^s(y)]=-[\vec{n}]\times [\Tree_1(w;y)]\times [\wb^s(w)].
    \end{align*}
\end{lemme}
\begin{proof}
    We fix $\lambda\in \mathcal{P}_{w\to y}X$ that is a parametrisation of a trajectory in $\mathcal{L}(w,y)$. We denote by $\phi_\lambda$ the parallel transport along $\lambda$. Parallel transport commutes with local orientations, that is, for a submanifold $A$ containing $\lambda$ that is invariant under $\phi_\lambda$, the following diagram commutes:
    \begin{center}
        \begin{tikzcd}[column sep=0.5 cm]
            \Or_w(A)\ar[rr, "\phi_\lambda"] \ar[rd, "\cong"] && \Or_y(A)\ar[ld, "\cong"']\\
            & \Or(A).
        \end{tikzcd}
    \end{center}
    The two cases we are interested in are $A=X$ and $A=\overline{\mathcal{L}}(x,y)$.
    
    We note that $\vec{m}\coloneqq\phi_\lambda(\vec{n})\in T_yX$ is a vector pointing out of $\wb^s(y)$ and thus into $\wb^u(y)$. See Figure \ref{fig:Ws_and_L}. 
    \begin{figure}[H]
        \centering
        \includegraphics[width=0.3\linewidth]{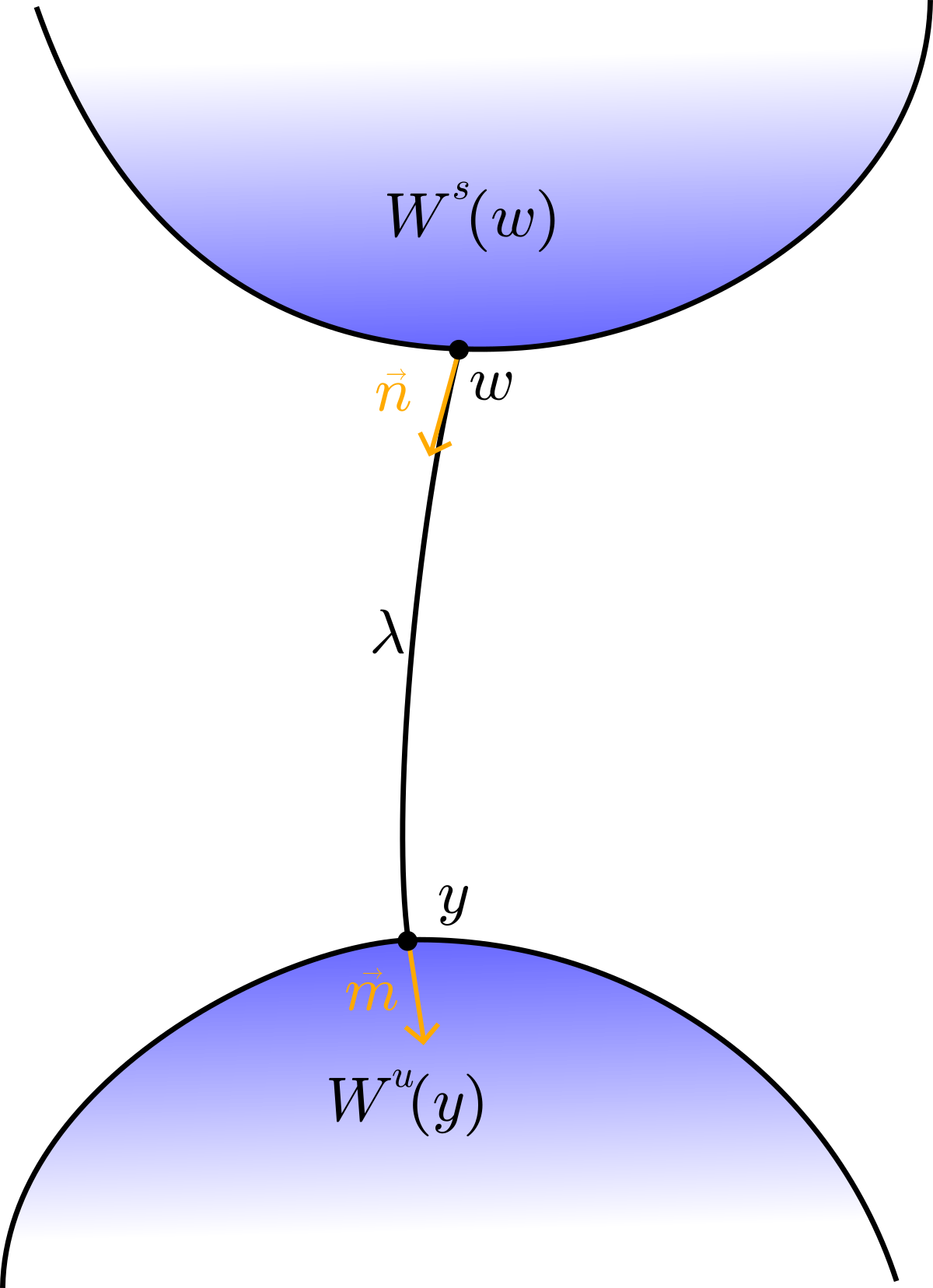}
        \caption{The trajectory $\lambda$ pushes the vector $\vec{n}\in T_w\wb^s(w)$ pointing out of $\wb^s(w)$ to $\vec{m}\in T_y\wb^u(y)$ pointing into $\wb^u(y)$.}
        \label{fig:Ws_and_L}
    \end{figure}
    
    We recall that by definition $[\Tree_1(w;y)]=(-1)^{\|w\|+\|w\|\|y\|+n\|w\|+n\|y\|}[\overline{\mathcal{L}}(w,y)]$ and that $[\vec{m}]\times [\overline{\mathcal{L}}(w,y)]\times [\wb^u(y)]=(-1)^{\|w\|+\|y\|+1}[\wb^u(w)]$ and thus as local orientations at $y$:
    \begin{align*}
        [\wb^u(w)]&=(-1)^{\|y\|+\|w\|\|y\|+n\|w\|+n\|y\|+1}[\vec{m}]\times [\Tree_1(w;y)]\times [\wb^u(y)]\\
        &=(-1)^{\|y\|+\|w\|\|y\|+n\|w\|+n\|y\|+1}[\phi_\lambda(\vec{n})]\times [\Tree_1(w;y)]\times \phi_\lambda([\wb^u(y)])\\
        &=(-1)^{\|y\|+\|w\|\|y\|+n\|w\|+n\|y\|+1}\phi_\lambda([\vec{n}]\times [\Tree_1(w;y)]\times [\wb^u(y)]).
    \end{align*}
    Applying $-\times [\wb^s(w)]$, we find by definition of $[\wb^s(w)]$ \eqref{eq:def of [ws]} that
    \begin{align*}
        [X]&=(-1)^{\|y\|+\|w\|\|y\|+n\|w\|+n\|y\|+1}\phi_\lambda([\vec{n}]\times [\Tree_1(w;y)]\times [\wb^u(y)])\times [\wb^s(w)]\\
        &=(-1)^{\|y\|+\|w\|\|y\|+n\|w\|+n\|y\|+1}\phi_\lambda([\vec{n}]\times [\Tree_1(w;y)]\times [\wb^u(y)]\times [\wb^s(w)]).
    \end{align*}
    Because $\phi_\lambda([X])=[X]$, we can deduce that, as local orientations at $w$, it holds:
    \begin{align*}
        [X]&=(-1)^{\|y\|+\|w\|\|y\|+n\|w\|+n\|y\|+1}[\vec{n}]\times [\Tree_1(w;y)]\times [\wb^u(y)]\times [\wb^s(w)]\\
        &=(-1)^{\|y\|+\|w\|\|y\|+n\|w\|+n\|y\|+1+|y|\|w\|}[\vec{n}]\times [\Tree_1(w;y)]\times [\wb^s(w)]\times [\wb^u(y)]\\
        &=(-1)^{\|y\|+n\|y\|+1}[\vec{n}]\times [\Tree_1(w;y)]\times [\wb^s(w)]\times [\wb^u(y)].
    \end{align*}
    By definition of $[\wb^s(y)]$, we find that 
    \begin{align*}
        [\wb^s(y)]\times [\wb^u(y)]=(-1)^{|y|\|y\|}[X]=(-1)^{\|y\|+n\|y\|}[X].
    \end{align*}
    Combining the last two equations we find $[\wb^s(y)]\times [\wb^u(y)]=-[\vec{n}]\times [\Tree_1(w;y)]\times [\wb^s(w)]\times [\wb^u(y)]$ which proves the lemma.
\end{proof}

We next note that both $f_{\ux,y}$ and $f_{\ux,w}$ are induced by the map $\Phi_{\ux}$. Therefore the following diagram commutes:
\begin{center}
    \begin{tikzcd}
        \Hom(\Or(\wb^s(y)),\Or(X^k))\ar[r, "f_{\ux,y}"] \ar[d] &\Hom(\Or(\Tree_k(\ux;y)),\Or(\wb^u(\ux)\times \overline{\mathcal{T}_k}))\ar[d]\\
        \Hom(\Or(\vec{n}\otimes\Or(\Tree_1(w;y))\otimes \Or(\wb^s(w)),\Or(X^k))\ar[d] &\makecell{\Hom(\Or(\vec{n})\otimes \Or(\Tree_1(w;y))\otimes \Or(\Tree_k(\ux;w)),\\\Or(\wb^u(\ux)\times \overline{\mathcal{T}_k}))}\ar[d]\\
        \makecell{\Hom(\Or(\vec{n})\otimes \Or(\Tree_1(w;y)), \mathbbm{1})
        \\\otimes\Hom(\Or(\wb^s(w)),\Or(X^k))}\ar[r, "\id\otimes f_{\ux,w}"] & \makecell{\Hom(\Or(\vec{n})\otimes \Or(\Tree_1(w;y)), \mathbbm{1})\\
        \otimes\Hom( \Or(\Tree_k(\ux;w)),\Or(\wb^u(\ux)\times \overline{\mathcal{T}_k}))}
    \end{tikzcd}
\end{center}

By Lemma \ref{lem:T1 in terms of Ws}, we thus have
\begin{align*}
    [X^k]&=[\nu_{\ux,y}]([\wb^s(y)])\\
    &=-[\nu_{\ux,y}]([\vec{n}]\times [\Tree_1(w;y)]\times [\wb^s(w)])\\
    &\mapsto (-1)^{1+|\nu_{\ux,w}||\vec{n}\times \Tree_1(w;y)|} ([\vec{n}]\times [\Tree_1(w;y)])^\lor\otimes [\nu_{\ux,w}]([\vec{n}]\times [\Tree_1(w;y)]\otimes [\wb^s(w)]).
\end{align*}
This shows that the left hand vertical map in the above diagram sends $[\nu_{\ux,y}]$ to 
\begin{align*}
    (-1)^{1+|\nu_{\ux,w}||\vec{n}\times \Tree_1(w;y)|} ([\vec{n}]\times [\Tree_1(w;y)])^\lor\otimes [\nu_{\ux,w}].
\end{align*}
By commutativity of the diagram, the right hand vertical maps send $f_{\ux,y}([\nu_{\ux,y}])$ to 
\begin{align*}
    (-1)^{1+|\nu_{\ux,w}||\vec{n}\times \Tree_1(w;y)|} ([\vec{n}]\times [\Tree_1(w;y)])^\lor\otimes f_{\ux,w}([\nu_{\ux,w}]).
\end{align*}

By Lemma \ref{lem:explicit orient k>1}, we deduce:
\begin{align*}
    f_{\ux,y}([\nu_{\ux,y}])_p\times [\vec{n}]\times [\Tree_1(w;y)]_p\times [\Tree_k(\ux;w)]_p\mapsto&\  (-1)^{1+|\nu_{\ux,w}||\vec{n}\times \Tree_1(w;y)|} ([\vec{n}]\times [\Tree_1(w;y)])^\lor\otimes f_{\ux,w}([\nu_{\ux,w}])\\
    &\ ( [\vec{n}]\times [\Tree_1(w;y)]_p\otimes [\Tree_k(\ux;w)]_p)\\
    =& - ([\vec{n}]\times [\Tree_1(w;y)])^\lor([\vec{n}]\times [\Tree_1(w;y)])\\
    &\ \otimes f_{\ux,w}([\nu_{\ux,w}])\times [\Tree_k(\ux;w)]\\
    \mapsto&\ (-1)^{1+\eta(\ux;w)}([\wb^u(\ux)]\times [\overline{\mathcal{T}_k}])_p.
\end{align*}
Using that 
\begin{align*}
    f_{\ux,y}([\nu_{\ux,y}])_p\times [\Tree_k(\ux;y)]_p=(-1)^{\eta(\ux;y)}([\wb^u(\ux)]\times [\overline{\mathcal{T}_k}])_p,
\end{align*}
we find
\begin{align*}
    [\vec{n}]\times [\Tree_1(w;y)]_p\times [\Tree_k(\ux;w)]_p=(-1)^{1+\eta(\ux;w)+\eta(\ux;y)}[\Tree_k(\ux;y)]_p.
\end{align*}
Finally by Lemma \ref{lem:explicit orient k>1}, we find
\begin{align*}
    1+\eta(\ux;w)+\eta(\ux;y)\equiv 1+\|w\|+\|y\|\equiv \eta_0,
\end{align*}
which proves this case.

\subsubsection{The case \texorpdfstring{$l,m>1$}{l,m>1}}
We finally assume that $l,m>1$ and thus no $\mathcal{L}(x,y)$ are involved. In this case, we need to relate $\alpha_{\ux,y}$ to $\alpha_{\uw,y}\otimes \alpha_{\uz,w_i}$. The key here is to relate 
\begin{align*}
    f_{\ux,y}\colon \Hom(\Or(\wb^s(y)),\Or(X^k))\to \Hom(\Or(\Tree_k(\ux;y)),\Or(\wb^u(\ux)\times \overline{\mathcal{T}_k}))
\end{align*}
to 
\begin{align*}
    f_{\uw,y}\otimes f_{\uz,w_i}\colon &\Hom(\Or(\wb^s(y)),\Or(X^l))\otimes \Hom(\Or(\wb^s(w_i)),\Or(X^m))\\
    &\to \Hom(\Or(\Tree_l(\uw;y)),\Or(\wb^u(\uw)\times \overline{\mathcal{T}_l}))\otimes \Hom(\Or(\Tree_m(\uz;w_i)),\Or(\wb^u(\uz)\times \overline{\mathcal{T}_m}))
\end{align*}
from Lemma \ref{lem:transport orient}. In fact, the rest of the definition of $\alpha_{\ux,y}$ is given by shuffling factors and cancelling them and these operations commute as explained in Remark \ref{rem:cancelling orient}.

The geometrical data underlying the isomorphism $\alpha_{\ux,y}$ is the map $\Phi_{\ux}$. In order to prove this case, we therefore need to relate $\Phi_{\ux}$ to $\Phi_{\uw}$ and $\Phi_{\uz}$. The problem is that $\Phi_{\ux}$ is not defined on the $\wb^u(\ux)\times \partial \mathcal{T}_k\subseteq \partial(\wb^u(\ux)\times \overline{\mathcal{T}_k})$. We circumvent this problem by extending $\Phi_{\uw}\times \Phi_{\uz}$ to a suitable neighbourhood of $\mathcal{T}_l(\uw;y)\times \mathcal{T}_m(\uz;w_i)\subseteq W^u(\uw)\times W^u(\uz)\times \overline{\mathcal{T}_k}$. 

We denote by $\widetilde{\mathcal{T}_k}\subseteq \overline{\mathcal{T}_k}$ the subspace of all trees that have an edge $e_0$ of length $>2$ and such that if the length of $e_0$ is changed to $\infty$, the resulting tree is in $\circ_i(\mathcal{T}_l\times \mathcal{T}_m)\subseteq \partial  \overline{\mathcal{T}_k}$. In particular, this is an open neighbourhood of $\mathcal{T}_l\times \mathcal{T}_m$ homeomorphic to $\mathcal{T}_l\times \mathcal{T}_m\times (2,\infty]$ and $\mathcal{T}_l\times \mathcal{T}_m$ is embedded as $\mathcal{T}_l\times \mathcal{T}_m\times \{\infty\}$.

We define $\widetilde{\mathcal{T}_k}(\ux;y)\subseteq \Tree_k(\ux;y)$ the subspace of perturbed Morse gradient trees that have underlying metric tree in $\widetilde{\mathcal{T}_k}$.

We now consider the diagram
\begin{equation}\label{diag:k,l>1}
    \begin{tikzcd}
        \mathcal{T}_l^\X(\uw;y)\times \mathcal{T}_m^\X(\uz;w_{i})\ar[r, hook, "\iota_1"] \ar[d, hook, "\iota_2"] &  W^u(\uw)\times W^u(\uz)\times \mathcal{T}_l\times \mathcal{T}_m\ar[d, hook, "\iota_3"]\\
        W^u(w_i)\times \widetilde{\mathcal{T}_k}^\X(\ux;y) \ar[r, hook, "\iota_4"] &  W^u(w_{i})\times W^u(\ux)\times\widetilde{\mathcal{T}_k}.
    \end{tikzcd}
\end{equation}

\begin{figure}[!ht]
    \centering
    \includegraphics[width=0.5\linewidth]{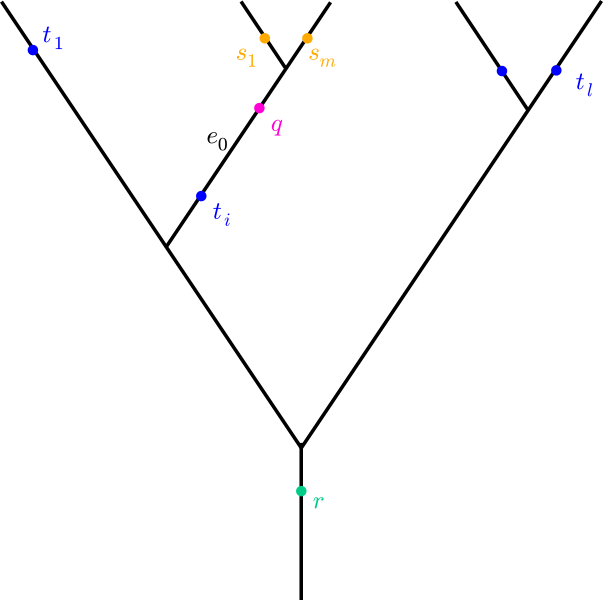}
    \caption{A tree in $\widetilde{\mathcal{T}_k}$ where $e_0$ has length $>2$. The times $t_i,s_i,r$ and $q$ are at distance one to the nearest vertex.}
    \label{fig:eval_at_broken_tree}
\end{figure}

The map $\iota_1$ evaluates both trees on all input legs at the point that has distance one to the first vertex. This corresponds to evaluating the tree at $t_1,\dots,t_l$ and $s_1,\dots, s_m$ in Figure \ref{fig:eval_at_broken_tree}.

The map $\iota_2$ is the inclusion of $\mathcal{T}_l(\uw;y)\times \mathcal{T}_m(\uz;w_{i})$ into $\widetilde{\mathcal{T}_k}^\X(\ux;y)$ and additionally evaluating the tree at $t_{i}$ as in Figure \ref{fig:eval_at_broken_tree}.

The map $\iota_3$ is given by the inclusion of $\mathcal{T}_l\times \mathcal{T}_m$ into $\widetilde{\mathcal{T}_k}$ and reordering the points $W^u(\uw)\times W^u(\uz)\cong W^u(w_i)\times W^u(\ux)$. 

The map $\iota_4$ evaluates the tree in $\widetilde{\mathcal{T}_k}^\X(\ux;y)$ at the points $t_1,\dots,t_{i-1},s_1,\dots,s_m,t_{i+1},\dots ,t_l$ as in Figure \ref{fig:eval_at_broken_tree}.

We recall from Definition \ref{defi : Phi_X} that for $s,t\in \gamma$, we define $\phi_{T,\X}^{s,t}\colon X\to X$ as the flow that pushes a point according to the perturbation of the path from $s$ to $t$ in $\gamma$. This is defined as long as the path from $s$ to $t$ does not pass through $e_0$ or $\gamma\notin \partial\widetilde{\mathcal{T}_k}$.

We define the map 
\begin{align*}
    \widetilde{\Phi}\colon  W^u(w_{i})\times W^u(\ux)\times  \widetilde{\mathcal{T}_k}\to&\  X^{k+1},\\
    (\eta_i,\eta_1,\dots,\eta_{i-1},\theta_1,\dots,\theta_m,\eta_{i+1},\dots,\eta_l,\gamma)\mapsto&\  (\phi_{\gamma,t_1\to r}(\eta_1),\dots,\phi_{\gamma,t_l\to r}(\eta_l),\\
    &\ \phi_{\gamma,s_1\to q}(\theta_1),\dots,\phi_{\gamma,s_m\to q}(\theta_m))
\end{align*}
for $t_j,s_j,r,q$ the points in $\gamma$ as in Figure \ref{fig:eval_at_broken_tree}. 

We note that $\widetilde{\Phi}\circ\iota_3=\Phi_{\uw}\times \Phi_{\uz}$. This shows that $\widetilde{\Phi}$ intersects $\Delta^l(W^s(y))\times \Delta^m(W^s(w_i))$ transversely and 
\begin{align*}
    \iota_1\colon \mathcal{T}^{\X}_l(\uw;y)\times \mathcal{T}_m^\X(\uz;w_i)\cong \widetilde{\Phi}^{-1}(\Delta^l(W^s(y))\times \Delta^m(W^s(w_i))).
\end{align*}

On the other hand, for $\gamma\in \mathcal{T}_k\cap \widetilde{\mathcal{T}_k}$, we have that $\phi_{\gamma,q\to r}$ is defined. We fix $(\underline{\eta},\underline{\theta})\in W^u(\uw)\times W^u(\uz)$ and compute that
\begin{align*}
    \sigma (\phi_{\gamma,t_i\to r}(\eta_i)&,\Phi_{\ux}(\eta_1,\dots,\eta_{i-1},\theta_1,\dots,\theta_m,\eta_{i+1},\dots,\eta_l))\\
    &=(\id^{\times l}\times \phi_{\gamma,q\to r}^{\times m})\circ \widetilde{\Phi}(\eta_i,\eta_1,\dots,\eta_{i-1},\theta_1,\dots,\theta_m,\eta_{i+1},\dots,\eta_l,\gamma)
\end{align*}
where $\sigma$ is a shuffle of coordinates according to the map $(w_i,\ux)\to (\uw;\uz)$. This shows that $\widetilde{\Phi}\circ \iota_4$ is smoothly isotopic to $\sigma\circ(\id\times \Phi_{\ux})$. We also note that 
\begin{align*}
    \sigma\circ (\id\times \Delta^k),\Delta^l\times \Delta^m\colon W^u(w_i)\times W^s(y)\to X^{l+m}
\end{align*}
are smoothly isotopic. 

Moreover, $\sigma\circ(\id\times \Phi_{\ux})$ intersects $\sigma\circ(\id\times \Delta^k)$ transversely and
\begin{align*}
    \iota_4\colon W^u(w_i)\times \mathcal{T}^{\X}_k(\ux;y)\cong W^u(w_i)\times \Phi_{\ux}^{-1}(\Delta^k(W^s(y)))= (\sigma\circ(\id\times \Phi_{\ux}))^{-1}\sigma\circ(\id\times \Delta^k)(W^u(w_i)\times W^s(y)). 
\end{align*}

This shows that $\widetilde{\Phi}$ and $\Phi_{\ux}$ induce the same map
\begin{align*}
    \Hom(\Or(\wb^u(w_i))\otimes \Or(\wb^s(y)),\Or(X^{l+m}))&\to\\ 
    \Hom(\Or(\wb^u(w_i))\otimes \Or(\Tree_k(\ux;y)),&\ \Or(\wb^u(w_i))\otimes \Or(\wb^u(\ux))\otimes \Or(\overline{\mathcal{T}_k})).
\end{align*}
The map that $\Phi_{\ux}$ induces is the same as given by evaluating at $[\wb^u(w_i)]$ and applying $f_{\ux,y}$.

The upshot of this is that the following diagram commutes:
\begin{center}
    \begin{tikzcd}
        \makecell{\Hom(\Or(\wb^s(y)),\Or(X^{k}))\\\otimes \Hom(\Or(\wb^s(w_i)),\Or(X))}\ar[r, "\sigma"] \ar[d] & \Hom(\Or(\wb^s(y)\times \wb^s(w_i)),\Or(X^{l+m})\ar[d]\\
        \Hom(\Or(\wb^s(y)),\Or(X^{k}))\otimes \Or(\wb^u(w_i))\ar[d, "f_{\ux,y}"] 
        &\makecell{\Hom(\Or(\wb^s(y)),\Or(X^l))\\\otimes \Hom(\Or(\wb^s(w_i)),\Or(X^m)) } \ar[d, "f_{\uw,y}\otimes f_{\uz,w_i}"]\\
        \makecell{\Hom ( \Tree_k(\ux;y)), \wb^u(\ux)\times \overline{\mathcal{T}_k}))\otimes \Or(\wb^u(w_i))}\ar[d]& \makecell{\Hom(\Or(\Tree_l(\uw;y)),\Or(\wb^u(\uw)\times \overline{\mathcal{T}_l}))\\
        \otimes \Hom(\Or(\Tree_m(\uz;w_i)),\Or(\wb^u(\uz)\times \overline{\mathcal{T}_m}))}\ar[d]\\
         \Hom (\Or(\Tree_k(\ux;y)),\Or(\wb^u(w_i)\times \wb^u(\ux)\times \overline{\mathcal{T}_k}))& \makecell{\Hom(\Or(\vec{n}\times \Tree_l(\uw;y)\times \Tree_m(\uz;w_i)),\\
        \Or(\vec{n}\times \wb^u(\uw)\times\wb^u(\uz)\times \overline{\mathcal{T}_l}\times \overline{\mathcal{T}_m})).}\ar[l, "\circ_i"]
    \end{tikzcd}
\end{center}
The top left vertical map is given by
\begin{align*}
    \Hom(\Or(\wb^s(w_i)),\Or(X))&\to \Or(\wb^u(w_i)),\\
    ([\wb^s(w_i)]\mapsto[\wb^u(w_i)]\times [\wb^s(w_i)]=[X])&\mapsto[\wb^u(w_i)] .
\end{align*}

This shows that $f_{\ux,y}$ is compatible with $f_{\uw,y}\otimes f_{\uz,w_i}$. Because $\Or$ commutes with shuffling and evaluation, this shows that $\alpha_{\ux,y}$ commutes with $\vec{n}\times -$, that is, the following diagram commutes:
\begin{center}
    \begin{tikzcd}[column sep=2cm ]
        \Or(\vec{n})\otimes \Or(\Tree_l(\uw;y))\otimes \Or(\Tree_m(\uz;w_i))\ar[d, "\circ_i"] \ar[r, "\id\otimes \alpha_{\uw,y}\otimes \alpha_{\uz,w_i}"] & \makecell{\Or(\vec{n})\otimes \underline{\Sigma}^{\|\cdot \|}(\uw;y)\otimes \Or(\overline{\mathcal{T}_l})\\
        \underline{\Sigma}^{\|\cdot \|}(\uz;w_i)\otimes \Or(\overline{\mathcal{T}_m})}\ar[d, "\circ_i"]\\
        \Or(\Tree_k(\ux;y))\ar[r, "\alpha_{\ux,y}"]& \underline{\Sigma}^{\|\cdot \|}(\ux;y)\otimes \Or(\overline{\mathcal{T}_k}).
    \end{tikzcd}
\end{center}
Because $[\Tree_k(\ux;y)]$ is defined by $\alpha_{\ux,y}([\Tree_k(\ux;y)]_p)=(\tau_{\ux,y}\otimes [\overline{\mathcal{T}_k}])_p$ this shows that 
\begin{align*}
    [\vec{n}]\times [\Tree_l(\uw;y)]_p\times [\Tree_m(\uz;w_i)]_p=\varepsilon_0 [\Tree_k(\ux;y)]_p
\end{align*}
for all $p\in \Tree_k(\ux;y)$. This concludes the proof of Lemma \ref{lem:existence of orientations appendix}.

\setcounter{secnumdepth}{-1}

\normalem
\bibliography{mybib}
\bibliographystyle{alpha}

\end{document}